\documentclass[3p,sort & compress]{elsarticle}
\usepackage{amssymb}
\usepackage{latexsym}
\usepackage{amsfonts}
\usepackage{amsthm}
\usepackage{amsbsy}
\usepackage{amsmath}
\usepackage{graphics,graphicx}
\usepackage{subfigure}
\usepackage{multirow,multicol}
\usepackage{bm}
\usepackage{bbm}
\usepackage{color}
\usepackage{float}
\usepackage{calc}
\usepackage{caption}
\usepackage{dsfont}
\usepackage{ifthen}
\usepackage{mathrsfs}
\usepackage[colorlinks,linkcolor=blue,citecolor=blue]{hyperref}
\usepackage{epstopdf}
\usepackage{epsfig}
\usepackage{algorithm}
\usepackage{algorithmic}
\usepackage{pifont}
\usepackage{natbib}
\usepackage{overpic}
\usepackage{psfrag}
\usepackage{rotating}
\usepackage{stmaryrd}
\usepackage{verbatim}
\usepackage{setspace}
\usepackage{tikz}
\usepackage{listings} 
\usepackage{xcolor} 
\usepackage{geometry}
\usepackage{mlmath}

\newdefinition{remark}{Remark}

\newdefinition{method}{Method}
\newdefinition{example}[theorem]{Example}

\numberwithin{equation}{section}

\numberwithin{theorem}{section}

\newcommand{\orcid}[1]{\href{https://orcid.org/#1}{\includegraphics[width=8pt]{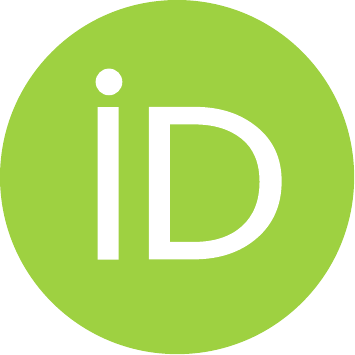}}}
\journal{Journal of \LaTeX\ Templates}
\date{June, 20, 2023}
\begin{document}
\begin{frontmatter}
    \title{Solving a class of multi-scale elliptic PDEs by Fourier-based mixed physics informed neural networks}
    \author[Ceyear]{Xi'an Li\orcid{0000-0002-1509-9328}\corref{cor1}}\ead{lixian9131@163.com}
    \author[acu]{Jinran Wu\orcid{0000-0002-2388-3614}} \ead{ryan.wu@acu.edu.au}
    \author[acu]{You-Gan Wang\orcid{0000-0003-0901-4671}}\ead{you-gan.wang@acu.edu.au}
    \author[Ceyear]{Xin Tai}\ead{taixin@ceyear.com}
    \author[41st]{Jianhua Xu}\ead{xujianhua@ei41.com}
    \cortext[cor1]{Corresponding author.}
    \address[Ceyear]{Ceyear Technologies Co., Ltd, Qingdao 266000, China}
    \address[acu]{The Institute for Learning Sciences and Teacher Education, Australian Catholic University, Brisbane 4000, Australia}
    \address[41st]{The 41st Institute of China Electronics Technology Group Corporation, Qingdao, 266000, China}

    \begin{abstract}
    Deep neural networks have garnered widespread attention due to their simplicity and flexibility in the fields of engineering and scientific calculation. In this study, we probe into solving a class of elliptic partial differential equations(PDEs) with multiple scales by utilizing Fourier-based mixed physics informed neural networks(dubbed FMPINN), its solver is configured as a multi-scale deep neural network. In contrast to the classical PINN method, a dual (flux) variable about the rough coefficient of PDEs is introduced to avoid the ill-condition of neural tangent kernel matrix caused by the oscillating coefficient of multi-scale PDEs. Therefore, apart from the physical conservation laws, the discrepancy between the auxiliary variables and the gradients of multi-scale coefficients is incorporated into the cost function, then obtaining a satisfactory solution of PDEs by minimizing the defined loss through some optimization methods. Additionally, a trigonometric activation function is introduced for FMPINN, which is suited for representing the derivatives of complex target functions. Handling the input data by Fourier feature mapping will effectively improve the capacity of deep neural networks to solve high-frequency problems.  Finally, to validate the efficiency and robustness of the proposed FMPINN algorithm, we present several numerical examples of multi-scale problems in various dimensional Euclidean spaces. These examples cover both low-frequency and high-frequency oscillation cases, demonstrating the effectiveness of our approach. All code and data accompanying this manuscript will be made publicly available at \href{https://github.com/Blue-Giant/FMPINN}{https://github.com/Blue-Giant/FMPINN}.
    \end{abstract}
\end{frontmatter}

\section{Introduction}\label{sec:01}
Multi-scale problems, governed by partial differential equations(PDEs) with multiple scales, are prevalent in diverse scientific and engineering fields like reservoir simulation, high-frequency scattering, and turbulence modeling. This paper focuses on solving the following type of multi-scale problem.
\begin{equation}\label{eq:multiscale}
\begin{cases}
-\textbf{div}\bigg{(}A^{\varepsilon}(\bm{x})\nabla u^{\varepsilon}(\bm{x})\bigg{)} = f(\bm{x}), ~~~\bm{x}\in \Omega,\\
~~~~~~~~~~~\mathcal{B}u^{\varepsilon}(\bm{x}) = g(\bm{x}),~~~~~~~~~~~~\bm{x}\in \partial\Omega.
\end{cases}
\end{equation}
where $\Omega$ is a bounded subset of $\mathbb{R}^{d}(d=1,2,3,\ldots)$ with piecewise Lipschitz boundary and satisfies the interior cone condition, $\varepsilon$ is a small positive parameter that signifies explicitly the multiscale nature of the rough coefficient $A^{\varepsilon}(\bm{x})$. 
$\mathcal{B}$ is a boundary operator in $\partial\Omega$ that imposes the boundary condition of $u^{\varepsilon}$, such as Dirchlete, Neumman and Robin. $\nabla $ and $\textbf{div}$ are the gradient and divergence operators, respectively. $f(\bm{x})\in L^2(\Omega)$ is a given function. In addition, $A^{\varepsilon}(\bm{x})$ is symmetric and uniformly elliptic on $\Omega$. It means that all eigenvalues of $A^{\varepsilon}$
are uniformly bounded by two strictly positive constants $\lambda_{\min}(A^{\varepsilon})$ and $\lambda_{\max}(A^{\varepsilon})$. In other word, for all $\bm{x}\in\Omega$ and $\bm{\xi}\in\mathbb{R}^d$, we have
\begin{equation}
     \lambda_{\min}(A^{\varepsilon})|\bm{\xi}|^2\leqslant \bm{\xi}^T A^{\varepsilon}(\bm{x})\bm{\xi}\leqslant \lambda_{\max}(A^{\varepsilon})|\bm{\xi}|^2.
\end{equation}

The multi-scale problem\eqref{eq:multiscale} frequently arise in the fields of physical simulations and engineering applications, including the study of flow in porous media and the analysis of mechanical properties in composite materials \cite{Ming2006,ming2005analysis,li2012efficient}. Generally, the analytical solutions of \eqref{eq:multiscale} are seldom available, then solving numerically this problem through approximation methods is necessary. Lots of numerical methods focus on efficient,
accurate and stable numerical schemes have gained favorable achievement, such as heterogeneous multi-scale methods \cite{ming2005analysis,li2012efficient,abdulle2014analysis}, numerical homogenization \cite{dur91,ab05,hellman2019numerical}, variational multi-scale methods \cite{hughes98,larson2007adaptive}, multi-scale finite element methods \cite{Arbogast_two_scale_04, eh09,ch03},
flux norm homogenization \cite{berlyand2010flux,owhadi2008homogenization}, rough polyharmonic splines (RPS)\cite{owhadi2014polyharmonic}, generalized multi-scale finite element methods \cite{Efendiev2013,chung2014adaptiveDG,chung2015residual}, localized orthogonal decomposition \cite{MalPet:2014,Henning2014}, etc.
In contrast to standard numerical methods including FEM and FDM, they alleviate substantially the computational complexity in handling all relevant scales, improve the numerical stabilities and expedite the convergence. However, they still will encounter the curse of complex domain and dimensionality in general.

Deep neural networks(DNN), an efficient meshfree method without the discretization for a given interested domain, have drawn more and more attention from researchers to solve numerically the ordinary and partial differential equations as well as the inverse problems for complex geometrical domain and high-dimensional cases \cite{e2018the,sirignano2018dgm,chen2018neural,raissi2019physics,khoo2019switchnet,zang2020weak,lyu2022mim}, due to their extraordinary universal approximation capacity\cite{hauptmann2020deep}. Among these methods, the physics-informed neural networks (PINN) dating  back to the early 1990s again attracted widespread attention of researchers  and have made remarkable achievements for approximating the solution of PDEs by embracing the physical laws with neural networks, on account of the rapid development of computer science and technology\cite{raissi2019physics,dissanayake1994neural}. This method skillfully  incorporates the residual of governing equations and the discrepancy of boundary/initial constraints, then formulates a cost function can be optimized easily via the automatic differentiation in DNN. Many efforts have been made to further enhance the performance of PINN are concluded as two aspects:  refining the selection of the residual term and designing the manner of initial/boundary constraints. In terms of the residual term, there are XPINN~\cite{jagtap2021extended}, cPINN~\cite{jagtap2020conservative}, two-stage PINN~\cite{lin2022two} and gPINN~\cite{yu2022gradient}, and so on. By subtly encoding the I/B constraints into DNN in a hard manner, the PINN can be easy to train with low computational complexity and obtain a high-precision solution of PDEs with complex boundary conditions\cite{berg2018a,sun2020surrogate,lu2021physics}.
Motivated by the reduction of order in conventional methods\cite{ch03}, some attempts have been made to solve the high-order PDEs by reframing them as some first-order systems, this will overcome the shortcomings of the computational burden for high-order derivatives in DNN. For example, the deep mixed residual method \cite{lyu2022mim}, the local deep learning method\cite{zhu2021local} and the deep FOSLS method\cite{cai2020deep, bersetche2023deep}.

Many studies and experiments have indicated that the general DNN-based algorithms are commonly used to solve a low-frequency problem in varying dimensional space, but will encounter tremendous challenge for high-frequency problems such as multi-scale PDEs\eqref{eq:multiscale}. The frequency principle\cite{Xu_2020} or spectral bias\cite{rahaman2018spectral} of DNN shows that neural networks are typically efficient for fitting objective functions with low-frequency modes but inefficient for high-frequency functions. Then, a series of multi-scale DNN(MscaleDNN) algorithms were proposed to overcome the shortcomings of normal DNN for high-frequency problems by converting high-frequency contents into low-frequency ones via a radial scale technique \cite{liu2020multi,wang2020multiscale,li2020elliptic,li2023subspace}. After that, some corresponding mechanisms were developed to explain this performance of DNN, such as the Neural Tangent Kernel (NTK)\cite{wang2020eigenvector,jacot2018neural}. Furthermore, many researchers attempted to utilize a Fourier feature mapping consisting of sine and cosine to improve the capacity of MscaleDNN, which will alleviate the pathology of spectral bias and let neural networks capture high frequencies component effectively\cite{wang2020eigenvector,ramabathiran2021spinn,li2023subspace,li2023deep,tancik2020fourier,han2021hierarchical}. 

Recently, some works\cite{leung2022nh,carney2022physics} have shown that general PINN architecture is unable to capture the multi-scale property of the solution due to the effect of rough coefficient in multi-scale PDEs. In \cite{leung2022nh}, Wing Tat Leung et.al proposed a Neural homogenization-based PINN(NH-PINN) method to solve \eqref{eq:multiscale}, it can well overcome the unconvergence of PINN for multi-scale problems. However, NH-PINN also will encounter the dilemma of dimensional and the burden of computation, because it will convert one low-dimensional problem into a high-dimensional case. By carefully analyzing the Neural Tangent Kernel matrix associated with the PINN, Sean P. Carney et. al\cite{carney2022physics} found that the Forbenius norm of the NTK matrix will become unbound as the oscillation factor $\varepsilon$ in $A^{\varepsilon}$ tends to zero. It means that the evolution of residual loss term in PINN will become increasingly stiff as $\varepsilon\rightarrow 0$, then lead to poor training behavior for PINN.

In this paper, a Fourier-based multi-scale mixed PINN(FMPINN) structure is proposed to solve the multi-scale problems \eqref{eq:multiscale} with rough coefficients. This method consists of the general PINN architecture and the aforementioned MscaleDNN model with subnetworks being used to capture different frequencies component. To overcome the weakness of the normal PINN that failed to capture the jumping gradient information of the oscillating coefficient when tackling the governed equation in multi-scale PDEs\eqref{eq:multiscale}, a (dual)flux variable is introduced to alleviate the adverse effect of the rough coefficient. Meantime, it can also reduce the computational burden of PINN for the second-order derivatives of space variables. In addition, the Fourier feature mapping is used in our model to learn each target frequency efficiently and express the derivatives of multi-frequency functions easily, it will remarkably improve the capacity for our FMPINN model to solve multi-scale problems. In a nutshell, the primary contributions of this paper are summarized as follows:
\begin{enumerate}
    \item {We propose a novel neural networks approach by combining normal PINN and MscaleDNN with subnetworks structure to address multi-scale problems, leveraging the Fourier theorem and the F-principle of DNN.}
    \item {Inspired by the reduced order scheme for high-order PDEs, a dual (flux) variable about the rough coefficient of multi-scale PDEs is introduced to address the gradient leakage about the rough coefficient for PINN.}
    \item {By introducing some numerical experiments, we show that the classical PINN method with MscaleDNN solver is still insufficient in providing accurate solutions for multi-scale equations.}
    \item {We showcase the exceptional performance of FMPINN in solving a class of multi-scale elliptic PDEs with essential boundaries in various dimensional spaces. Our method outperforms existing approaches and demonstrates its superiority in addressing these complex problems.}
\end{enumerate}

The remaining parts of our work are organized as follows. In Section \ref{sec:02}, we briefly introduce the underlying conceptions and formulations for MscaleDNN and the structure of PINN. Section \ref{sec:03} provides a unified architecture of the FMPINN to solve the elliptic multi-scale problem \eqref{eq:multiscale} based on its equivalent reduced order scheme, and gives the option of activation function as well as the error analysis of our proppsed method. Section~\ref{algori2FMPINN} details the FMPINN algorithm for approximating the solution of multi-scale PDEs, then provide the option of activation function and the simple error analysis for FMPINN method. In Section \ref{sec:04}, some scenarios of multi-scale PDEs  are performed to evaluate the feasibility and effectiveness of our proposed method. Finally, some conclusions of this paper are made in Section \ref{sec:05}.

\section{Multi-scale Physics Informed Neural Networks}\label{sec:02}
\subsection{Multi-scale Deep Neural Networks with ResNet technique}
The basic concept and formulation of DNN are described briefly in this section, which helps audiences to understand the DNN structure through functional terminology. Mathematically, a deep neural network defines the following mapping
\begin{equation}
   \mathcal{F}: \bm{x}\in\mathbb{R}^{d}\Longrightarrow \bm{y}=\mathcal{F}(x)\in\mathbb{R}^{c}
\end{equation}
with $d$ and $c$ being the dimensions of input and output, respectively. In fact, the DNN functional $\mathcal{F}$ is a nested composition of the following single-layer neural unit:
\begin{equation}
    \bm{y}=\{y_1, y_2, \cdots, y_m\}~~\textup{and}~~  y_l= \sigma\left(\sum_{n=1}^{d} w_{ln}*x_n  + b_l\right)
\end{equation}
where $w_{ln}$ and $b_l$ are called weight and bias of $l_{th}$ neuron, respectively. $\sigma(\cdot)$ is an element-wise non-linear operator, generally referred as the activation function. Then, we have the following formulation of DNN:
\begin{equation}
		\bm{y}^{[\ell]} = \sigma\circ(\bm{W}^{[\ell]}\bm{y}^{[\ell-1]}+\bm{b}^{[\ell]}), ~~\text{for}~~\ell =1, 2, 3, \cdots\cdots, L
\end{equation}
and $\bm{y}^{[0]} = \bm{x}$, where $\bm{W}^{[\ell]} \in  \mathbb{R}^{n_{\ell+1}\times n_{\ell}}, \bm{b}^{[\ell]}\in\mathbb{R}^{n_{\ell+1}}$ stand for the weight matrix and bias vector of $\ell$-th hidden layer, respectively, $n_0=d$ and $n_{L+1}$ is the dimension of output, and $``\circ"$ stands for the elementary-wise operation. For convenience, the output of DNN is denoted by $\bm{y}(\bm{x};\bm{\theta})$ with $\bm{\theta}$ standing for its all weights and biases.

Residual neural network (ResNet) \cite{he2016deep} as a common skillful technique by introducing skip connections between adjacent or nonadjacent hidden layers can overcome effectively the vanishing gradient of parameters in the backpropagation for DNN, then make the network much easier to train and improve well the performance of DNN. Many experiment results showed that the ResNet can also improve the performance of DNNs to approximate high-order derivatives and solutions of PDEs \cite{e2018the,lyu2022mim}. 
We utilize the one-step skip connection scheme of ResNet in this work.  Except for the normal data flow, the data will also flow along with the skip connection if the two consecutive layers in DNN have the same number of neurons, otherwise, the data flows directly from one to the next layer. The filtered $\bm{y}^{[\ell+1]}(\bm{x};\bm{\theta})$ produced by the input $\bm{y}^{[\ell]}(\bm{x};\bm{\theta})$ is expressed as 
\begin{equation*}
\bm{y}^{[\ell+1]}(\bm{x};\bm{\theta}) = \bm{y}^{[\ell]}(\bm{x};\bm{\theta})+\sigma    \circ\bigg{(}\bm{W}^{[\ell+1]}\bm{y}^{[\ell]}(\bm{x};\bm{\theta})+\bm{b}^{[\ell+1]}\bigg{)}.
\end{equation*}

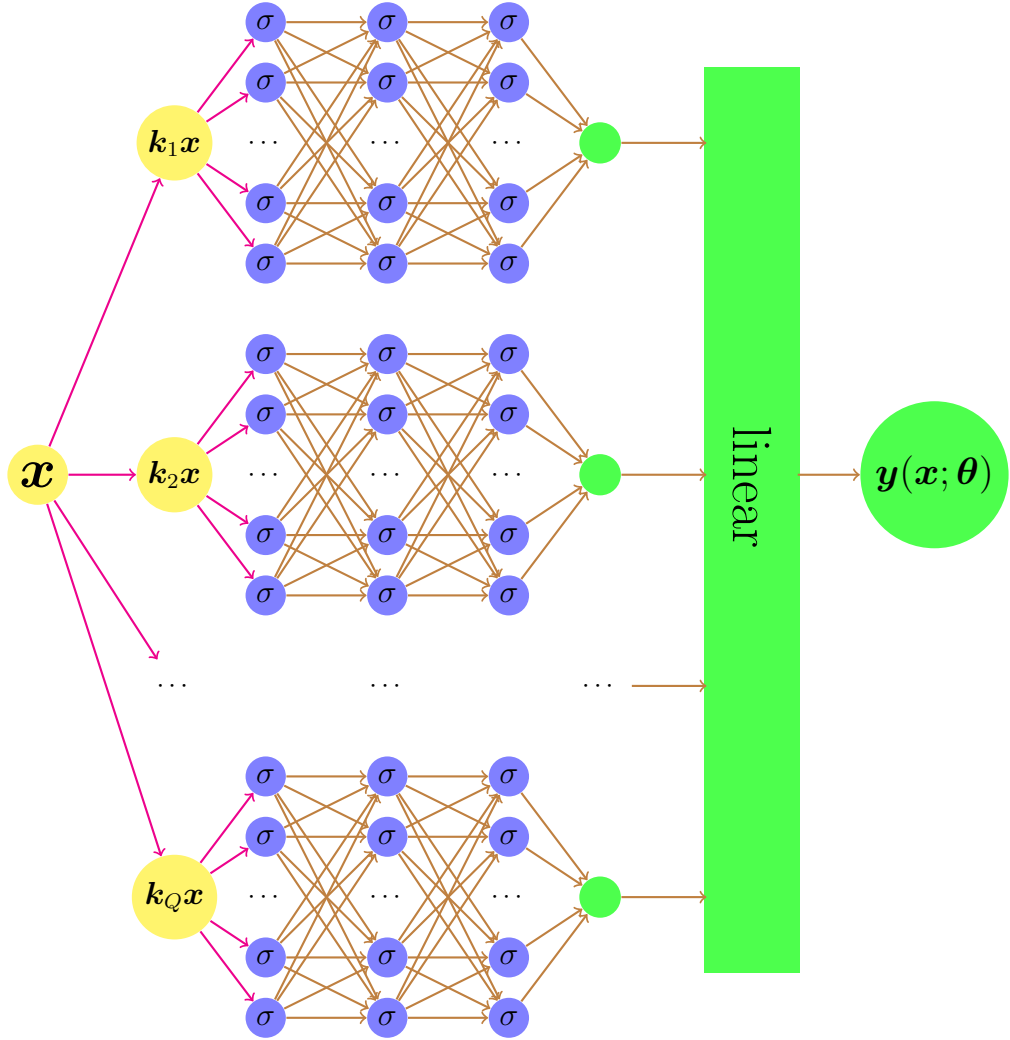
\begin{figure}[!htp]
    \begin{center}
    \begin{tikzpicture}[scale=0.8]	
        \node[circle, fill=yellow!70,inner sep=2.25pt](in) at (-2.25,8){\huge$\bm{x}$};
        \node[circle, fill=yellow!70,inner sep=2.5pt] (xh) at  (0.0,1){\large$\bm{k}_Q\bm{x}$};
        \draw[line width=0.8pt,color=magenta,->] (in) -- (xh);
        
        \node[circle, fill=blue!50,inner sep=2.25pt] (h01) at (1.5,-1){\large$\sigma$};
        \node[circle, fill=blue!50,inner sep=2.25pt] (h02) at (1.5,0){\large$\sigma$};
        \node[circle, fill=blue!0] (h03) at (1.5,1){$\cdots$};
        \node[circle, fill=blue!50,inner sep=2.25pt] (h04) at (1.5,2){\large$\sigma$};
        \node[circle, fill=blue!50,inner sep=2.25pt] (h05) at (1.5,3){\large$\sigma$};
        
        \draw[line width=0.8pt,color=magenta,->] (xh) -- (h01);
        \draw[line width=0.8pt,color=magenta,->] (xh) -- (h02);
        \draw[line width=0.8pt,color=magenta,->] (xh) -- (h04);
        \draw[line width=0.8pt,color=magenta,->] (xh) -- (h05);
        
        \node[circle, fill=blue!50,inner sep=2.25pt] (h11) at (3.5,-1.0){\large$\sigma$};
        \node[circle, fill=blue!50,inner sep=2.25pt] (h12) at (3.5,0){\large$\sigma$};
        \node[circle, fill=blue!0] (h13) at (3.5,1){$\cdots$};
        \node[circle, fill=blue!50,inner sep=2.25pt] (h14) at (3.5,2){\large$\sigma$};
        \node[circle, fill=blue!50,inner sep=2.25pt] (h15) at (3.5,3){\large$\sigma$};
        
        \draw[line width=0.8pt,color=brown,->] (h01) -- (h11);
        \draw[line width=0.8pt,color=brown,->] (h01) -- (h12);
        \draw[line width=0.8pt,color=brown,->] (h01) -- (h14);
        \draw[line width=0.8pt,color=brown,->] (h01) -- (h15);
        
        \draw[line width=0.8pt,color=brown,->] (h02) -- (h11);
        \draw[line width=0.8pt,color=brown,->] (h02) -- (h12);
        \draw[line width=0.8pt,color=brown,->] (h02) -- (h14);
        \draw[line width=0.8pt,color=brown,->] (h02) -- (h15);
        
        \draw[line width=0.8pt,color=brown,->] (h04) -- (h11);
        \draw[line width=0.8pt,color=brown,->] (h04) -- (h12);
        \draw[line width=0.8pt,color=brown,->] (h04) -- (h14);
        \draw[line width=0.8pt,color=brown,->] (h04) -- (h15);
        
        \draw[line width=0.8pt,color=brown,->] (h05) -- (h11);
        \draw[line width=0.8pt,color=brown,->] (h05) -- (h12);
        \draw[line width=0.8pt,color=brown,->] (h05) -- (h14);
        \draw[line width=0.8pt,color=brown,->] (h05) -- (h15);
        
        \node[circle, fill=blue!50,inner sep=2.25pt] (h21) at (5.5,-1.0){\large$\sigma$};
        \node[circle, fill=blue!50,inner sep=2.25pt] (h22) at (5.5,0){\large$\sigma$};
        \node[circle, fill=blue!0] (h23) at (5.5,1){$\cdots$};
        \node[circle, fill=blue!50,inner sep=2.25pt] (h24) at (5.5,2){\large$\sigma$};
        \node[circle, fill=blue!50,inner sep=2.25pt] (h25) at (5.5,3){\large$\sigma$};
        \node[circle, fill=green!70,inner sep=5.5pt] (uh) at (7.0,1){};
        
        \draw[line width=0.8pt,color=brown,->] (h11) -- (h21);
        \draw[line width=0.8pt,color=brown,->] (h11) -- (h22);
        \draw[line width=0.8pt,color=brown,->] (h11) -- (h24);
        \draw[line width=0.8pt,color=brown,->] (h11) -- (h25);
        
        \draw[line width=0.8pt,color=brown,->] (h12) -- (h21);
        \draw[line width=0.8pt,color=brown,->] (h12) -- (h22);
        \draw[line width=0.8pt,color=brown,->] (h12) -- (h24);
        \draw[line width=0.8pt,color=brown,->] (h12) -- (h25);
        
        \draw[line width=0.8pt,color=brown,->] (h14) -- (h21);
        \draw[line width=0.8pt,color=brown,->] (h14) -- (h22);
        \draw[line width=0.8pt,color=brown,->] (h14) -- (h24);
        \draw[line width=0.8pt,color=brown,->] (h14) -- (h25);
        
        \draw[line width=0.8pt,color=brown,->] (h15) -- (h21);
        \draw[line width=0.8pt,color=brown,->] (h15) -- (h22);
        \draw[line width=0.8pt,color=brown,->] (h15) -- (h24);
        \draw[line width=0.8pt,color=brown,->] (h15) -- (h25);
        
        \draw[line width=0.8pt,color=brown,->] (h21) -- (uh);
        \draw[line width=0.8pt,color=brown,->] (h22) -- (uh);
        \draw[line width=0.8pt,color=brown,->] (h24) -- (uh);
        \draw[line width=0.8pt,color=brown,->] (h25) -- (uh);

        \node[circle, fill=blue!0,inner sep=3.5pt] (xi) at (0.0,4.5){$\cdots$};
	\node[circle, fill=blue!0,inner sep=3.5pt] (ihidden) at (3.5,4.5){$\cdots$};
	\node[circle, fill=blue!0,inner sep=3.5pt] (ui) at (7.0,4.5){$\cdots$};
        \draw[line width=0.8pt,color=magenta,->] (in) -- (xi);

        \node[circle, fill=yellow!70,inner sep=2.25pt] (xj) at  (0.0,8){\large$\bm{k}_2\bm{x}$};
        \draw[line width=0.8pt,color=magenta,->] (in) -- (xj);

        \node[circle, fill=blue!50,inner sep=2.25pt] (j01) at (1.5,6){\large$\sigma$};
        \node[circle, fill=blue!50,inner sep=2.25pt] (j02) at (1.5,7){\large$\sigma$};
        \node[circle, fill=blue!0] (j03) at (1.5,8){$\cdots$};
        \node[circle, fill=blue!50,inner sep=2.25pt] (j04) at (1.5,9){\large$\sigma$};
        \node[circle, fill=blue!50,inner sep=2.25pt] (j05) at (1.5,10){\large$\sigma$};

        \draw[line width=0.8pt,color=magenta,->] (xj) -- (j01);
        \draw[line width=0.8pt,color=magenta,->] (xj) -- (j02);
        \draw[line width=0.8pt,color=magenta,->] (xj) -- (j04);
        \draw[line width=0.8pt,color=magenta,->] (xj) -- (j05);

        \node[circle, fill=blue!50,inner sep=2.25pt] (j11) at (3.5,6){\large$\sigma$};
        \node[circle, fill=blue!50,inner sep=2.25pt] (j12) at (3.5,7){\large$\sigma$};
        \node[circle, fill=blue!0] (j13) at (3.5,8){$\cdots$};
        \node[circle, fill=blue!50,inner sep=2.25pt] (j14) at (3.5,9){\large$\sigma$};
        \node[circle, fill=blue!50,inner sep=2.25pt] (j15) at (3.5,10){\large$\sigma$};

        \draw[line width=0.8pt,color=brown,->] (j01) -- (j11);
        \draw[line width=0.8pt,color=brown,->] (j01) -- (j12);
        \draw[line width=0.8pt,color=brown,->] (j01) -- (j14);
        \draw[line width=0.8pt,color=brown,->] (j01) -- (j15);
        
        \draw[line width=0.8pt,color=brown,->] (j02) -- (j11);
        \draw[line width=0.8pt,color=brown,->] (j02) -- (j12);
        \draw[line width=0.8pt,color=brown,->] (j02) -- (j14);
        \draw[line width=0.8pt,color=brown,->] (j02) -- (j15);
        
        \draw[line width=0.8pt,color=brown,->] (j04) -- (j11);
        \draw[line width=0.8pt,color=brown,->] (j04) -- (j12);
        \draw[line width=0.8pt,color=brown,->] (j04) -- (j14);
        \draw[line width=0.8pt,color=brown,->] (j04) -- (j15);
        
        \draw[line width=0.8pt,color=brown,->] (j05) -- (j11);
        \draw[line width=0.8pt,color=brown,->] (j05) -- (j12);
        \draw[line width=0.8pt,color=brown,->] (j05) -- (j14);
        \draw[line width=0.8pt,color=brown,->] (j05) -- (j15);

        \node[circle, fill=blue!50,inner sep=2.25pt] (j21) at (5.5,6){\large$\sigma$};
        \node[circle, fill=blue!50,inner sep=2.25pt] (j22) at (5.5,7){\large$\sigma$};
        \node[circle, fill=blue!0] (j23) at (5.5,8){$\cdots$};
        \node[circle, fill=blue!50,inner sep=2.25pt] (j24) at (5.5,9){\large$\sigma$};
        \node[circle, fill=blue!50,inner sep=2.25pt] (j25) at (5.5,10){\large$\sigma$};

        \draw[line width=0.8pt,color=brown,->] (j11) -- (j21);
        \draw[line width=0.8pt,color=brown,->] (j11) -- (j22);
        \draw[line width=0.8pt,color=brown,->] (j11) -- (j24);
        \draw[line width=0.8pt,color=brown,->] (j11) -- (j25);
        
        \draw[line width=0.8pt,color=brown,->] (j12) -- (j21);
        \draw[line width=0.8pt,color=brown,->] (j12) -- (j22);
        \draw[line width=0.8pt,color=brown,->] (j12) -- (j24);
        \draw[line width=0.8pt,color=brown,->] (j12) -- (j25);
        
        \draw[line width=0.8pt,color=brown,->] (j14) -- (j21);
        \draw[line width=0.8pt,color=brown,->] (j14) -- (j22);
        \draw[line width=0.8pt,color=brown,->] (j14) -- (j24);
        \draw[line width=0.8pt,color=brown,->] (j14) -- (j25);
        
        \draw[line width=0.8pt,color=brown,->] (j15) -- (j21);
        \draw[line width=0.8pt,color=brown,->] (j15) -- (j22);
        \draw[line width=0.8pt,color=brown,->] (j15) -- (j24);
        \draw[line width=0.8pt,color=brown,->] (j15) -- (j25);

        \node[circle, fill=green!70,inner sep=5.5pt] (uj) at (7.0,8){};
        \draw[line width=0.8pt,color=brown,->] (j21) -- (uj);
        \draw[line width=0.8pt,color=brown,->] (j22) -- (uj);
        \draw[line width=0.8pt,color=brown,->] (j24) -- (uj);
        \draw[line width=0.8pt,color=brown,->] (j25) -- (uj);

        \node[circle, fill=yellow!70,inner sep=2.25pt] (xk) at  (0.0,13.5){\large$\bm{k}_1\bm{x}$};
        \draw[line width=0.8pt,color=magenta,->] (in) -- (xk);
        
        \node[circle, fill=blue!50,inner sep=2.25pt] (k01) at (1.5,11.5){\large$\sigma$};
        \node[circle, fill=blue!50,inner sep=2.25pt] (k02) at (1.5,12.5){\large$\sigma$};
        \node[circle, fill=blue!0,inner sep=2.25pt] (k03) at (1.5,13.5){$\cdots$};
        \node[circle, fill=blue!50,inner sep=2.25pt] (k04) at (1.5,14.5){\large$\sigma$};
        \node[circle, fill=blue!50,inner sep=2.25pt] (k05) at (1.5,15.5){\large$\sigma$};
        
        \draw[line width=0.8pt,color=magenta,->] (xk) -- (k01);
        \draw[line width=0.8pt,color=magenta,->] (xk) -- (k02);
        \draw[line width=0.8pt,color=magenta,->] (xk) -- (k04);
        \draw[line width=0.8pt,color=magenta,->] (xk) -- (k05);
        
        \node[circle, fill=blue!50,inner sep=2.25pt] (k11) at (3.5,11.5){\large$\sigma$};
        \node[circle, fill=blue!50,inner sep=2.25pt] (k12) at (3.5,12.5){\large$\sigma$};
        \node[circle, fill=blue!0] (k13) at (3.5,13.5){$\cdots$};
        \node[circle, fill=blue!50,inner sep=2.25pt] (k14) at (3.5,14.5){\large$\sigma$};
        \node[circle, fill=blue!50,inner sep=2.25pt] (k15) at (3.5,15.5){\large$\sigma$};
        
        \draw[line width=0.8pt,color=brown,->] (k01) -- (k11);
        \draw[line width=0.8pt,color=brown,->] (k01) -- (k12);
        \draw[line width=0.8pt,color=brown,->] (k01) -- (k14);
        \draw[line width=0.8pt,color=brown,->] (k01) -- (k15);
        
        \draw[line width=0.8pt,color=brown,->] (k02) -- (k11);
        \draw[line width=0.8pt,color=brown,->] (k02) -- (k12);
        \draw[line width=0.8pt,color=brown,->] (k02) -- (k14);
        \draw[line width=0.8pt,color=brown,->] (k02) -- (k15);
        
        \draw[line width=0.8pt,color=brown,->] (k04) -- (k11);
        \draw[line width=0.8pt,color=brown,->] (k04) -- (k12);
        \draw[line width=0.8pt,color=brown,->] (k04) -- (k14);
        \draw[line width=0.8pt,color=brown,->] (k04) -- (k15);
        
        \draw[line width=0.8pt,color=brown,->] (k05) -- (k11);
        \draw[line width=0.8pt,color=brown,->] (k05) -- (k12);
        \draw[line width=0.8pt,color=brown,->] (k05) -- (k14);
        \draw[line width=0.8pt,color=brown,->] (k05) -- (k15);
        
        \node[circle, fill=blue!50,inner sep=2.25pt] (k21) at (5.5,11.5){\large$\sigma$};
        \node[circle, fill=blue!50,inner sep=2.25pt] (k22) at (5.5,12.5){\large$\sigma$};
        \node[circle, fill=blue!0] (k23) at (5.5,13.5){$\cdots$};
        \node[circle, fill=blue!50,inner sep=2.25pt] (k24) at (5.5,14.5){\large$\sigma$};
        \node[circle, fill=blue!50,inner sep=2.25pt] (k25) at (5.5,15.5){\large$\sigma$};
        \node[circle, fill=green!70,inner sep=5.5pt] (uk) at (7.0,13.5){};
        
        \draw[line width=0.8pt,color=brown,->] (k11) -- (k21);
        \draw[line width=0.8pt,color=brown,->] (k11) -- (k22);
        \draw[line width=0.8pt,color=brown,->] (k11) -- (k24);
        \draw[line width=0.8pt,color=brown,->] (k11) -- (k25);
        
        \draw[line width=0.8pt,color=brown,->] (k12) -- (k21);
        \draw[line width=0.8pt,color=brown,->] (k12) -- (k22);
        \draw[line width=0.8pt,color=brown,->] (k12) -- (k24);
        \draw[line width=0.8pt,color=brown,->] (k12) -- (k25);
        
        \draw[line width=0.8pt,color=brown,->] (k14) -- (k21);
        \draw[line width=0.8pt,color=brown,->] (k14) -- (k22);
        \draw[line width=0.8pt,color=brown,->] (k14) -- (k24);
        \draw[line width=0.8pt,color=brown,->] (k14) -- (k25);
        
        \draw[line width=0.8pt,color=brown,->] (k15) -- (k21);
        \draw[line width=0.8pt,color=brown,->] (k15) -- (k22);
        \draw[line width=0.8pt,color=brown,->] (k15) -- (k24);
        \draw[line width=0.8pt,color=brown,->] (k15) -- (k25);
        
        \draw[line width=0.8pt,color=brown,->] (k21) -- (uk);
        \draw[line width=0.8pt,color=brown,->] (k22) -- (uk);
        \draw[line width=0.8pt,color=brown,->] (k24) -- (uk);
        \draw[line width=0.8pt,color=brown,->] (k25) -- (uk);
        
        \node[rectangle, fill=green!70, minimum width=1.25cm, minimum height=12cm, inner sep=0.5pt] (lin) at (9.5,7.25){\rotatebox{-90}{\huge{linear~~~~}}};
        
        \draw[line width=0.8pt,color=brown,->] (uh) -- (8.75,1.0);
        \draw[line width=0.8pt,color=brown,->] (ui) -- (8.75,4.5);
        \draw[line width=0.8pt,color=brown,->] (uj) -- (8.75,8);
        \draw[line width=0.8pt,color=brown,->] (uk) -- (8.75,13.5);

        \node[circle, fill=green!70,inner sep=3.5pt] (y) at (12.5,8.0){\Large{$\bm{y}(\bm{x};\bm{\theta})$}};

        \draw[line width=0.8pt,color=brown,->] (10.25,8.0) -- (y);
    \end{tikzpicture}
    \end{center}
    \caption{A schematic diagram of MscaleDNN with $Q$ subnetworks, $\sigma$ stands for the activation function}
    \label{fig2mscalednn}
\end{figure}

As we are aware, a normal DNN model is capable of providing a satisfactory solution for general problems. However, it will encounter troublesome difficulty to solve multi-scale problems with high-frequency components. Recently, a MscaleDNN architecture has shown its remarkable performance to deal with high-frequency problems by converting original data to a low-frequency space \cite{liu2020multi,wang2020multiscale,li2020elliptic,wang2020eigenvector}. A schematic diagram of MscaleDNN with $Q$ subnetworks is depicted in Fig. \ref{fig2mscalednn}.

The detailed procedure of MscaleDNN is described in the following.
\begin{enumerate}
    \item Generating a scale vector or matrix with $Q$ parts
    \begin{equation}\label{separate}
        \Lambda = (\bm{k}_1, \bm{k}_2, \bm{k}_3\cdots,\bm{k}_{Q-1},\bm{k}_Q)^T,
    \end{equation}
    where $\bm{k}_i(i=1,2,\ldots,Q)$ is a scalar or matrix (trainable or untrainable). 
    \item Converting the input data $\bm{x}$ into $\tilde{\bm{x}} = \Lambda\odot \bm{x}$ with $\odot$ being the Hadamard product, then feeding $\tilde{\bm{x}}$ into the pipeline of MscaleDNN. It is
    \begin{equation}
    \begin{cases}
     \hat{\bm{x}}= \bm{k}_i\bm{x}\\
     \boldsymbol{F}_{i}(\bm{x})  =\mathcal{F C N}_{i}\left(\hat{\bm{x}}\right)
    \end{cases} \quad i=1,2, \ldots, Q,
    \end{equation}
    where $\mathcal{F C N}_{i}$ stands for the $i_{th}$ fully connected subnetwork and $\boldsymbol{F}_{i}$ is its output.
    \item Obtaining the result of MscaleDNN by aggregating linearly the output of all subnetworks, each scale input goes through a subnetwork. It is
    \begin{equation}\label{out2mscalednn}
     \boldsymbol{NN}(\bm{x}) =\boldsymbol{W}_O \cdot\left[\boldsymbol{F}_{1}(\bm{x}), \boldsymbol{F}_{2}(\bm{x}), \cdots, \boldsymbol{F}_{Q}(\bm{x})\right]+\boldsymbol{b}_O,
    \end{equation}
    where $\boldsymbol{W}_O$ and $\boldsymbol{b}_O$ stand for the weights and biases of the last linear layer, respectively.
\end{enumerate}
From the perspective of Fourier transformation and decomposition, the first layer of the MscaleDNN model will be treated as a series of basis in Fourier space and its output is the combination of basis functions \cite{liu2020multi,li2020elliptic,wang2020eigenvector}.

\subsection{Overview of Physics-Informed Neural Networks}
In the scope of PINN, a type of PDE governed by parameters as the toy to show its implementation, it is
\begin{equation}\label{eq2PINN}
\begin{aligned}
&\mathcal{N}_{\bm{\lambda}}[\hat{u}(\bm{x})]=\hat{f}(\bm{x}), \quad \bm{x} \in \Omega\\
&\mathcal{B}\hat{u}\left(\bm{x}\right)=\hat{g}(\bm{x}), \quad \quad\bm{x} \in \partial\Omega
\end{aligned}
\end{equation}
in which $\mathcal{N}_{\bm{\lambda}}$ stands for the linear or nonlinear differential operator with parameters $\bm{\lambda}$, $\mathcal{B}$ is the boundary operator, such as Dirichlet,
Neumann, periodic boundary conditions, or a mixed form of them.
$\Omega $ and  $\partial \Omega $ respectively illustrate the zone of interest and its border.
For approximating the solution of the multi-scale PDE, a multi-scale deep neural network is used.
In classical PINN, the ideal parameters of the DNN can be obtained by
minimizing the following composite loss function
\begin{equation}\label{loss2PINN}
Loss=Loss_{R} +\gamma Loss_{B}
\end{equation}
with
\begin{equation}\label{subloss2PINN}
\begin{aligned}
&Loss_{R} = \frac{1}{N_R}\sum_{i=1}^{N_R}\left\| \mathcal{N}_{\bm{\lambda}}[u_{NN}(\bm{x}^i_I)]-\hat{f}(\bm{x}^i_I)\right\|^2 \\
&Loss_{B} = \frac{1}{N_B}\sum_{j=1}^{N_B}\bigg{\|}\mathcal{B}u_{NN}\left(\bm{x}^j_B\right)-\hat{g}(\bm{x}^j_B)\bigg{\|}^2
\end{aligned}
\end{equation}
where $\gamma>0$ is used to control the contribution for the corresponding loss term. $Loss_R$ and $Loss_B$ depict the residual of the governing equations and the loss on the boundary condition, respectively. If some additional observed data are available inside the interested domain, then a loss term indicating the mismatch between the predictions produced by DNN and the observations can be taken into account
\begin{equation}
    Loss_{D} = \frac{1}{N_D}\sum_{i=1}^{N_D}\bigg{\|}u_{NN}(\bm{x}^i)-u_{Data}^i\bigg{\|}^2.
\end{equation}

\section{Fourier-based mixed PINN to solve multi-scale problem}\label{sec:03}
In this section, the unified architecture of FMPINN is proposed to overcome the adverse effect of derivative for rough coefficient $A^{\varepsilon}$ by embracing a multi-output neural network with an equivalent reduced-order formulation of the multi-scale problem \eqref{eq:multiscale}.

\subsection{Failure of classical PINN}
Despite the success of various PINN models in studying ordinary and partial differential equations, it has been observed in \cite{leung2022nh} that the classical PINN approach fails to provide accurate predictions for multi-scale PDEs\eqref{eq:multiscale}. Furthermore, we find that a direct application of the PINN with multi-scale DNN framework on solving \eqref{eq:multiscale} still cannot provide a satisfactory solution, because of the ill-posed NTK matrix caused by rough coefficient $A^{\varepsilon}$. 
For example, let us consider the following one-dimensional elliptic equation  with a homogeneous Dirichlet boundary in $\Omega=[0,1]$:
\begin{equation*}
\begin{cases}
\displaystyle-\frac{d}{dx}\bigg{(}A^{\varepsilon}(x)\frac{d}{dx} u^{\varepsilon}(x)\bigg{)} = 5\cos(\pi x)\\
~~~~~~~~~~~u^{\varepsilon}(0)=u^{\varepsilon}(1) = 0
\end{cases},
\end{equation*}
in which $\displaystyle A^{\varepsilon}(x)=\frac{1+x^2}{2+\sin(2\pi x/\varepsilon)}$ with $\varepsilon>0$ being a small constant. 

\begin{figure}[H]
	\centering    
	\subfigure[$A^{\varepsilon}$ for $\varepsilon=1/32$] {
		\label{Aeps32}     
		\includegraphics[scale=0.33]{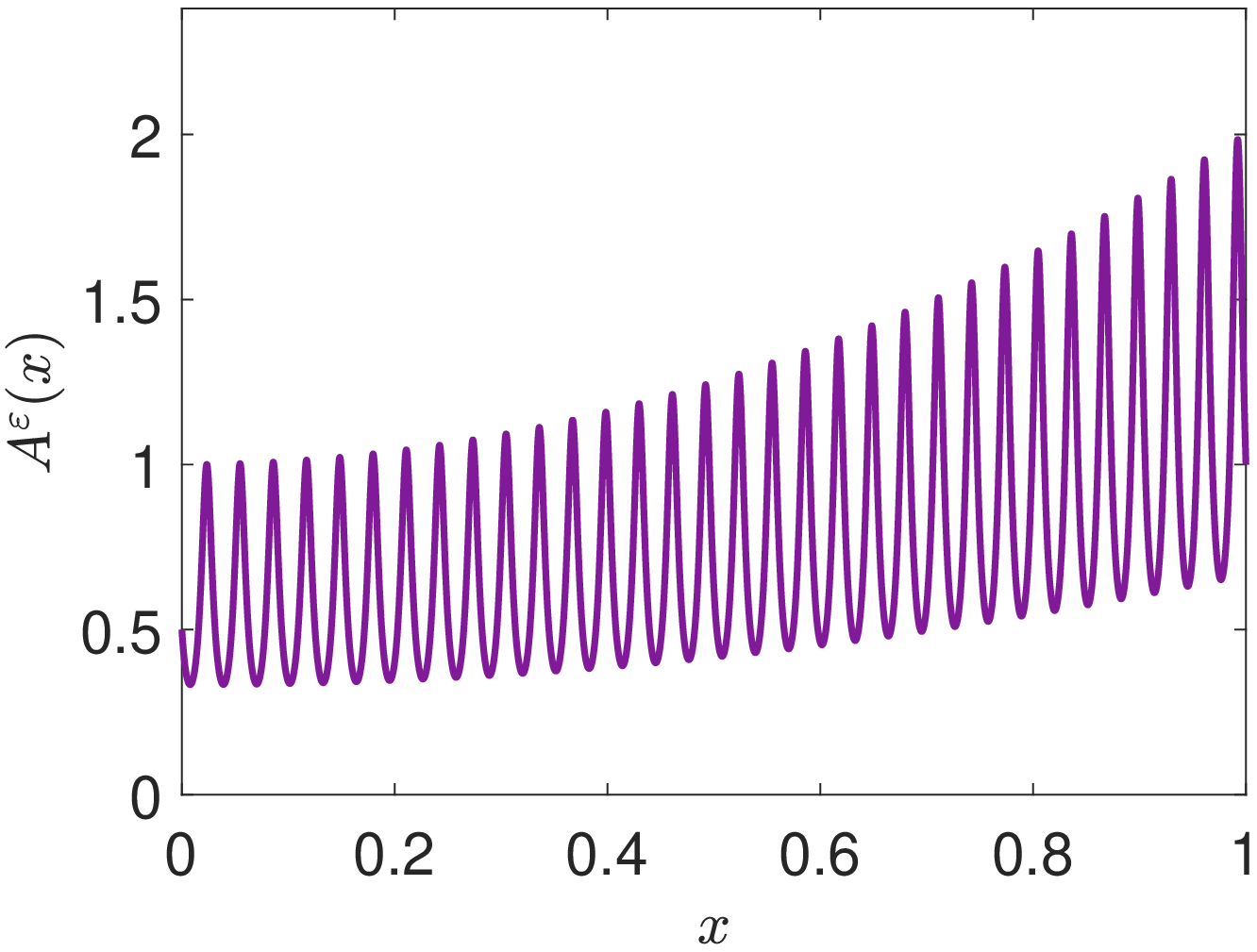}  
	}          
	\subfigure[$u^{\varepsilon}$ for $\varepsilon=1/32$] {
		\label{Ueps32}     
		\includegraphics[scale=0.33]{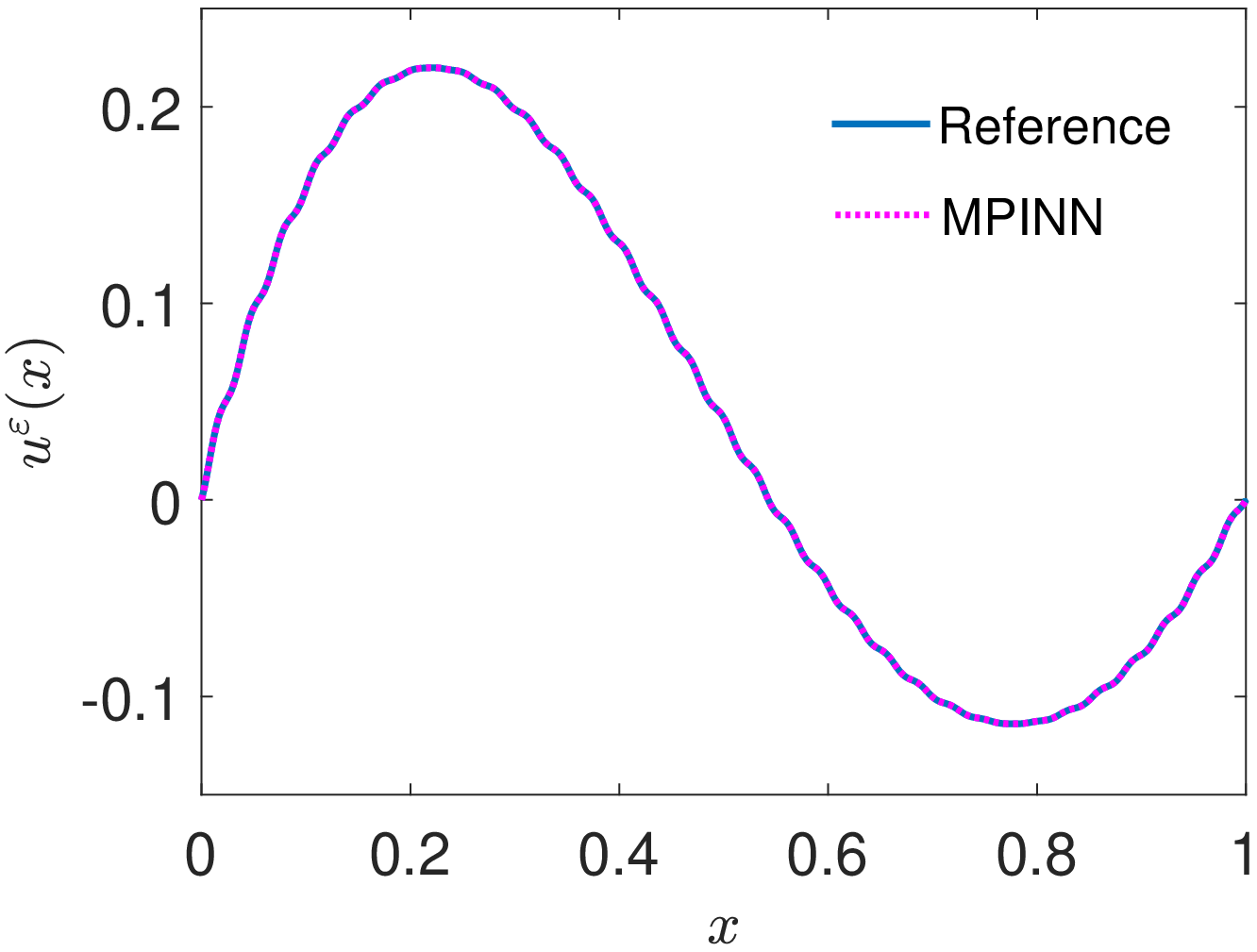}  
	}
	\subfigure[$l^2$ relative error for $\varepsilon=1/32$] {
		\label{rel32}     
		\includegraphics[scale=0.33]{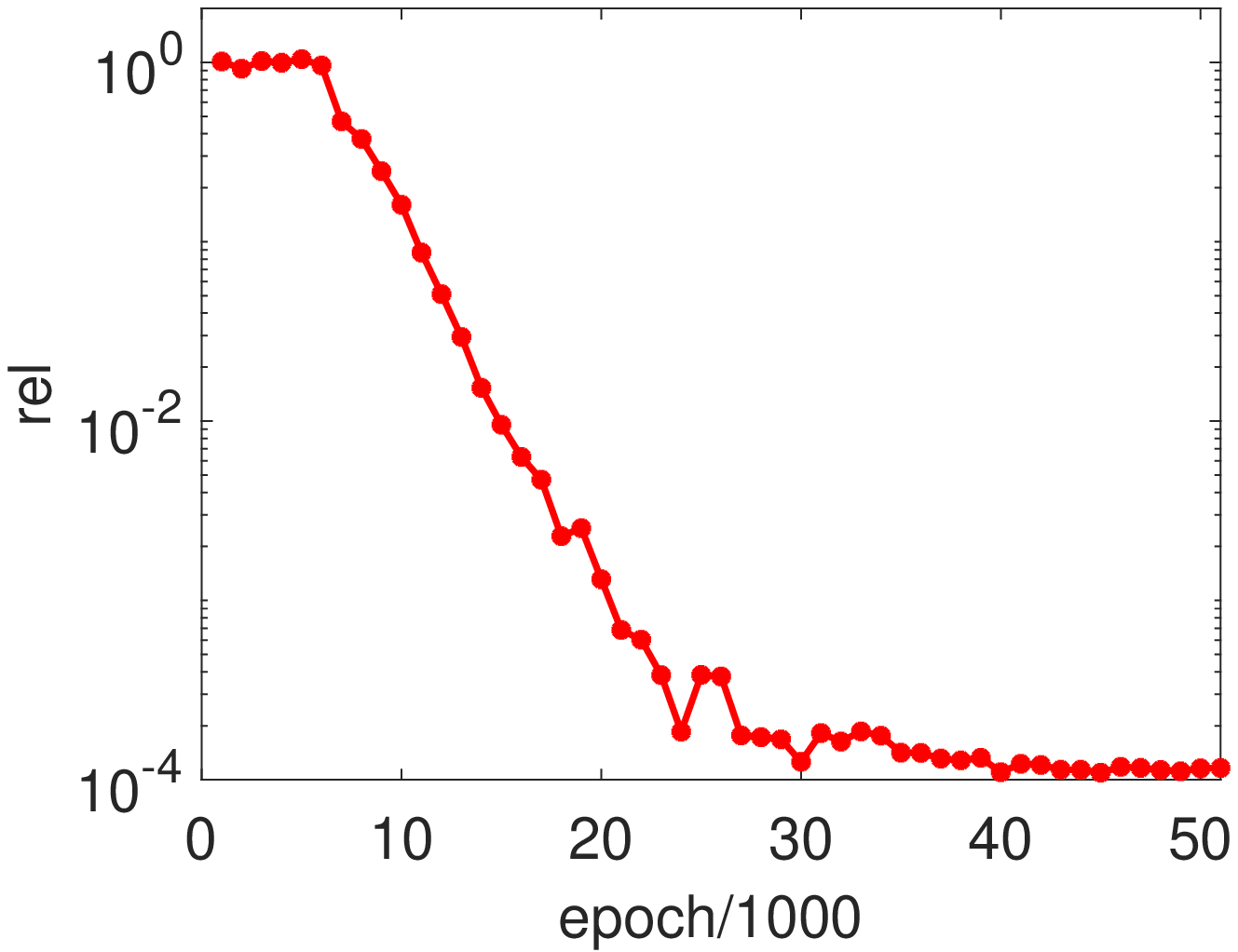}  
	}
	\subfigure[$A^{\varepsilon}$ for $\varepsilon=1/64$] {
		\label{Aeps64}     
		\includegraphics[scale=0.33]{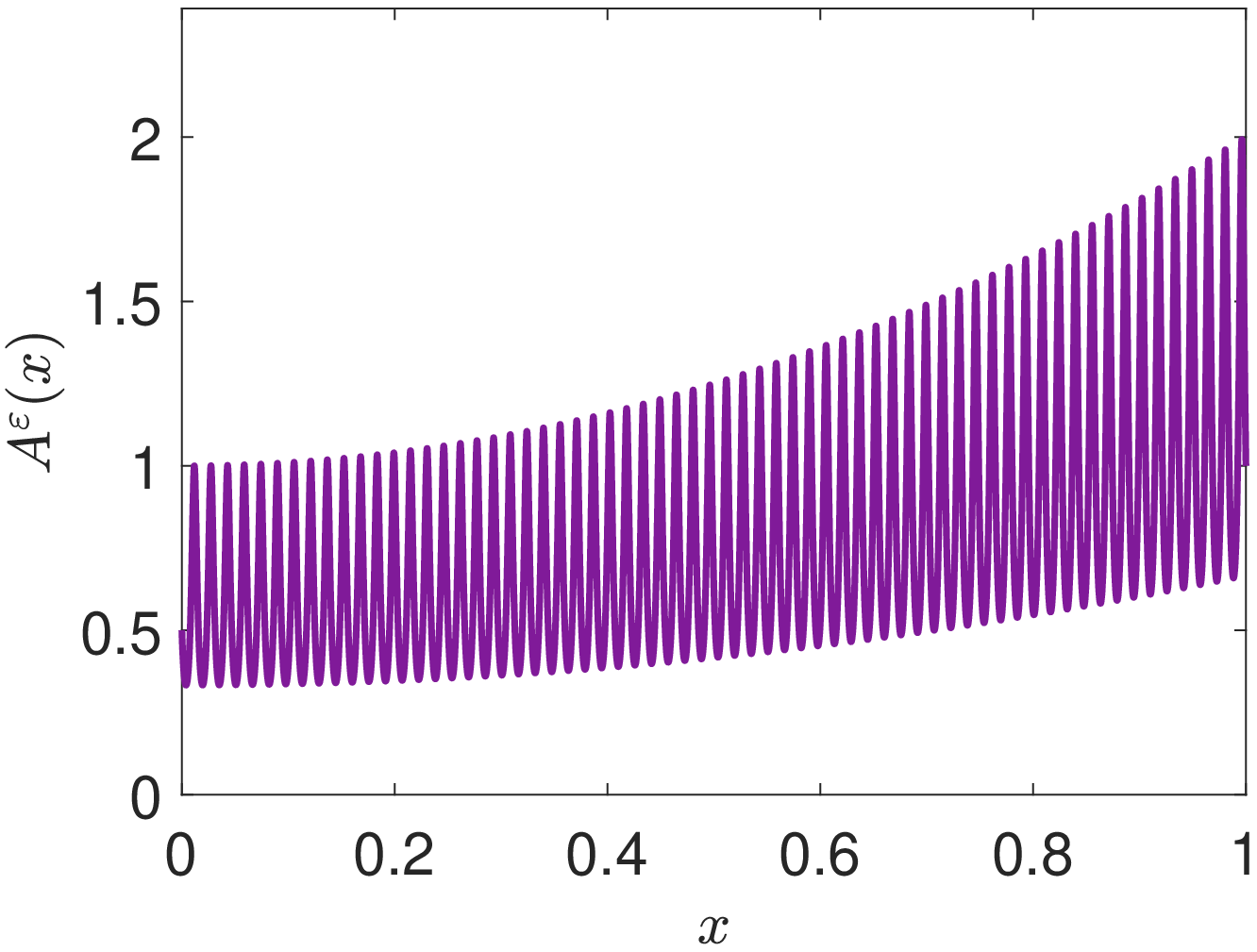}  
	}          
	\subfigure[$u^{\varepsilon}$ for $\varepsilon=1/64$] {
		\label{Ueps64}     
		\includegraphics[scale=0.33]{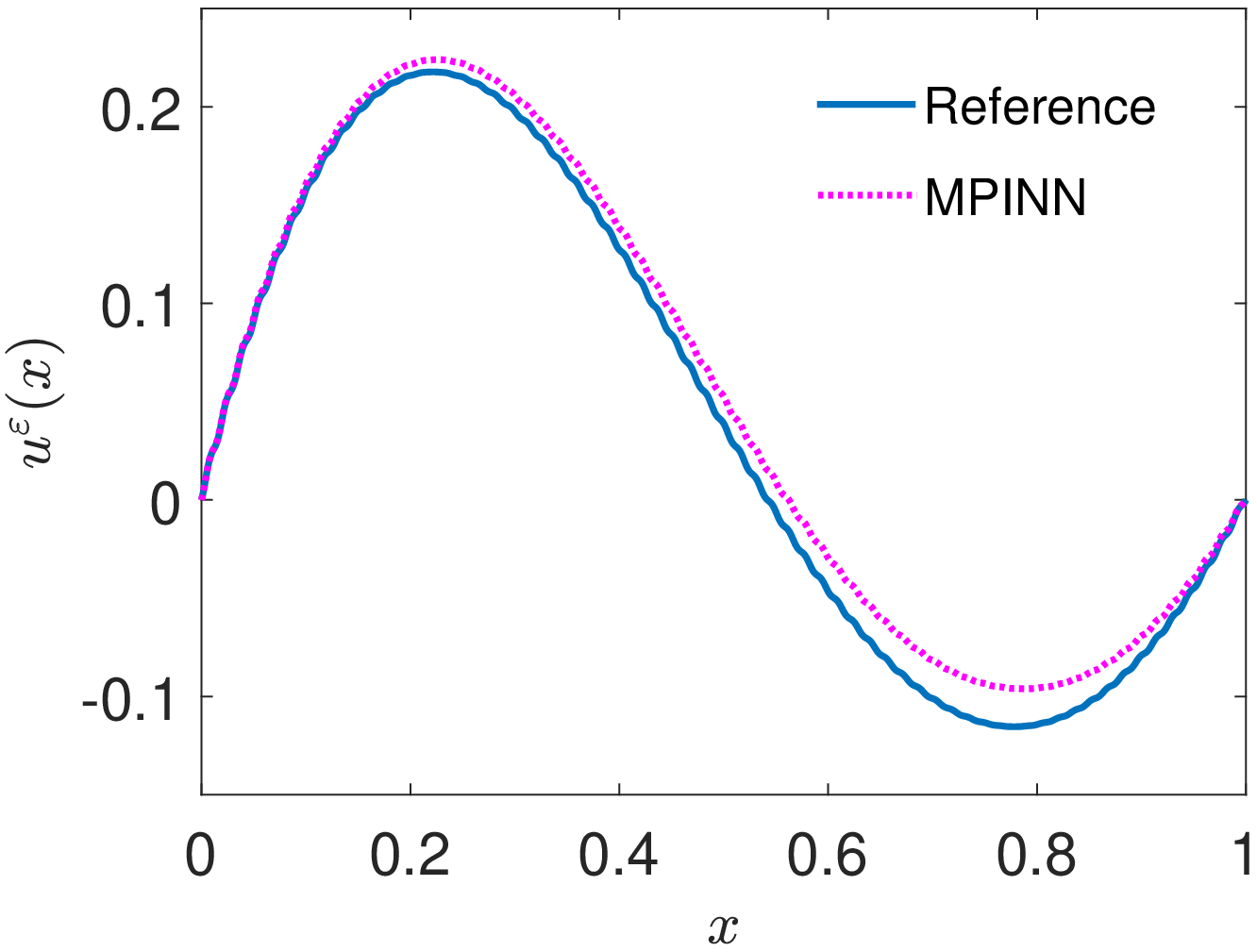}  
	}
	\subfigure[$l^2$ relative error for $\varepsilon=1/64$] {
		\label{rel64}     
		\includegraphics[scale=0.33]{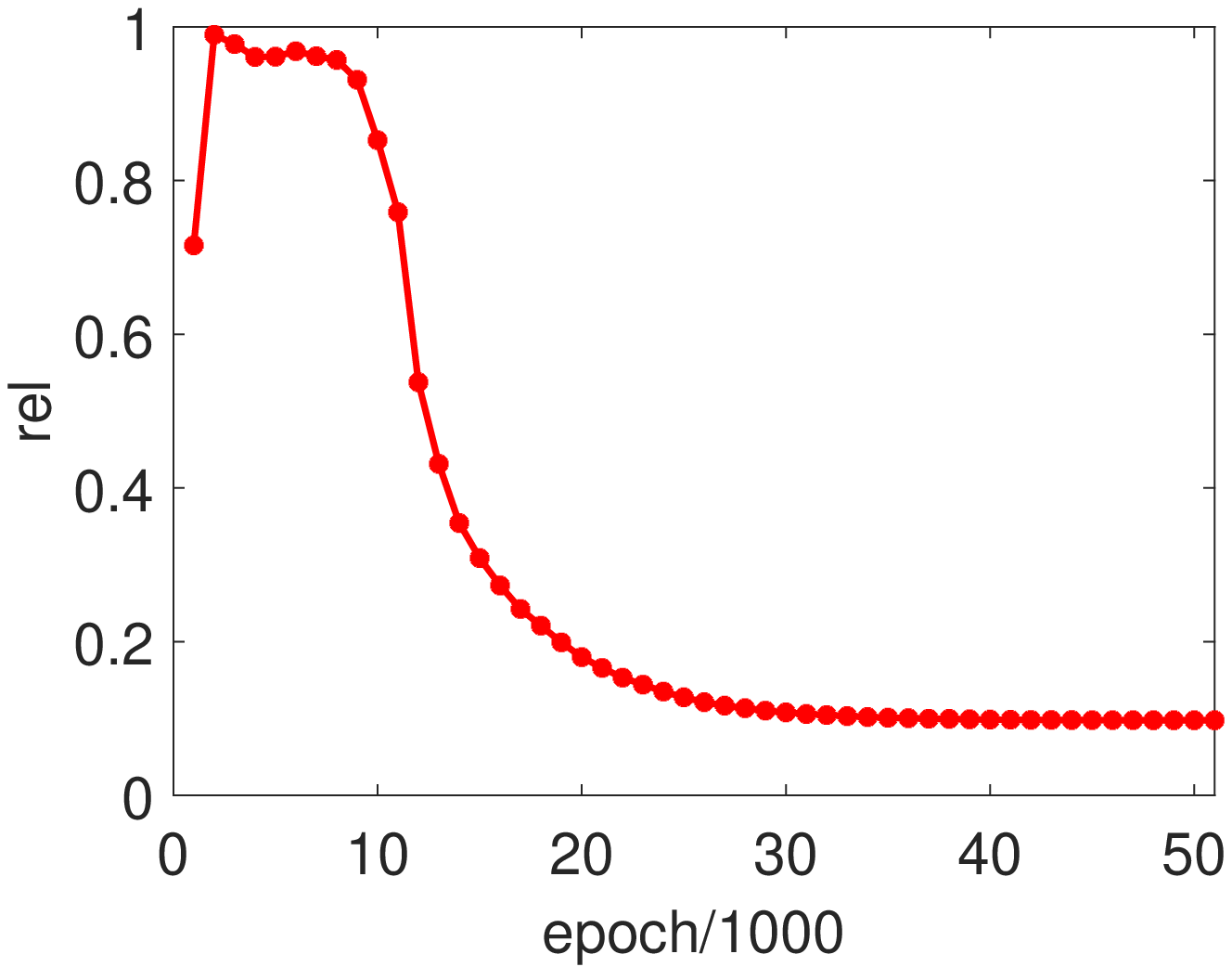}  
	}
	\subfigure[$A^{\varepsilon}$ for $\varepsilon=1/128$] {
		\label{Aeps128}     
		\includegraphics[scale=0.33]{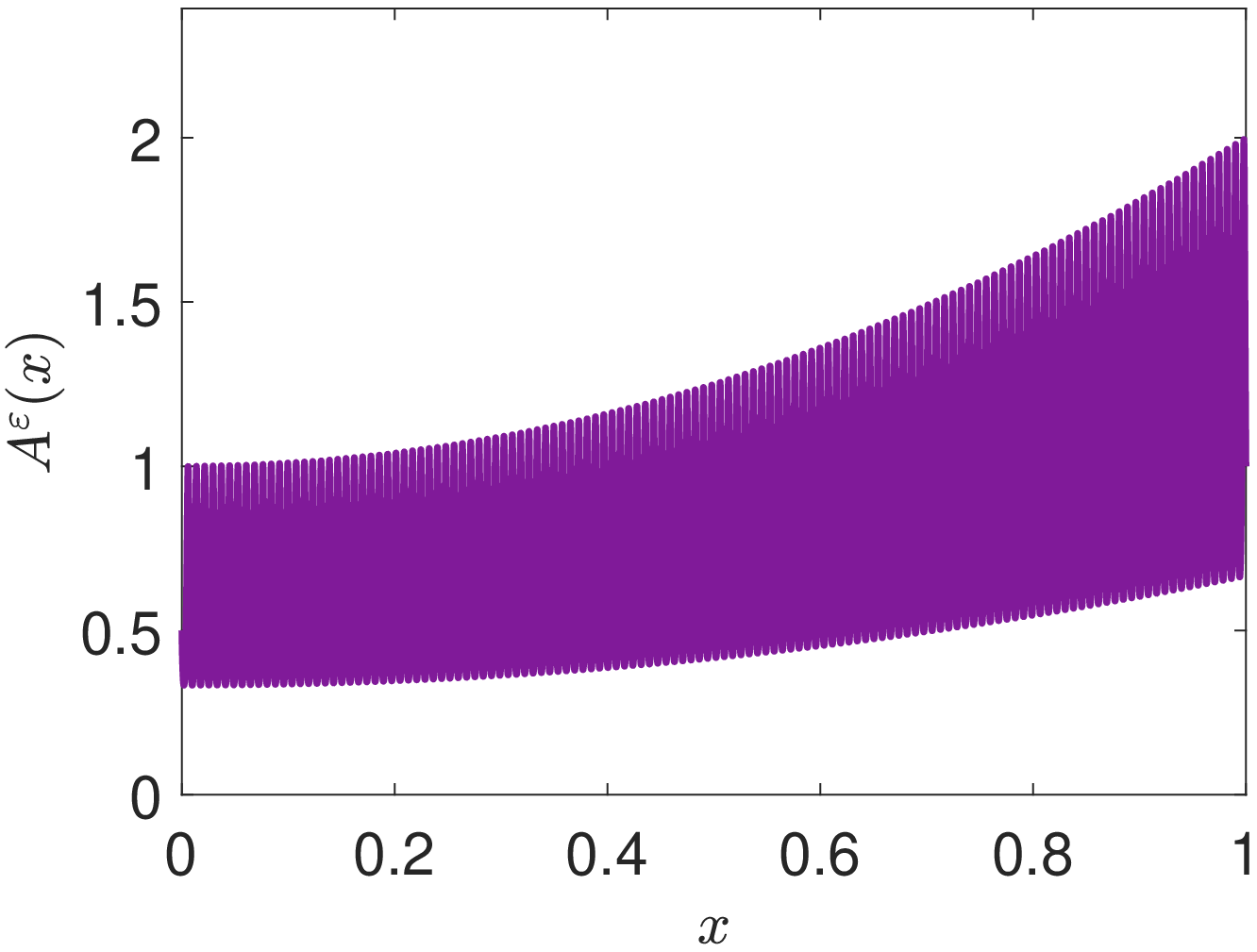}  
	}          
	\subfigure[$u^{\varepsilon}$ for $\varepsilon=1/128$] {
		\label{Ueps128}     
		\includegraphics[scale=0.33]{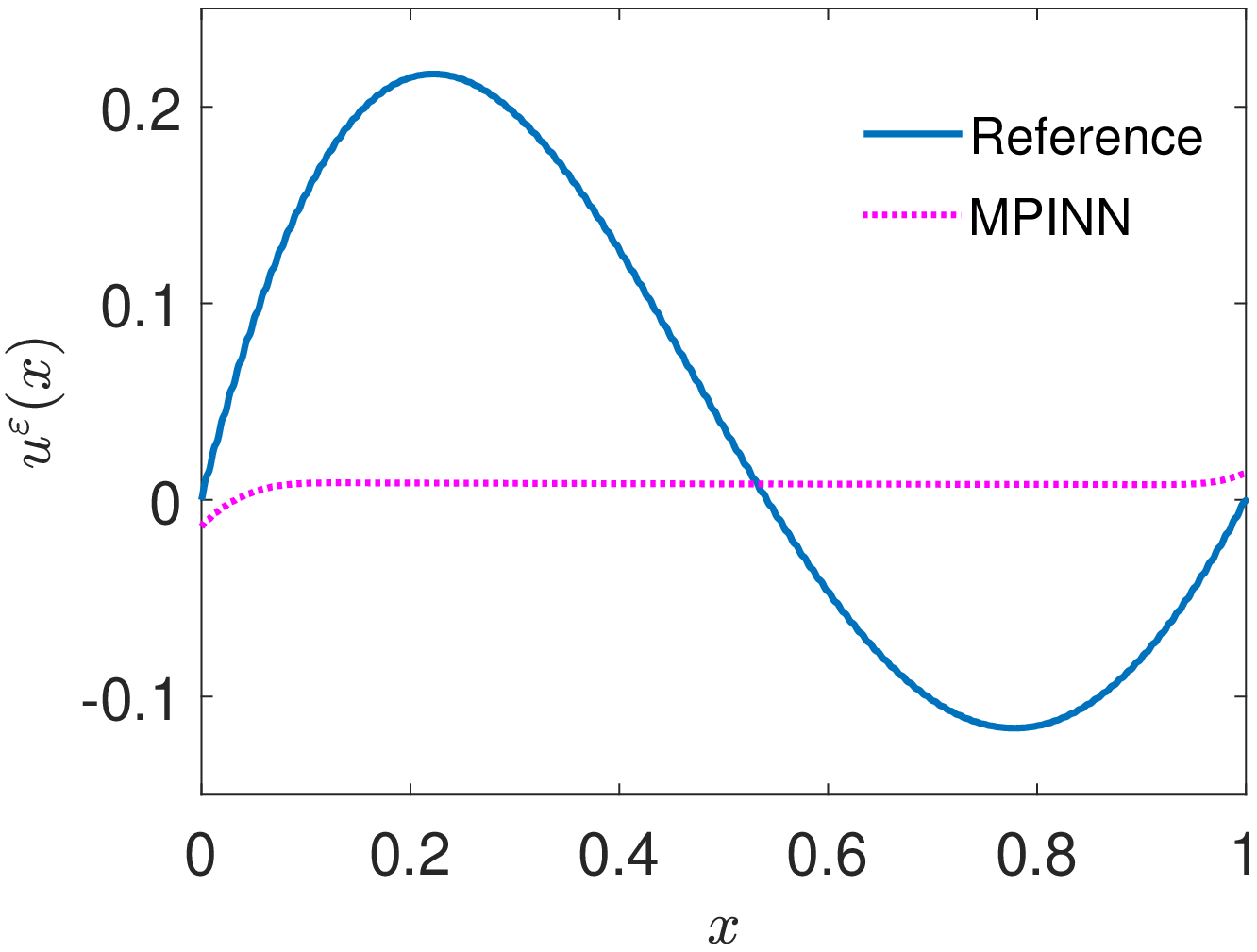}  
	}
	\subfigure[$l^2$ relative error for $\varepsilon=1/128$] {
		\label{rel128}     
		\includegraphics[scale=0.33]{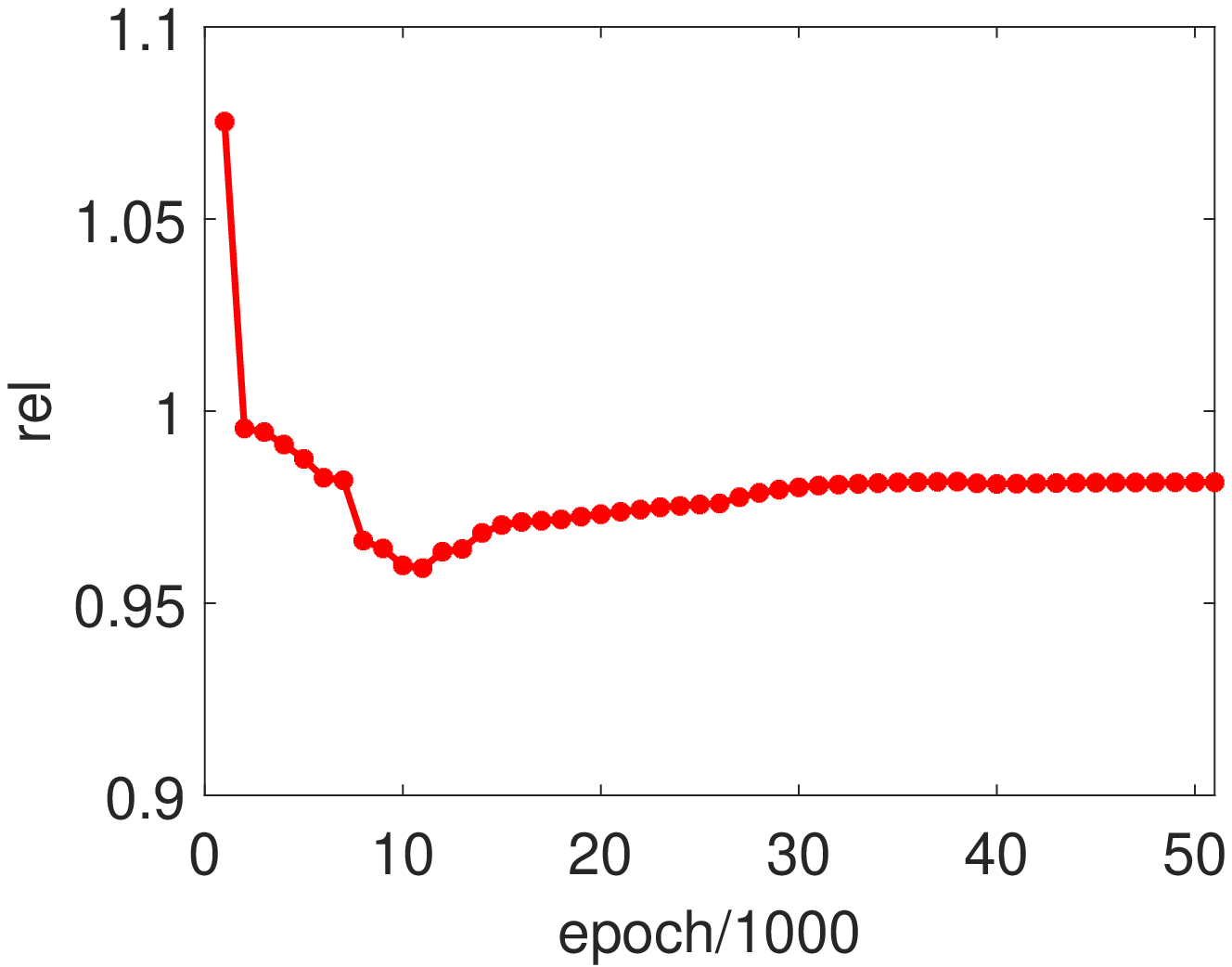}  
	}
\caption{Left: the rough coefficient $A^{\varepsilon}$. Middle: the MPINN approximated solution vs the reference solution. Right: $l^2$ relative error varies with the testing epoch.} 
	\label{PINNTest}         
\end{figure}

We employ the classical PINN method with the MscaleDNN framework(see Fig. \ref{fig2mscalednn}) to solve \eqref{eq:multiscale}, called this method as MPINN. The scale factors $\Lambda$ for MscaleDNN is set as $(1,2,3,4,5,10,\cdots, 90, 95, 100)$ and the size of each subnetwork is chosen as $(1, 30, 40, 30, 30, 30, 1)$. The activation function of the first hidden layer for all subnetworks is set as Fourier feature mapping(see Section \ref{opt2activation}) and the other activation functions(except for their output layer) are chosen as $\frac{1}{2}\sin(x)+\frac{1}{2}\cos(x)$ \cite{chen2022adaptive}, their output layers are all linear. For $\varepsilon=\frac{1}{32}, \frac{1}{64}$ and $\frac{1}{128}$,We train the aforementioned MPINN model for 50000 epochs and conduct testing every 1000 epochs within the training cycle. The optimizer is set as Adam with an initial learning rate of 0.01 and the learning rate will decay by $2.5\%$ for every 100 epochs. Finally, the results are demonstrated in Fig. \ref{PINNTest}. 

As $\varepsilon=1/32$, the coefficient $A^{\varepsilon}(x)$ possesses a little multi-scale information, the MPINN performs quite well. However, the permeability $A^{\varepsilon}(x)$ will exhibits various multi-scale properties for $\varepsilon=1/64$, the performance of MPINN deteriorates with a low relative error and the MPINN fails to converge for $\varepsilon=1/128$. In addition, we perform the MPINN with different setups of the hyperparameters such as the learning rate and the $\omega_{B}$ for $Loss_{B}$ in \eqref{loss2PINN} as well as the network size, but we still cannot obtain a satisfactory result.

\subsection{Unified architecture of FMPINN}

Based on the above observation, it is necessary to seek some extra techniques to improve the accuracy of the PINN. Inspired by the mixed finite element method \cite{ch03,araya2013multiscale} and the mixed residual method\cite{lyu2022mim}, we can leverage a mixed scheme to solve \eqref{eq:multiscale} by replacing the flux term $A^{\varepsilon}\nabla u$ in \eqref{eq:multiscale} with an auxiliary variable. This strategy not only can avoid the unfavorable effect of the oscillating coefficient $A^{\varepsilon}$, but also can reduce the computation burden of second-order derivatives in cost function when utilizing a multi-scale deep neural network to approximate the solution of \eqref{eq:multiscale}. Therefore, we introduce a flux variable $\bm{\phi}(\bm{x})=\big{(}\phi_1(\bm{x}), \ldots,\phi_d(\bm{x})\big{)}=A^{\varepsilon}(\bm{x})\nabla u^{\varepsilon}(\bm{x})$ and rewrite the first equation in \eqref{eq:multiscale} as the following expressions:
\begin{equation}\label{EqvilentForm2MS}
   \begin{aligned}
   	&-\mathbf{div}\bm{\phi}(\bm{x}) = f(\bm{x})\\
   	&\bm{\phi}(\bm{x}) - A^{\varepsilon}(\bm{x})\nabla u^{\varepsilon}(\bm{x}) = \bm{0}
   \end{aligned}
\end{equation}
Then we turn to search a couple of functions $(u^{\varepsilon},\bm{\phi})$ in admissible space, rather than approximating a unique solution of the original problem \eqref{eq:multiscale}. Here and thereafter, $(u^{\varepsilon},\bm{\phi})\in \mathcal{A}=\mathcal{H}^1(\Omega)\times \mathcal{H}(\mathbf{div};\Omega)$ with $\mathcal{H}^1(\Omega)=\big{\{} v\in L^2(\Omega):\nabla v\in L^2(\Omega)\big{\}}$ and $\mathcal{H}(\mathbf{div};\Omega)=\big{\{} \bm{\psi}\in [L^2(\Omega)]^d:\mathbf{div} \bm{\psi}\in L^2(\Omega)\big{\}}.$

When utilizing numerical solvers to address the equation \eqref{EqvilentForm2MS}, one can obtain the optimum solution by minimizing the following least-squares formula in the domain $\Omega$:
\begin{equation}\label{loss2continous}
     u^*,\bm{\phi}^* = \underset{( u,\bm{\phi})\in \mathcal{H}^1(\Omega)\times \mathcal{H}(\mathbf{div};\Omega)}{\arg\min}\mathcal{L}(u,\bm{\phi})
\end{equation}
with

\begin{equation}\label{intergalForm}
\mathcal{L}(u,\bm{\phi}) = \int_{\Omega}\big{|}-\mathbf{div}\bm{\phi}(\bm{x}) - f(\bm{x})\big{|}^2d\bm{x} + \beta\int_{\Omega}\big{|}\bm{\phi}(\bm{x}) - A^{\varepsilon}(\bm{x})\nabla u^{\varepsilon}(\bm{x})\big{|}^2d\bm{x}
\end{equation}
where $\beta>0$ is used to adjust the approximation error of the flux variable and flux term.

\begin{figure}[H]
	\centering
	\begin{tikzpicture}[scale=0.75]
	\node[] (input) at (2, -3.75) {Input layer};
		
	\node[circle, fill=green!60,inner sep=4.5pt] (x) at (2.0, 0) {$\bm{x}$};		
	
	\node[circle, fill=cyan!70,inner sep=4.5pt] (h10) at (4, 3) {};
	\node[circle, fill=cyan!70,inner sep=4.5pt] (h11) at (4, 2) {};
	\node[circle, fill=cyan!70,inner sep=4.5pt] (h12) at (4, 1) {};
	\node[circle, fill=cyan!70,inner sep=4.5pt] (h13) at (4, 0) {};
	\node[circle, fill=cyan!70,inner sep=4.5pt] (h14) at (4, -1) {};
	\node[circle, fill=cyan!70,inner sep=4.5pt] (h15) at (4, -2) {};
	\node[circle, fill=cyan!70,inner sep=4.5pt] (h16) at (4, -3) {};

	\draw[line width=1.0pt,->] (x) -- (h10);
	\draw[line width=1.0pt,->] (x) -- (h11);
	\draw[line width=1.0pt,->] (x) -- (h12);
	\draw[line width=1.0pt,->] (x) -- (h13);
	\draw[line width=1.0pt,->] (x) -- (h14);
	\draw[line width=1.0pt,->] (x) -- (h15);
	\draw[line width=1.0pt,->] (x) -- (h16);

	\node[circle, fill=cyan!70,inner sep=4.5pt] (h20) at (6.75, 2) {};
	\node[circle, fill=cyan!70,inner sep=4.5pt] (h21) at (6.75, 1) {};
	\node[circle, fill=cyan!70,inner sep=4.5pt] (h22) at (6.75, 0) {};
	\node[circle, fill=cyan!70,inner sep=4.5pt] (h23) at (6.75, -1) {};
	\node[circle, fill=cyan!70,inner sep=4.5pt] (h24) at (6.75, -2) {};
	\node[] (input) at (7, -3.75) {Hidden layers};
	
	\draw[line width=1.0pt,->] (h10) -- (h20);
	\draw[line width=1.0pt,->] (h10) -- (h21);
	\draw[line width=1.0pt,->] (h10) -- (h22);
	\draw[line width=1.0pt,->] (h10) -- (h23);
	\draw[line width=1.0pt,->] (h10) -- (h24);
	
	\draw[line width=1.0pt,->] (h11) -- (h20);
	\draw[line width=1.0pt,->] (h11) -- (h21);
	\draw[line width=1.0pt,->] (h11) -- (h22);
	\draw[line width=1.0pt,->] (h11) -- (h23);
	\draw[line width=1.0pt,->] (h11) -- (h24);
	
	\draw[line width=1.0pt,->] (h12) -- (h20);
	\draw[line width=1.0pt,->] (h12) -- (h21);
	\draw[line width=1.0pt,->] (h12) -- (h22);
	\draw[line width=1.0pt,->] (h12) -- (h23);
	\draw[line width=1.0pt,->] (h12) -- (h24);
	
	\draw[line width=1.0pt,->] (h13) -- (h20);
	\draw[line width=1.0pt,->] (h13) -- (h21);
	\draw[line width=1.0pt,->] (h13) -- (h22);
	\draw[line width=1.0pt,->] (h13) -- (h23);
	\draw[line width=1.0pt,->] (h13) -- (h24);;
	
	\draw[line width=1.0pt,->] (h14) -- (h20);
	\draw[line width=1.0pt,->] (h14) -- (h21);
	\draw[line width=1.0pt,->] (h14) -- (h22);
	\draw[line width=1.0pt,->] (h14) -- (h23);
	\draw[line width=1.0pt,->] (h14) -- (h24);
	
	\draw[line width=1.0pt,->] (h15) -- (h20);
	\draw[line width=1.0pt,->] (h15) -- (h21);
	\draw[line width=1.0pt,->] (h15) -- (h22);
	\draw[line width=1.0pt,->] (h15) -- (h23);
	\draw[line width=1.0pt,->] (h15) -- (h24);
	
	\draw[line width=1.0pt,->] (h16) -- (h20);
	\draw[line width=1.0pt,->] (h16) -- (h21);
	\draw[line width=1.0pt,->] (h16) -- (h22);
	\draw[line width=1.0pt,->] (h16) -- (h23);
	\draw[line width=1.0pt,->] (h16) -- (h24);
	
	\node[circle, fill=cyan!70,inner sep=4.5pt] (h30) at (9.5, 3) {};
	\node[circle, fill=cyan!70,inner sep=4.5pt] (h31) at (9.5, 2) {};
	\node[circle, fill=cyan!70,inner sep=4.5pt] (h32) at (9.5, 1) {};
	\node[circle, fill=cyan!70,inner sep=4.5pt] (h33) at (9.5, 0) {};
	\node[circle, fill=cyan!70,inner sep=4.5pt] (h34) at (9.5, -1) {};
	\node[circle, fill=cyan!70,inner sep=4.5pt] (h35) at (9.5, -2) {};
	\node[circle, fill=cyan!70,inner sep=4.5pt] (h36) at (9.5, -3) {};	
	
	\draw[line width=1.0pt,->] (h20) -- (h30);
	\draw[line width=1.0pt,->] (h20) -- (h31);
	\draw[line width=1.0pt,->] (h20) -- (h32);
	\draw[line width=1.0pt,->] (h20) -- (h33);
	\draw[line width=1.0pt,->] (h20) -- (h34);
	\draw[line width=1.0pt,->] (h20) -- (h35);
	\draw[line width=1.0pt,->] (h20) -- (h36);
	
	\draw[line width=1.0pt,->] (h21) -- (h30);
	\draw[line width=1.0pt,->] (h21) -- (h31);
	\draw[line width=1.0pt,->] (h21) -- (h32);
	\draw[line width=1.0pt,->] (h21) -- (h33);
	\draw[line width=1.0pt,->] (h21) -- (h34);
	\draw[line width=1.0pt,->] (h21) -- (h35);
	\draw[line width=1.0pt,->] (h21) -- (h36);
	
	\draw[line width=1.0pt,->] (h22) -- (h30);
	\draw[line width=1.0pt,->] (h22) -- (h31);
	\draw[line width=1.0pt,->] (h22) -- (h32);
	\draw[line width=1.0pt,->] (h22) -- (h33);
	\draw[line width=1.0pt,->] (h22) -- (h34);
	\draw[line width=1.0pt,->] (h22) -- (h35);
	\draw[line width=1.0pt,->] (h22) -- (h36);
	
	\draw[line width=1.0pt,->] (h23) -- (h30);
	\draw[line width=1.0pt,->] (h23) -- (h31);
	\draw[line width=1.0pt,->] (h23) -- (h32);
	\draw[line width=1.0pt,->] (h23) -- (h33);
	\draw[line width=1.0pt,->] (h23) -- (h34);
	\draw[line width=1.0pt,->] (h23) -- (h35);
	\draw[line width=1.0pt,->] (h23) -- (h36);
	
	\draw[line width=1.0pt,->] (h24) -- (h30);
	\draw[line width=1.0pt,->] (h24) -- (h31);
	\draw[line width=1.0pt,->] (h24) -- (h32);
	\draw[line width=1.0pt,->] (h24) -- (h33);
	\draw[line width=1.0pt,->] (h24) -- (h34);
	\draw[line width=1.0pt,->] (h24) -- (h35);
	\draw[line width=1.0pt,->] (h24) -- (h36);
	
	\node[circle, fill=yellow!100,inner sep=4.5pt] (o0) at (12, 1) {};	
	\node[circle, fill=yellow!100,inner sep=4.5pt] (o1) at (12, 0) {};	
	\node[circle, fill=yellow!100,inner sep=4.5pt] (o2) at (12, -1) {};	
	
	\draw[line width=1.0pt,->] (h30) -- (o0);
	\draw[line width=1.0pt,->] (h30) -- (o1);
	\draw[line width=1.0pt,->] (h30) -- (o2);
	
	\draw[line width=1.0pt,->] (h31) -- (o0);
	\draw[line width=1.0pt,->] (h31) -- (o1);
	\draw[line width=1.0pt,->] (h31) -- (o2);
	
	\draw[line width=1.0pt,->] (h32) -- (o0);
	\draw[line width=1.0pt,->] (h32) -- (o1);
	\draw[line width=1.0pt,->] (h32) -- (o2);
	
	\draw[line width=1.0pt,->] (h33) -- (o0);
	\draw[line width=1.0pt,->] (h33) -- (o1);
	\draw[line width=1.0pt,->] (h33) -- (o2);
	
	\draw[line width=1.0pt,->] (h34) -- (o0);
	\draw[line width=1.0pt,->] (h34) -- (o1);
	\draw[line width=1.0pt,->] (h34) -- (o2);
	
	\draw[line width=1.0pt,->] (h35) -- (o0);
	\draw[line width=1.0pt,->] (h35) -- (o1);
	\draw[line width=1.0pt,->] (h35) -- (o2);
	
	\draw[line width=1.0pt,->] (h36) -- (o0);
	\draw[line width=1.0pt,->] (h36) -- (o1);
	\draw[line width=1.0pt,->] (h36) -- (o2);
	
	\node[circle] (y0) at (14, 1) {};
	\node[circle] (y1) at (14, 0) {};
	\node[circle] (y2) at (14, -1) {};
	\node[] (input) at (12, -3.75) {Out layer};
	
	\draw[line width=1.0pt,->] (o0) -- node[above]{$u$}(y0);
	\draw[line width=1.0pt,->] (o1) -- node[above]{$\phi_1$}(y1);
	\draw[line width=1.0pt,->] (o2) -- node[above]{$\phi_2$}(y2);	
	\end{tikzpicture}
	\caption{ \small The multi-output neural network for approximating the state and flux variables}
	\label{fig2multioutput}
\end{figure}
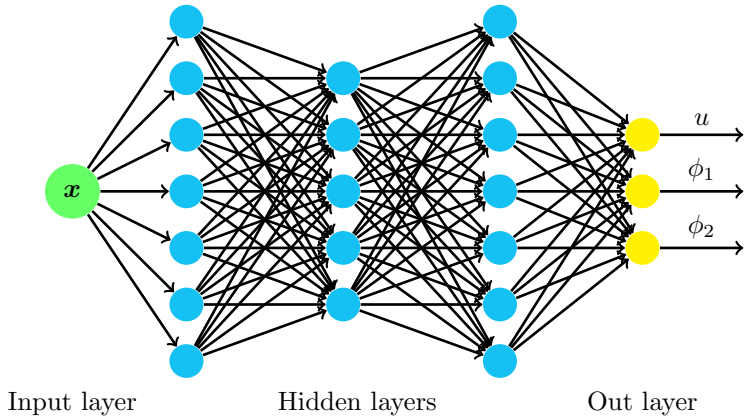

Generally, two independent neural networks are necessary to approximate the
flux variable  $\bm{\phi}$ and solution $u$, but $\bm{\phi}$ is unconstrained without any coercive boundary condition. Based on the potentiality of DNN for approximating any linear and non-linear complex functions, we take a DNN with multi outputs to model ansatzes $\bm{\phi}$ and $u$, denoted by $\bm{\phi}_{NN}$ and $u_{NN}$, respectively. 
Fig. \ref{fig2multioutput} describes the  multi-output neural network for input $\bm{x}\in \mathbb{R}^2$.

Once the expressions of auxiliary functions $\bm{\phi}$ and solution $u$ have been determined, we can discretize \eqref{intergalForm} by the Monte Carlo method \cite{robert1999monte}, then employ the PINN conception and obtain the following form
\begin{equation}\label{loss2ellptic}
    \mathscr{L}_{in}(S_I;\bm{\theta})=\frac{|\Omega|}{N_{in}}\sum_{i=1}^{N_{in}}\bigg{[}\big{|}-\textbf{div} \bm{\phi}_{NN}(\bm{x}_I^i;\bm{\theta}) - f(\bm{x}_I^i)\big{|}^2 + \beta\big{|} \bm{\phi}_{NN}(\bm{x}^i_I,\bm{\theta}) -A^{\varepsilon}(\bm{x}^{i}_I)\nabla u_{NN}(\bm{x}_I^i,\bm{\theta})\big{|}^2\bigg{]},
\end{equation}
for $\bm{x}_I^i \in S_I$, here and hereinafter $S_I$ stands for the collection sampled from $\Omega$ with prescribed probability density.

Same to the traditional numerical methods such as FDM and FEM for addressing PDEs, boundary conditions play a crucial role in DNN representation as well. They serve as important constraints that ensure the uniqueness and accuracy of the solution. Consequently, the output $u_{NN}$ of DNN should also satisfy the boundary conditions of \eqref{eq:multiscale}, which means 
\begin{equation}\label{loss2boundary}
\mathscr{L}_{bd}(S_B;\bm{\theta})=\frac{1}{N_{bd}}\sum_{j=1}^{N_{bd}} \bigg{[}\mathcal{B}u_{NN}\big(\bm{x}_B^j;\bm{\theta}\big)-g(\bm{x}_B^j)\bigg{]}^2\rightarrow 0   ~~\text{for}~~\bm{x}_B^j\in S_B.
\end{equation}
here and hereinafter $S_B$ represents the collection sampled on $\partial \Omega$ with prescribed probability density. 

According to the above results, the weights and biases of the DNN model are updated by optimizing gradually the following cost function:
\begin{equation}\label{totalloss}
\mathscr{L}({S_I,S_B};\bm{\theta}) = \mathscr{L}_{in}(S_I;\bm{\theta}) + \gamma\mathscr{L}_{bd}(S_B;\bm{\theta})
\end{equation}
where $S_I=\{\bm{x}_I^i\}_{i=1}^{N_{in}}$ and $S_B=\{\bm{x}_B^j\}_{j=1}^{N_{bd}}$ stand for the train data of $\Omega$ and $\partial \Omega$, respectively. The term of $\mathscr{L}_{in}$ composed of the residual governed by differential equations and the discrepancy with respect to flux, minimizes
the residual of the PDE, whereas the term of $\mathscr{L}_{bd}$ pushes the DNN solver to match the given boundary conditions. In addition, a constant parameter $\gamma>0$ is introduced to forces well the $\mathscr{L}_{bd}(S_B;\bm{\theta})\rightarrow 0$ in the loss function, it is increasing gradually with training process going on. 

Based on the analysis in \cite{bersetche2023deep}, a nonconstant continuous activation function $\sigma$ can guarantee the mapping $\bm{\theta}\mapsto (u_{NN},\bm{\phi}_{NN})$ is continuous, then the distance between approximation functions $\bm{q}_{NN}=(u_{NN},\bm{\phi}_{NN})$ and exact solution $\bm{q}^{*}=(u^*,\bm{\phi}^*)$ will decrease by adjusting gradually the parameters of DNN, i.e.,
\begin{equation*}
d(\bm{q}^*,\mathcal{A}_k)=\underset{\bm{q}_{NN}\in\mathcal{A}_k}{\inf}\|\bm{q}^*-\bm{q}_{NN}\| \rightarrow 0 ~~\text{as}~~ k\rightarrow 
  \infty.
\end{equation*}
It means the loss function $\mathscr{L}({S_I,S_B};\bm{\theta})$ will attain the corresponding minimum when $d\rightarrow0$.

Hence,  Our purpose is to find an optimal set of parameter $\bm{\theta}^*$ such that the approximations $u_{NN}$ and $\bm{\phi}_{NN}$ minimize the loss function $\mathscr{L}({S_I,S_B};\bm{\theta})$. 
In order to obtain the ideal  $\bm{\theta}^*$, one can update the weights and biases of DNN through the optimization methods such as gradient descent (GD) or stochastic gradient descent (SGD) during the training process. In this context, the SGD method with a ”mini-batch” of training data is given by:
\begin{equation}\label{optimize}
\bm{\theta}^{k+1}=\bm{\theta}^{k}-\alpha^k\nabla_{\bm{\theta}^k}\mathscr{L}(\bm{x};\bm{\theta}^{k})~~\text{with}~~\bm{x}\in S_I~\text{or}~\bm{x}\in S_B
\end{equation}
where the ``learning rate'' $\alpha^k$ decreases with $k$ increasing.

\subsection{Option of activation function for FMPINN and its explanation}\label{opt2activation}
Choosing a suitable and effective activation function is a critical concern when aiming to enhance the performance of DNN in computer vision, natural language processing, and scientific computation. Generally, an activation function such as rectified linear unit $\text{ReLU}(\bm{z})$ and hyperbolic tangent function $\tanh(\bm{z})$, can obviously improve the capacity and nonlinearity of neural networkS to address various nonlinear problems, such as the solution of various PDEs and classification. Recently, the works \cite{Xu_2020,rahaman2018spectral} manifested that the DNN often captures firstly the low-frequency component for target functions, then match the high-frequency component, they called it as the spectral bias or frequency preference of DNN. Under this phenomenon, many researchers attempt to utilize a Fourier feature mapping consisting of sine and cosine as the activation function to  improve the capacity of MscaleDNN, it will mitigate the pathology of spectral bias and enable networks to learn high frequencies more effectively\cite{rahaman2018spectral,wang2020eigenvector, tancik2020fourier,li2023deep}. It is expressed as follows:
\begin{equation}
\zeta(\bm{x}) = 
\left[\begin{array}{c}
\cos(\bm{\kappa} \bm{x})\\
\sin(\bm{\kappa} \bm{x}) 
\end{array}
\right],
\label{fourier}
\end{equation}
where $\bm{\kappa}$ is a user-specified vector or matrix (trainable or untrainable) which is consistent with the number of neural units in the first hidden layer for DNNs. Further, the work \cite{li2023subspace} designed a soften Fourier mapping by introducing a relaxing parameter $s\in(0,1]$ in $\zeta(\bm{x})$,  numerical results show that this modification will improve the performance of $\zeta(\bm{x})$. Actually, this activation function is used in the first hidden layer of DNN and maps the input data in $\Omega$ into a range of $[-1,1]$, then enhances the ability of DNN and expedites its convergence. 

Therefore, a real function  $\mathcal{P}(x)$ represented by DNN can be expressed as follows
\begin{equation*}
\mathcal{P}(\bm{x}) =  \sum_{n=0}^{\tilde{N}}\bigg{(}S\left(\cos(\bm{k}_n \bm{x});\bm{\bar{\theta}}_n\right)+T\left(\sin(\bm{k}_n\bm{x});\bm{\tilde{\theta}}_n\right)\bigg{)},
\end{equation*}
where $S(\cdot,\bm{\bar{\theta}}), T(\cdot,\bm{\tilde{\theta}})$ are the DNNs or the sub-modules of DNNs, respectively, $\{\bm{k}_0,\bm{k}_1,\bm{k}_2,\cdots\}$ are the frequencies of interest for the objective function. Obviously, the first hidden layer performed by Fourier feature mapping mimics the Fourier basis function, and the remaining blocks with different activation functions are used to learn the coefficients of these functions. After performing the Fourier mapping for input points with a given scale factor, the neural network can well capture the fine varying information for multi-scale problems. 

\begin{remark}
    (\textbf{Lipschitz continuous})~If an activation function $\sigma$ is continuous(i.e., $\sigma\in C^1$) and satisfies the following boundedness condition:
    \begin{equation*}
        |\sigma(x)|<1~~~~\textup{and}~~~~|\sigma'(x)|<1
    \end{equation*}
    for any $x\in\mathbb{R}$. Then, we have
    \begin{equation*}
        |\sigma(x)-\sigma(y)|\leqslant |x-y|~~~~\textup{and}~~~~|\sigma'(x)-\sigma'(y)|<|x-y|
    \end{equation*}
    for any $x,y\in\mathbb{R}$. Obviously, the activation functions $\tanh(x)$, sigmoid(x), Fourier feature mapping $\zeta(a)$ and $\frac{1}{2}sin(x)+\frac{1}{2}cos(x)$ are all satisfy the above condition and have a good regularity, they will overcome the gradient explosion of parameter in the backpropagation for DNN and improve the capacity of DNN.
\end{remark}

\subsection{Simple error analysis for FMPINN}



In recent times, there have been endeavors to rigorously analyze the convergence rate of the deep mixed residual method and compare it with the deep Galerkin method (DGM) and deep Ritz method (DRM) across different scenarios\cite{bersetche2023deep,gu2023error,li2022priori}. In this study, we investigate those results of convergence again, then provide the expression of generalization error for FMPINN and some remarks of errors.

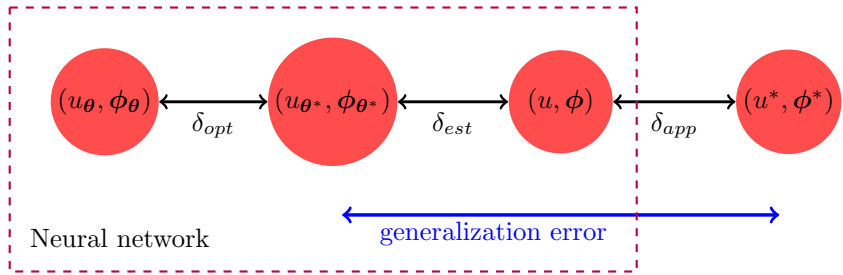
\begin{figure}[H]
	\centering
	\begin{tikzpicture}
	    	\node[circle, fill=red!70,inner sep=1.5pt] (n1) at (0, 3) {$(u_{\bm{\theta}},\bm{\phi}_{\bm{\theta}})$};
	    	
	    	\node[circle, fill=red!70,inner sep=1.5pt] (n2) at (3, 3) {$(u_{\bm{\theta}^*},\bm{\phi}_{\bm{\theta}^*})$};
	    	\node[circle, fill=red!70,inner sep=4.5pt] (n3) at (6, 3) {$(u,\bm{\phi})$};
	    	\node[circle, fill=red!70,inner sep=1.5pt] (n4) at (9, 3) {$(u^*,\bm{\phi}^*)$};
	    	\draw[line width=1.1pt,color=black,<->] (n1)--node[below]{$\delta_{opt}$}(n2);
	    	\draw[line width=1.1pt,color=black,<->] (n2)--node[below]{$\delta_{est}$}(n3);
	    	\draw[line width=1.1pt,color=black,<->] (n3)--node[below]{$\delta_{app}$}(n4);
	    	
	    	\node[] (n5) at (3, 1.5) {};
	    	\node[] (n6) at (9, 1.5) {};
	    	\draw[line width=1.3pt,color=blue,<->] (n5)--node[left,xshift=0.75cm, yshift=-0.25cm]{generalization error}(n6);
	    	
	    	\draw[thick, dashed, draw = purple](-1.25, 4.25)-- (7.0, 4.25) -- (7.0, 0.75) -- (-1.25, 0.75) -- (-1.25, 4.25);
	    	
	    	\node[rectangle,inner sep=1.5pt] (n2) at (0.2, 1.2) {Neural network};	
	\end{tikzpicture}
	\caption{Illustration of the total error for FMPINN.}
	\label{fig:total error}
\end{figure}

For convenience, let $\bm{q}^*=(u^*,\bm{\phi}^*)$ be the exact solution of equation \eqref{EqvilentForm2MS} or the minimum of  cost function \eqref{loss2continous} with \eqref{intergalForm} for coercive boundary constraints. Meantime, the $\bm{q}_{\bm{\theta}^*}=(u_{\bm{\theta}^*},\bm{\phi}_{\bm{\theta}^*})$ stands for the final output of DNN optimized by SGD optimizer(such as Adam or LBFGS) that attains the local minimum of \eqref{totalloss}. Further, we let $\widetilde{\mathscr{L}}(u,\bm{\phi})$ be the cost function evaluated on $N$ points sampled from $\Omega$ and denote the output of DNN as $\bm{q}_{\bm{\theta}}=(u_{\bm{\theta}},\bm{\phi}_{\bm{\theta}})$. Finally, $
\mathcal{S}_{NN}$ represents the function space sapnned by the output of DNN.
Then, the total error(or generalization error) between the exact solution $\bm{q}^*$ and the output of DNN $\bm{q}_{\bm{\theta}}$ can be expressed as
\begin{equation}\label{boundary2error}
    \big{\|}u_{\bm{\theta}}-u^*\big{\|}_{\mathcal{H}^{1}(\Omega)}+\big{\|}\bm{\phi}_{\bm{\theta}}-\bm{\phi}^*\big{\|}_{\mathcal{H}^{1}(\bm{div},\Omega)}\leqslant C(coe)\sqrt{\delta_{app}+\delta_{est}+\delta_{opt}}
\end{equation}
with
\begin{equation*}
    \begin{cases}
        \delta_{app}= \underset{(u,\bm{\phi})\in\mathcal{S}_{NN}}{\inf}\|u-u^*\|^2_{\mathcal{H}^1(\Omega)}+\|\bm{\phi}-\bm{\phi}^*\|^2_{\mathcal{H}(\bm{div},\Omega)}\\
        \delta_{est} = \underset{(u,\bm{\phi})\in\mathcal{S}_{NN}}{\sup}[\mathcal{L}(u,\bm{\phi})-\widetilde{\mathscr{L}}(u,\bm{\phi})] + \underset{(u,\bm{\phi})\in\mathcal{S}_{NN}}{\sup}[\widetilde{\mathscr{L}}(u,\bm{\phi})-\mathcal{L}(u,\bm{\phi})]\\
        \delta_{opt} = \widetilde{\mathscr{L}}(u_{\bm{\theta}^*},\bm{\phi}_{\bm{\theta}^*}) - \widetilde{\mathscr{L}}(u_{\bm{\theta}},\bm{\phi}_{\bm{\theta}})
    \end{cases}
\end{equation*}
In which, the approximated error $\delta_{app}$ indicates the difference between $(u^*,\bm{\phi}^*)$ and its projection onto $\mathcal{S}_{NN}$, the estimation error $\delta_{est}$ measures the difference between the continuous cost function $\mathcal{L}$ and discrete cost function $\widetilde{\mathscr{L}}$, the optimization error $\delta_{opt}$ stands for the discrepancy between the output of DNN with optimizing and the output of DNN without optimizing. In Fig. \ref{fig:total error}, we depict the diagram of error for FMPINN. 

\begin{remark}\label{error2app}
For the approximating error, it is generally dependent on the architectural design of the neural network and the choice of the activation function. Classical radial basis network\cite{orr1996introduction}, the vanilla DNN and extreme learning machine(ELM)\cite{ding2014extreme} are the common meshless method for approximating the solution of PDEs. To address the spatio-temporal problems, some hybrid network frameworks  have been designed by combining PINN with traditional numerical methods to solve PDE, such as FDM-PINN and Runge-Kutta PINN\cite{raissi2019physics,xiang2022hybrid}. Moreover, instead of soft constraints by a hard manner for the boundary or initial conditions in those methods, the approximation will automatically meet the boundary and initial conditions of PDEs, then reduce the complexity and improve the precision of NN\cite{sun2020surrogate}. On the other hand, a powerful activation function, such as the hyperbolic tangent activation function and Fourier feature mapping, not only enhance the nonlinearity of DNN, but also improve its approximating capacity and accuracy. In addition, some available data are generally considered as a loss term to reduce the approximating error.
\end{remark}

\begin{remark}\label{error2sta}
Generally, the proposed FMPINN surrogate can provide more accurate approximations as the number of random collocation points increases. However, it will lead to heavy computational costs for lots of samplings. Then, it is worthwhile to take into account the trade-off between accuracy and computational cost  when designing a DNN surrogate and determining its training mode. Alternatively, one can employ some effective low-discrepancy sampling approaches to decrease the statistical error, such as the Latin hypercube sampling  method\cite{viana2016tutorial},
quasi-random sampling\cite{shaw1988quasirandom} and multilevel Monte Carlo method\cite{giles2015multilevel}. 

\end{remark}

\begin{remark}\label{error2opt}


Since the cost function generally is non-convex and has several local minima, then the gradient-based optimizer will almost certainly become caught in one of them. Therefore, choosing a good optimizer is important to reduce the optimization error and get a better minimum. In many scenarios of optimizing DNN, the Adam optimization method has shown its good performance including efficiency and accuracy, it can dynamically adjust the learning rates of each parameter by using the first and second moments estimation of the gradients\cite{kingma2015adam}. BFGS is a quasi-Newton method and numerically stable, it may provide a higher-precision approximated solution\cite{yuan2011active}. In an implementation, the limited memory version of BFGS(L-BFGS) is the common choice to decrease the optimization error and accelerate convergence for cases with a little amount of training data and/or residual points. Further, by combining the merits  of the above two approaches, one can optimize the cost function firstly by the Adam algorithm with a predefined stop criterion, then obtain a better result by the L-BFGS optimizer.


\end{remark}

\section{FMPINN algorithm}\label{algori2FMPINN}
For the FMPINN method with the MscaleDNN model composed of $Q$ subnetworks as in Fig.\ref{fig2mscalednn} being its solver, the input data for each subnetwork will be transformed by the following operation 
\begin{equation*}
    \widehat{\bm{x}}= a_i * \bm{x}, \quad i=1,2, \ldots, Q
\end{equation*}
with $a_i\geqslant 1$ being a positive scalar factor, it means the scale vector $\Lambda=(a_1, a_2,\ldots, a_Q)$ as in \eqref{separate}. Denoting the output of each subnetwork as $\boldsymbol{F}_{i}(i=1,2,\ldots,Q)$, then the overall output of the MscaleDNN model is obtained by 
\begin{equation*}\label{eq:initial solution}
    \boldsymbol{y}(\bm{x};\bm{\theta}) =\frac{1}{Q}\sum_{i=1}^{Q}\frac{\boldsymbol{F}_{i}}{a_i}.
\end{equation*}

According to the above discussions, the procedure of the FMPINN algorithm for addressing the multi-scale problem \eqref{eq:multiscale} in finite-dimensional spaces is described in the following.

\begin{algorithm}[H]
\caption{FMPINN algorithm for solving multi-scale PDEs\eqref{eq:multiscale}}
1. Generating the $k_{th}$ training set $\mathcal{S}^k$ includes interior points $S_{I}^k=\{\bm{x}^i_I\}_{i=1}^{N_{in}}$ with $\bm{x}_I^i\in\mathbb{R}^d$ and boundary points $S_{B}^k=\{\bm{x}^j_B\}_{j=1}^{N_{bd}}$ with $\bm{x}_B^j\in\mathbb{R}^d$. Here, we draw the random points $\bm{x}_I^i $ and $\bm{x}_B^j$ from $\mathbb{R}^d$ with positive probability density $\nu$, such as uniform distribution.
	
2. Calculating the objective function $\mathscr{L}(\mathcal{S}^k;\bm{\theta}^{k})$ for train set $\mathcal{S}^k$:
\begin{equation*}
    \mathscr{L}(\mathcal{S}^k;\bm{\theta}^{k}) = \mathscr{L}_{in}(S_{I}^k;\bm{\theta}^k) + \gamma\mathscr{L}_{bd}(S_{B}^k;\bm{\theta}^k)
\end{equation*}
with $\mathscr{L}_{in}(\cdot;\bm{\theta}^k)$ being defined in \eqref{loss2ellptic} and $\mathscr{L}_{bd}(\cdot;\bm{\theta}^k)$ being defined in \eqref{loss2boundary}.
	
3. Take a descent step at the random point of $\bm{x}^k$:
\begin{equation*}
    \bm{\theta}^{k+1}=\bm{\theta}^{k}-\alpha^k\nabla_{\bm{\theta}^k}\mathscr{L}(\tilde{\bm{x}}^k;\bm{\theta}^{k})~~\text{with}~~\tilde{\bm{x}}^k\in\mathcal{S}^k,
\end{equation*}
where the ``learning rate'' $\alpha^k$ decreases with $k$ increasing.
	
4. Repeat steps 1-3 until the convergence criterion is satisfied or the objective function tends to be stable.
	\label{algor:FMPINN}
\end{algorithm}

\section{Numerical experiments}\label{sec:04}
The goal of our experiments is to show that our Fourier-based mixed physics-informed neural networks are indeed capable of approximating the analytical solution given in \eqref{eq:multiscale}. For comparison purposes, the PINN method with MscaleDNN being its solver and the local deep learning method(LDLM) with normal DNN being its solver are as the baseline to solve \eqref{eq:multiscale} in varying-dimensional spaces. 

\subsection{Model and training setup}\label{sec:model}

In the aforementioned FMPINN and MPINN models, a standard MscaleDNN with multi subnetworks that stretch the input data via various scale factors is configurated as their solver. The MscaleDNN consists of 25 subnetworks according to the manually defined frequencies vector $\Lambda=(1, 2, 3, 4, 5, 10,\cdots,90, 95, 100)$. Each subnetwork contains 5 hidden layers with proper size and the activation function of the first hidden layer for each subnetwork is set as Fourier feature mapping and the other activation functions(except for the output layer) are set as $\frac{1}{2}\sin(x)+\frac{1}{2}\cos(x)$, its output layer is linear. The overall output is a weighted sum of the outputs of all subnetworks through the relevant scale factors. In terms of the LDLM\cite{zhu2021local}, two activation functions are considered for this model: LDLM1 with $ReQU=\max\{0,x\}^2$ being its activation for hidden layers and LDLM2 with $\frac{1}{2}\sin(x)+\frac{1}{2}\cos(x)$ being its activation function for hidden layers, their output are all linear.






In our numerical experiments, all training data are sampled from the domain(including its boundaries) of interest in Euclidean space $\mathbb{R}^d$, the sampling probability densities are assigned as the uniform distribution. We train all neural networks by an Adam optimizer with an initial learning rate of 0.01, and the learning rate will be decayed by 2.5\% for every 100 training epochs~\cite{kingma2015adam}. Here, the following $l^2$ relative error is used to evaluate our models:
\begin{equation*}
REL = \sqrt{\frac{\sum_{i=1}^{N'}|\tilde{u}(\bm{x}^i)-u^*(\bm{x}^i)|^2}{\sum_{i=1}^{N'}|u^*(\bm{x}^i)|^2}}
\end{equation*}
where $\tilde{u}(\bm{x}^i)$ and $u^*(\bm{x}^i)$ are the approximate solution of deep neural network and exact solution for testing points $\{\bm{x}^i\}(i=1,2,\cdots,N')$, respectively, and $N'$ represents the number of sample points for testing. In order to visualize the training process, our model will be evaluated once for every 1000 iterations in the whole training cycle and recorded the result at the end. In our codes, the penalty parameter $\gamma$ is set as
\begin{equation}
\gamma=\left\{
\begin{aligned}
\gamma_0, \quad &\textup{if}~~i_{\textup{epoch}}<M_{\max}*0.1\\
10\gamma_0,\quad &\textup{if}~~M_{\max}*0.1<=i_{\textup{epoch}}<M_{\max}*0.2\\
50\gamma_0, \quad&\textup{if}~~ M_{\max}*0.2<=i_{\textup{epoch}}<M_{\max}*0.25\\
100\gamma_0, \quad&\textup{if}~~ M_{\max}*0.25<=i_{\textup{epoch}}<M_{\max}*0.5\\
200\gamma_0, \quad&\textup{if}~~ M_{\max}*0.5<=i_{\textup{epoch}}<M_{\max}*0.75\\
500\gamma_0, \quad&\textup{otherwise}
\end{aligned}
\right.
\end{equation}
where the $\gamma_0=10$ in all our tests and $M_{\max}$ represents the total number of epochs. We implement and perform all neural network models by means of the package of Pytorch (version 1.14.0) on a workstation (64-GB RAM, single NVIDIA GeForce RTX 4090 24-GB).


\subsection{Performance of FMPINN for solving multi-scale elliptic PDEs}
\begin{example}\label{DiffusionEq_1d_01}
	Firstly, we consider the one-dimensional case for \eqref{eq:multiscale} with Dirichlet boundary in interval $[0,1]$, in which $A^{\varepsilon}(x)$ is given by
	\begin{equation}\label{DiffusionEq_1d_01_aeps}
	A^{\varepsilon}(x)=\left(2+\cos\left(2\pi\frac{x}{\varepsilon}\right)\right)^{-1}
	\end{equation}
with a small parameter $\varepsilon>0$ such that $\varepsilon^{-1}\in\mathbb{N}^+$ and the force term $f(x)=1$. Under these conditions, a unique solution is given by
	\begin{equation}\label{DiffusionEq_1d_01_ueps}
	u^\varepsilon(x) = x-x^2+\varepsilon\left(\frac{1}{4\pi}\sin\left(2\pi\frac{x}{\varepsilon}\right)-\frac{1}{2\pi}x\sin\left(2\pi\frac{x}{\varepsilon}\right)-\frac{\varepsilon}{4\pi^2}\cos\left(2\pi\frac{x}{\varepsilon}\right)+\frac{\varepsilon}{4\pi^2}\right).
	\end{equation}
Clearly, the analytical solution induces its boundary condition $u(0)=u(1)=0$.
\end{example}

\begin{figure}
	\centering
	\subfigure[Rough coefficient]{
		\label{0p1Aeps}
		\includegraphics[scale=0.4]{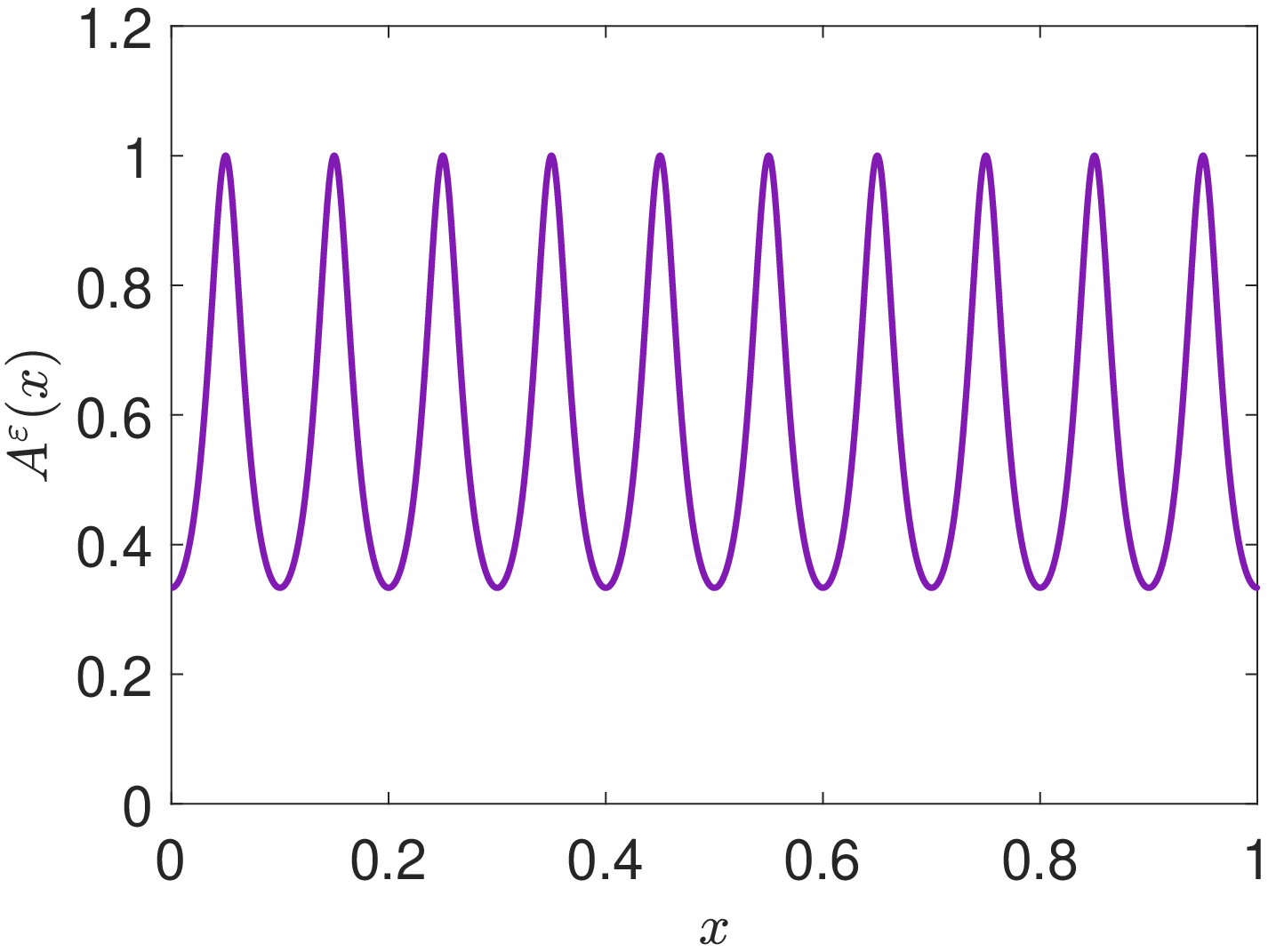}
	}
	\subfigure[loss of flux term for FMPINN]{
		\label{0p1loss2flux}
		\includegraphics[scale=0.4]{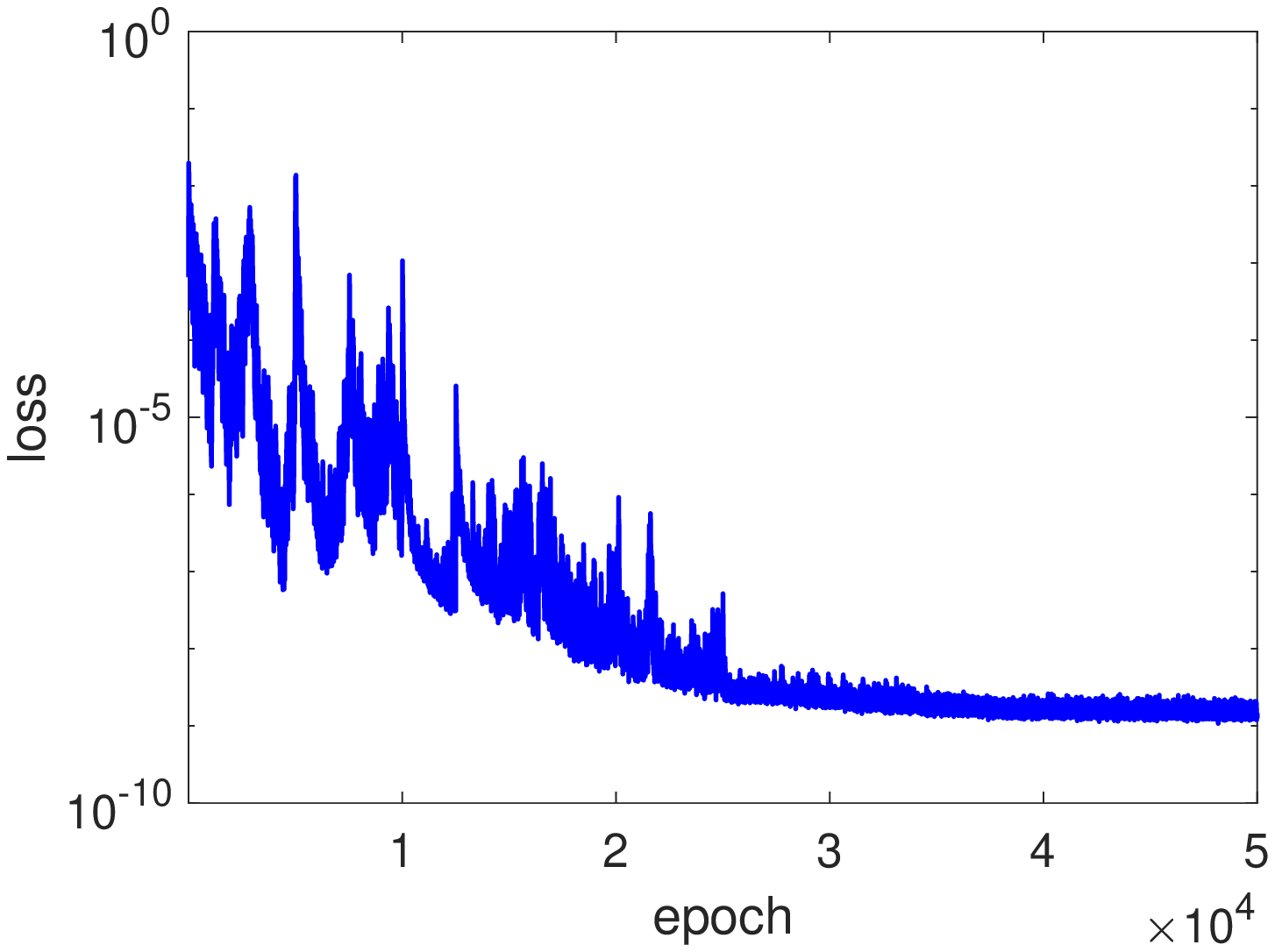}
	}
	\subfigure[Solutions]{
		\label{0p1Solus}
		\includegraphics[scale=0.4]{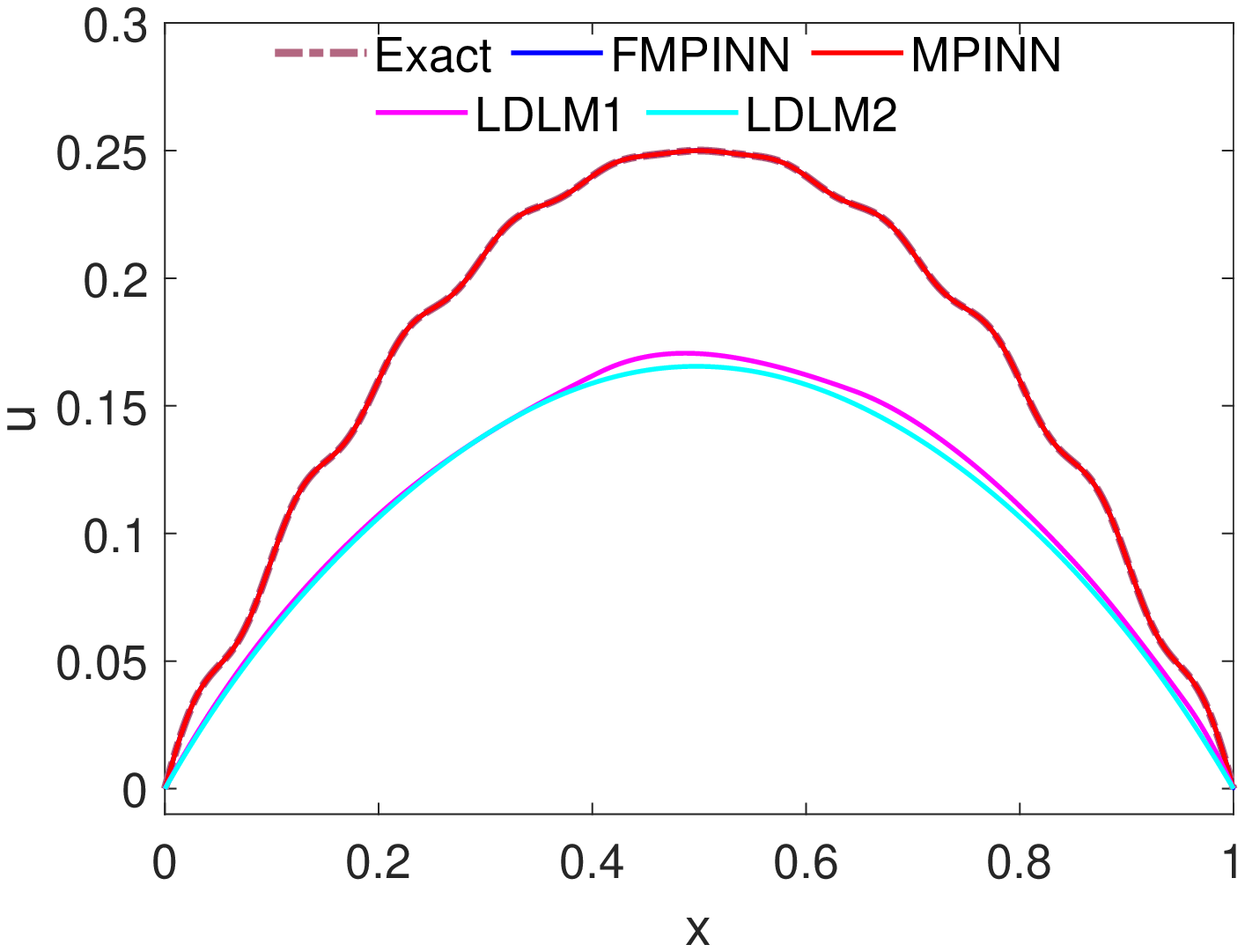}
	}
	\subfigure[point-wise error of FMPINN]{
		\label{0p1ABS2FMPINN}
		\includegraphics[scale=0.4]{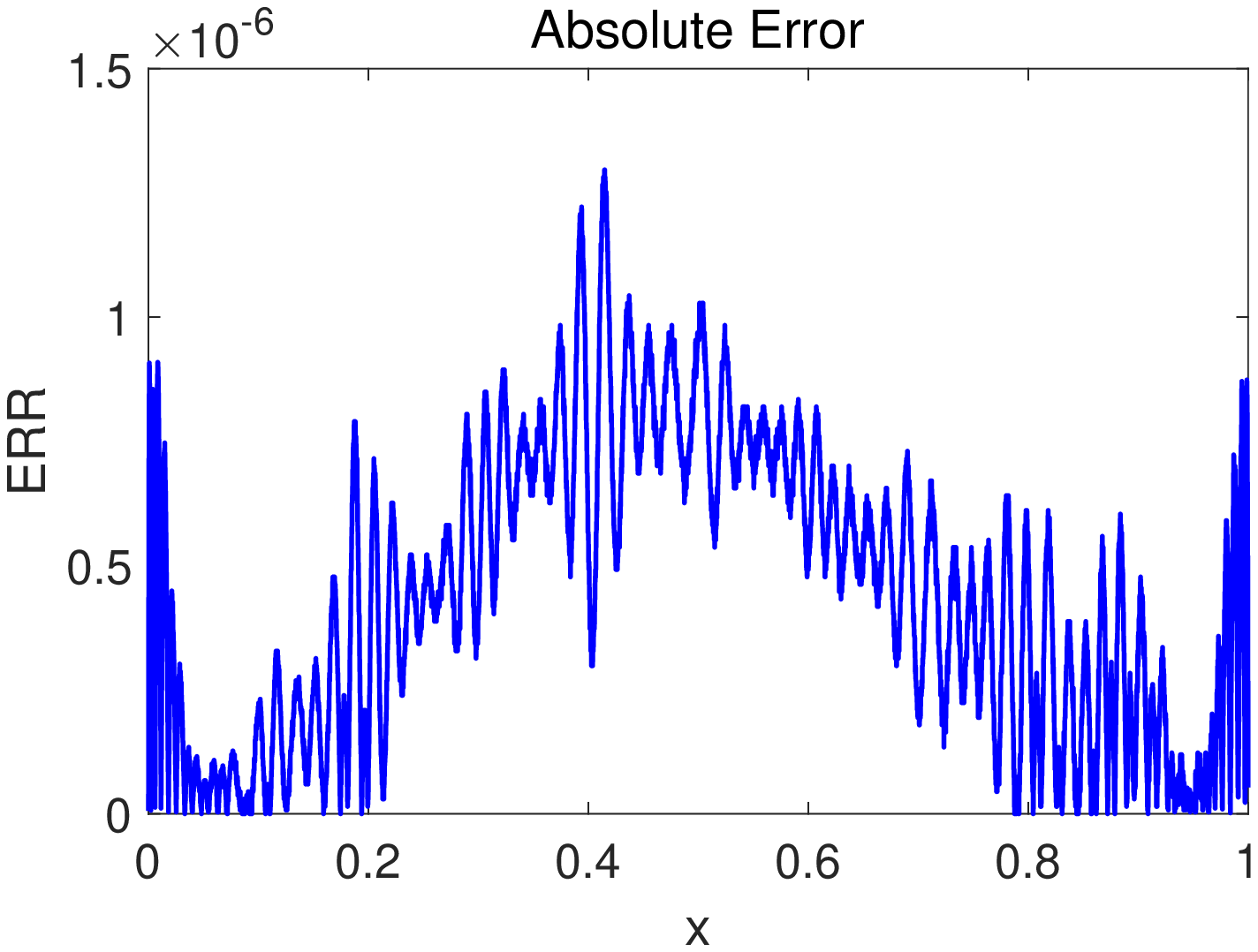}
	}
	\subfigure[point-wise error of MPINN]{
		\label{0p1ABS2MPINN}
		\includegraphics[scale=0.4]{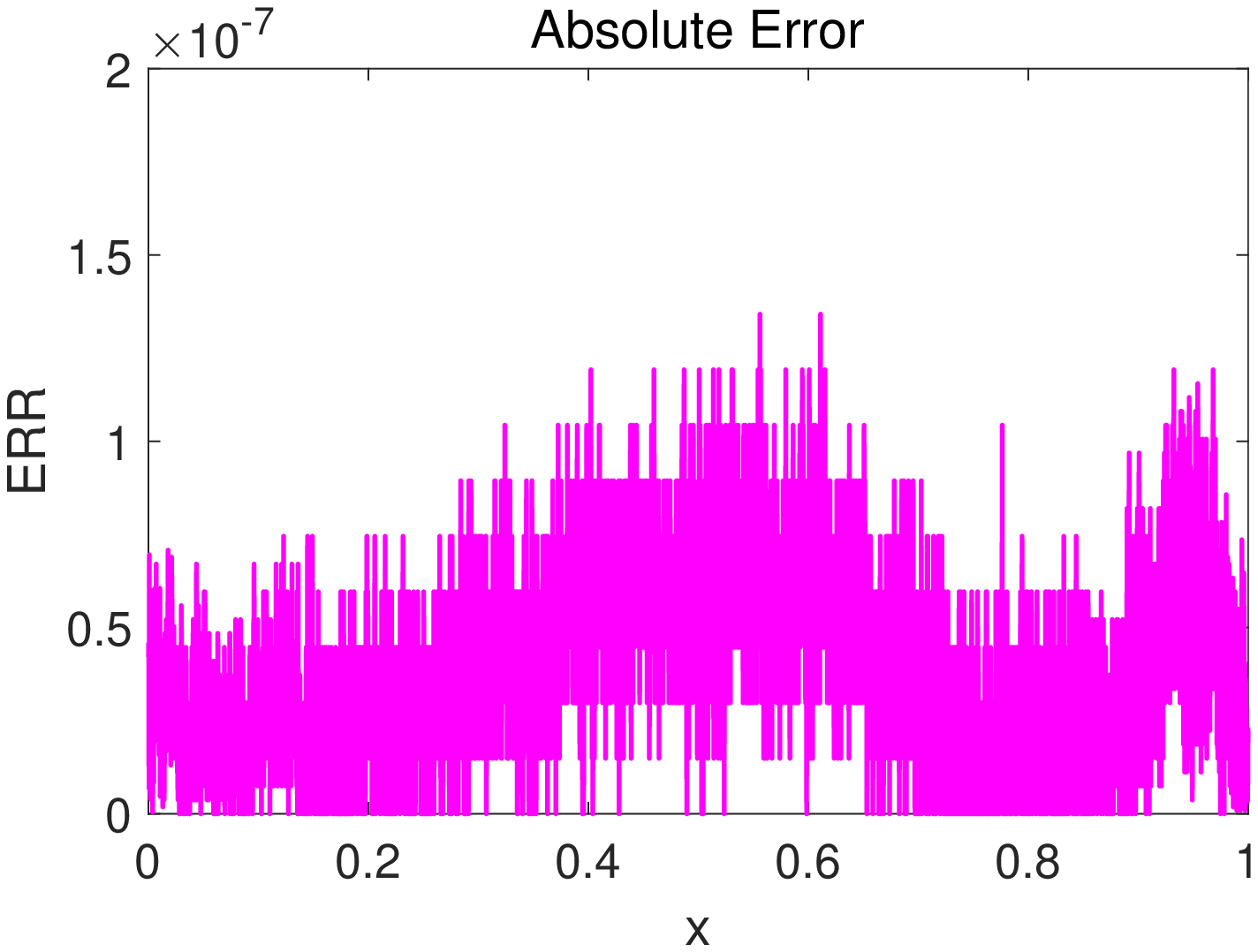}
	}
	\subfigure[point-wise error of LDLM1]{
		\label{0p1ABS2LDLM1}
		\includegraphics[scale=0.4]{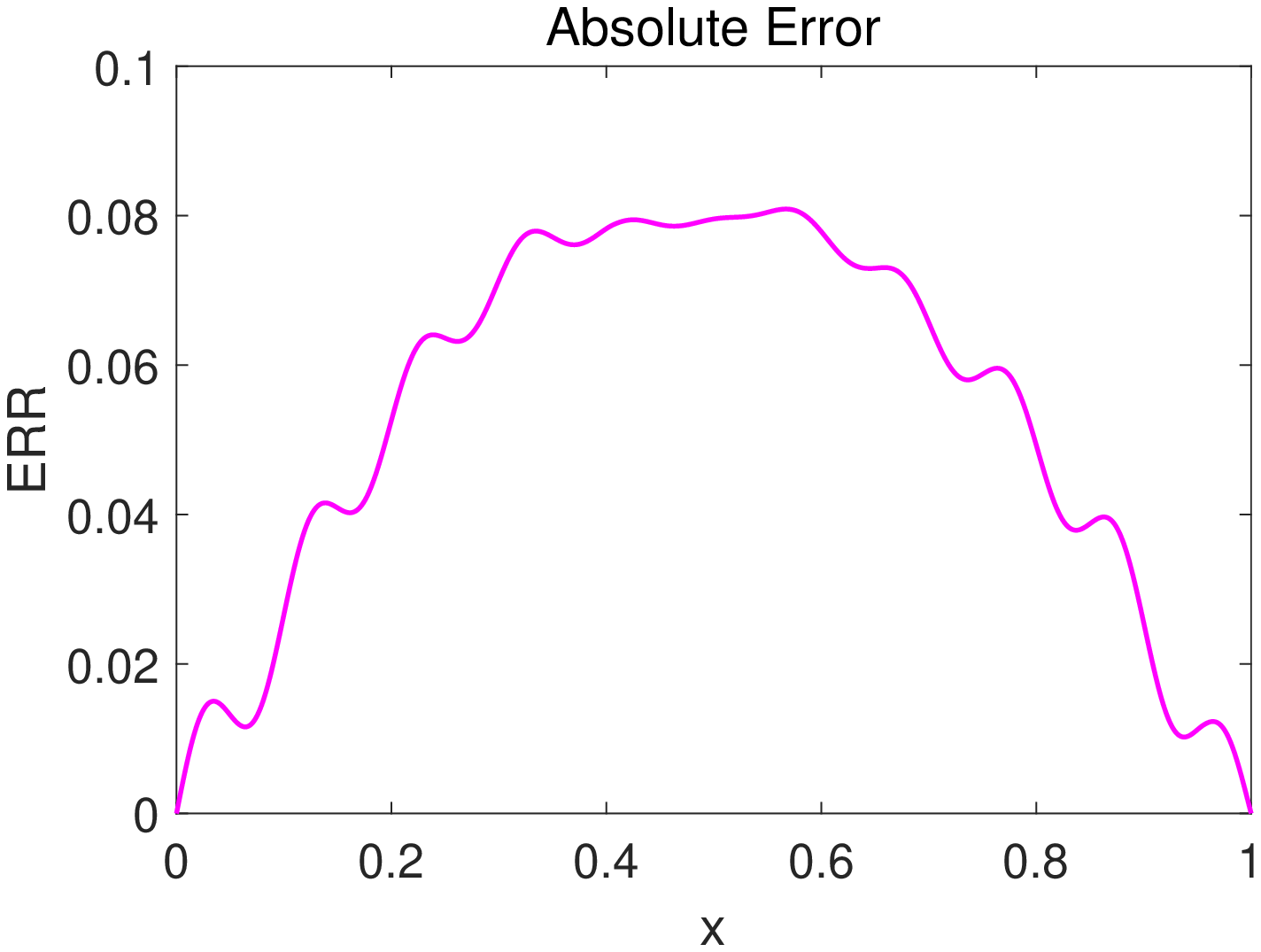}
	}
	\subfigure[point-wise error of LDLM2]{
		\label{0p1ABS2LDLM2}
		\includegraphics[scale=0.4]{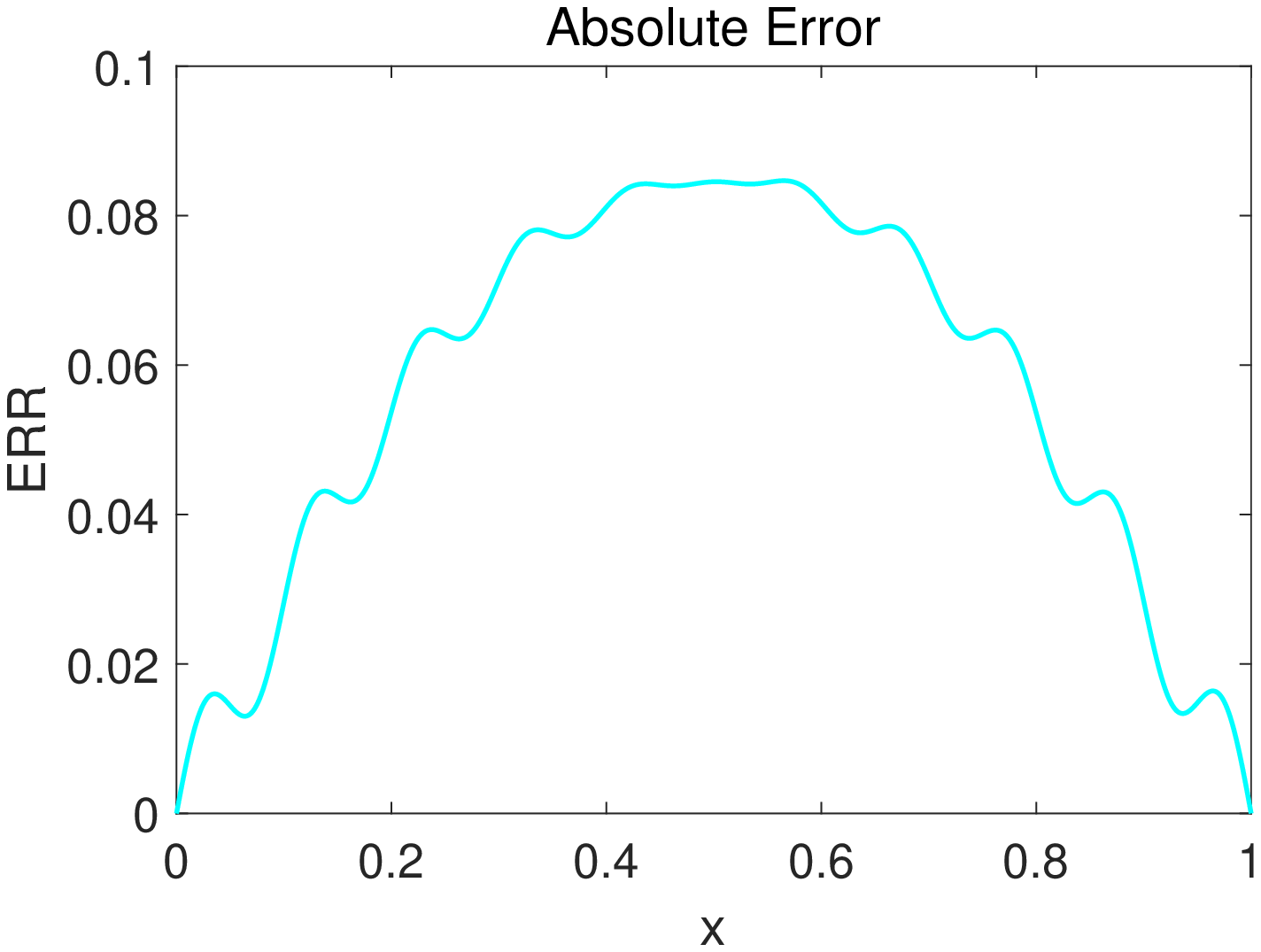}
	}
	\subfigure[REL]{
		\label{0p1err}
		\includegraphics[scale=0.38]{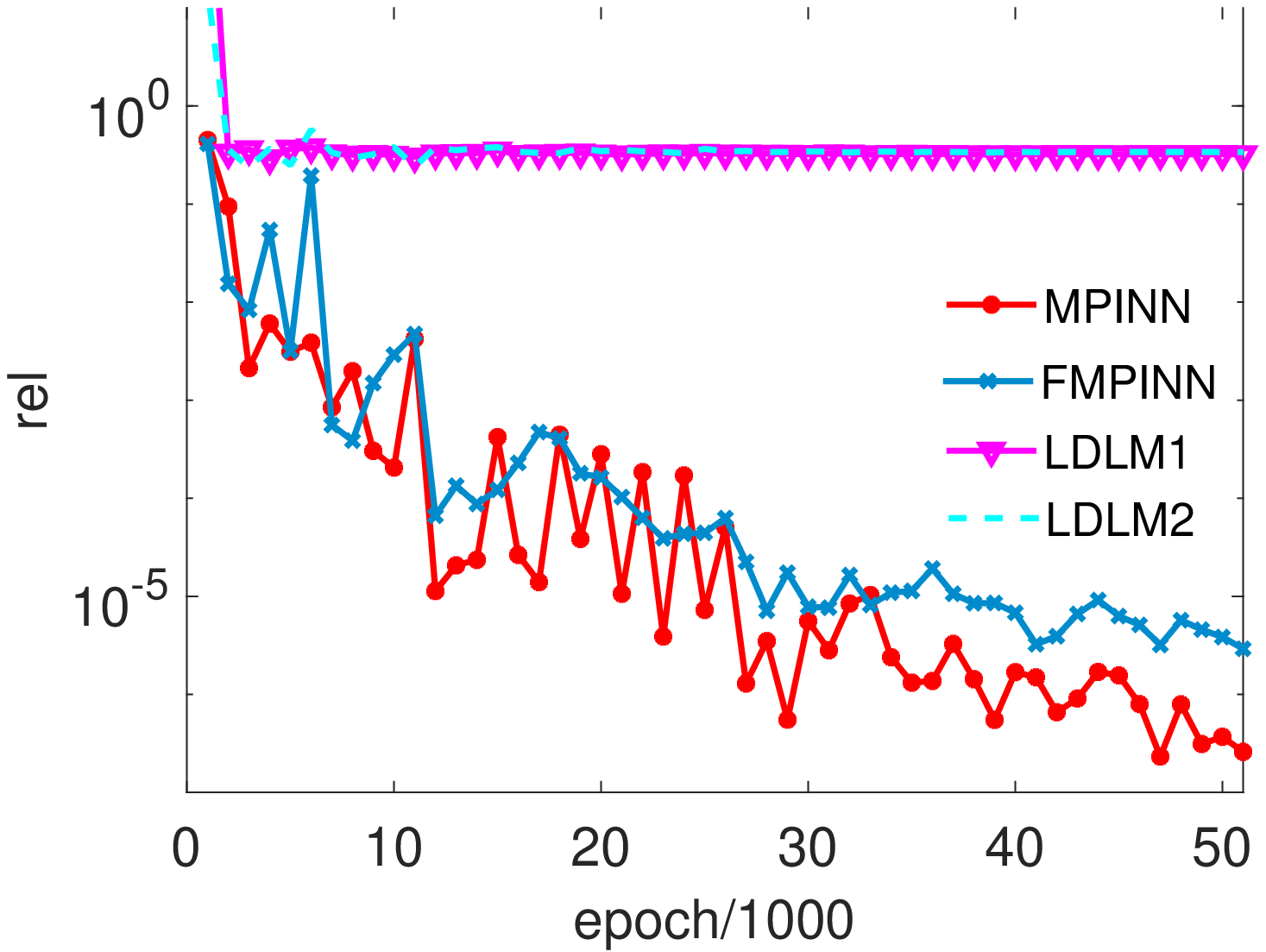}
	}
	\caption{Rough coefficient, loss of flux term and testing results for Example \ref{DiffusionEq_1d_01} when $\varepsilon=0.1$.}
	\label{E1_0p1}
\end{figure}

\begin{figure}
	\centering
	\subfigure[Rough coefficient]{
		\label{0p01Aeps}
		\includegraphics[scale=0.4]{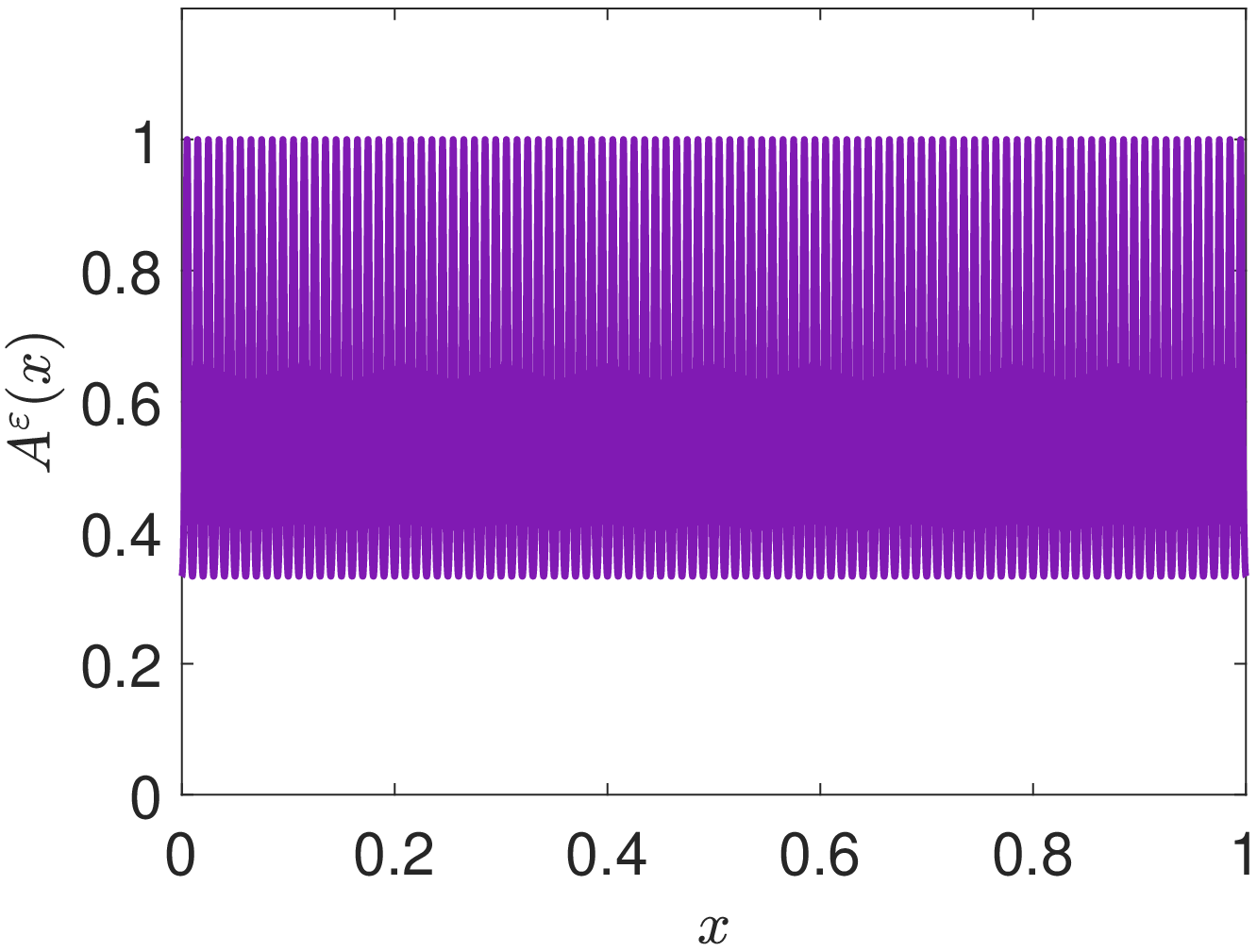}
	}
	\subfigure[loss of flux term for FMPINN]{
		\label{0p01loss2flux}
		\includegraphics[scale=0.4]{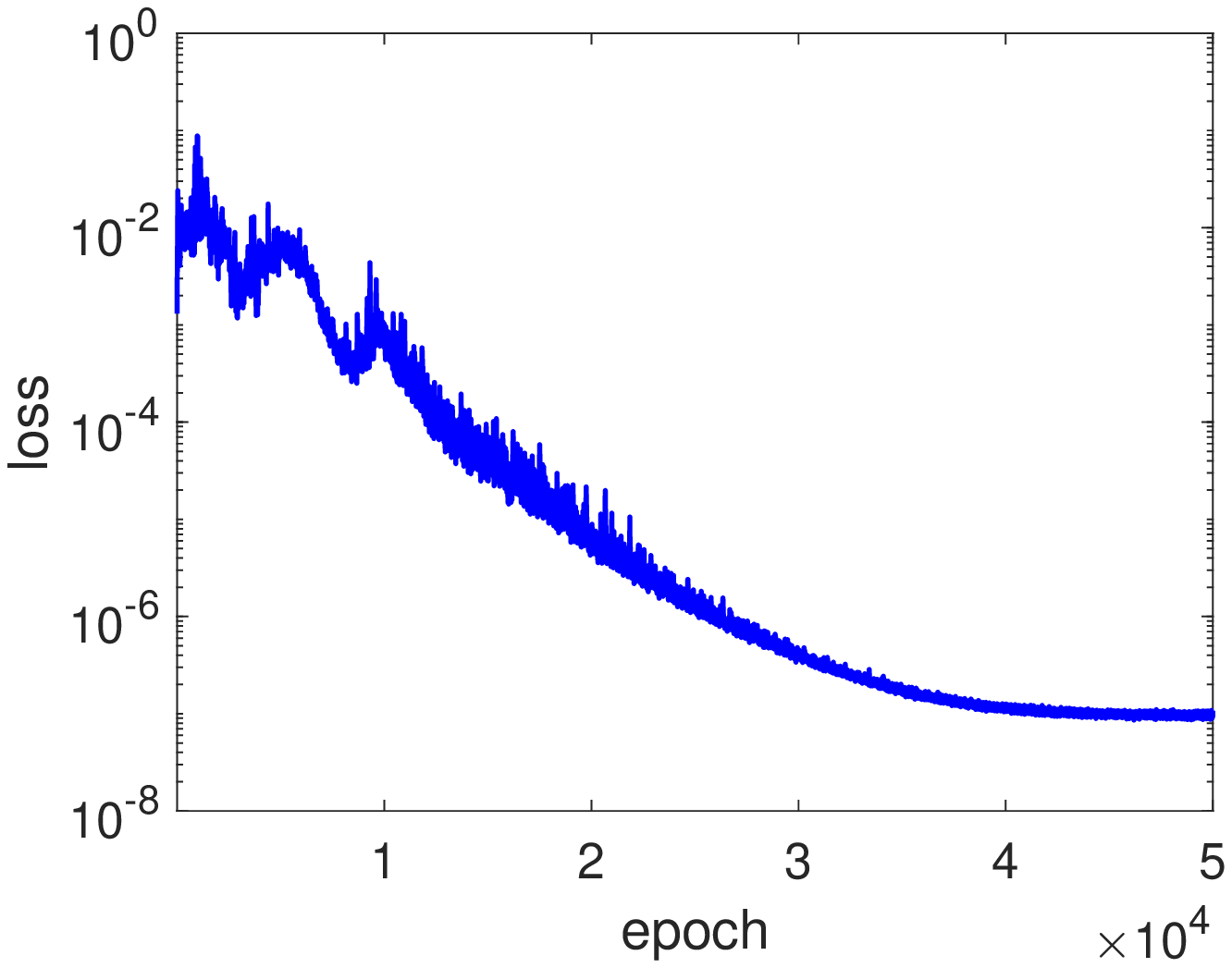}
	}
	\subfigure[Solutions]{
		\label{0p01Solus}
		\includegraphics[scale=0.4]{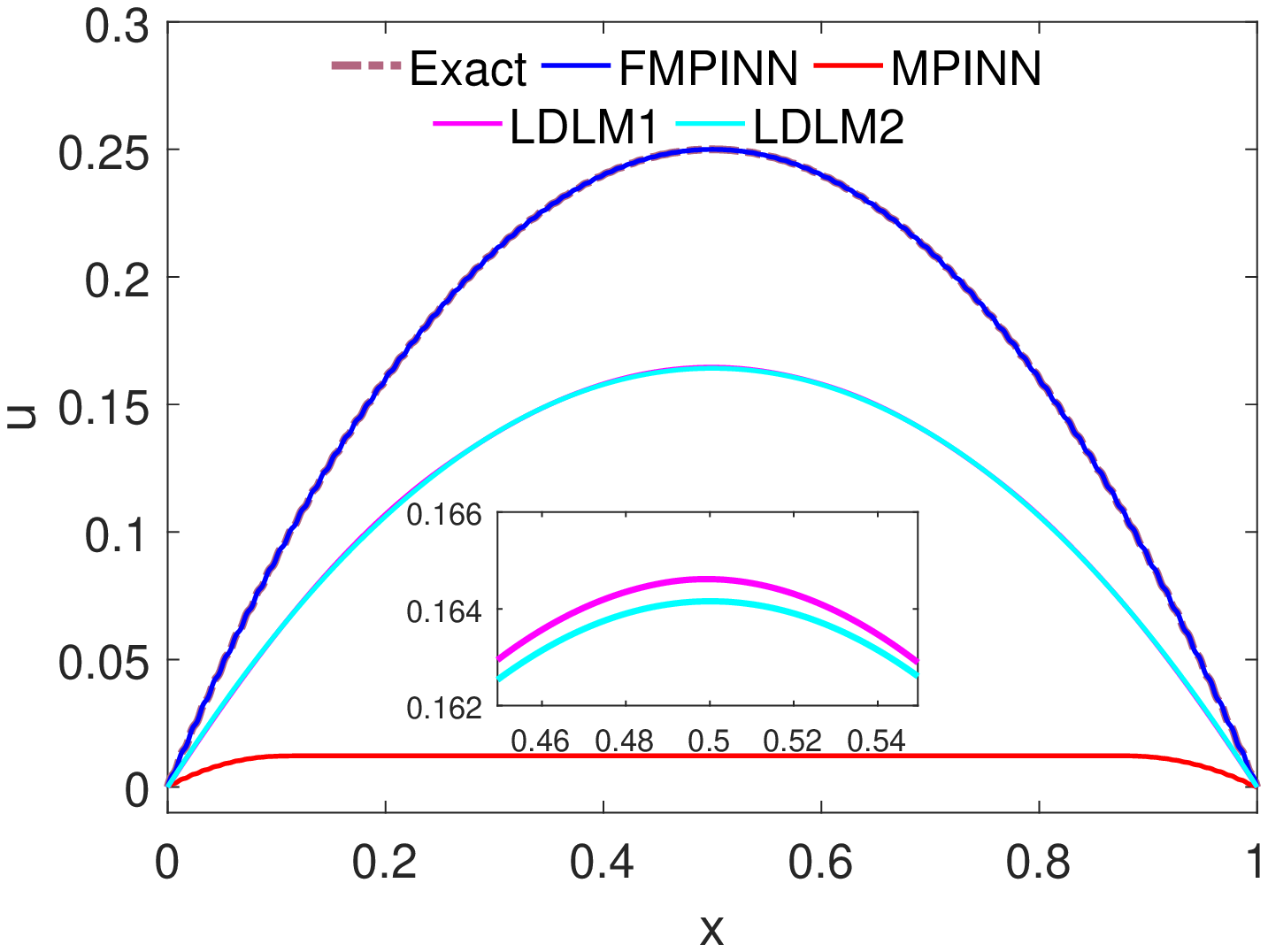}
	}
        \subfigure[point-wise error of FMPINN]{
	    \label{0p01ABS2FMPINN}
		\includegraphics[scale=0.4]{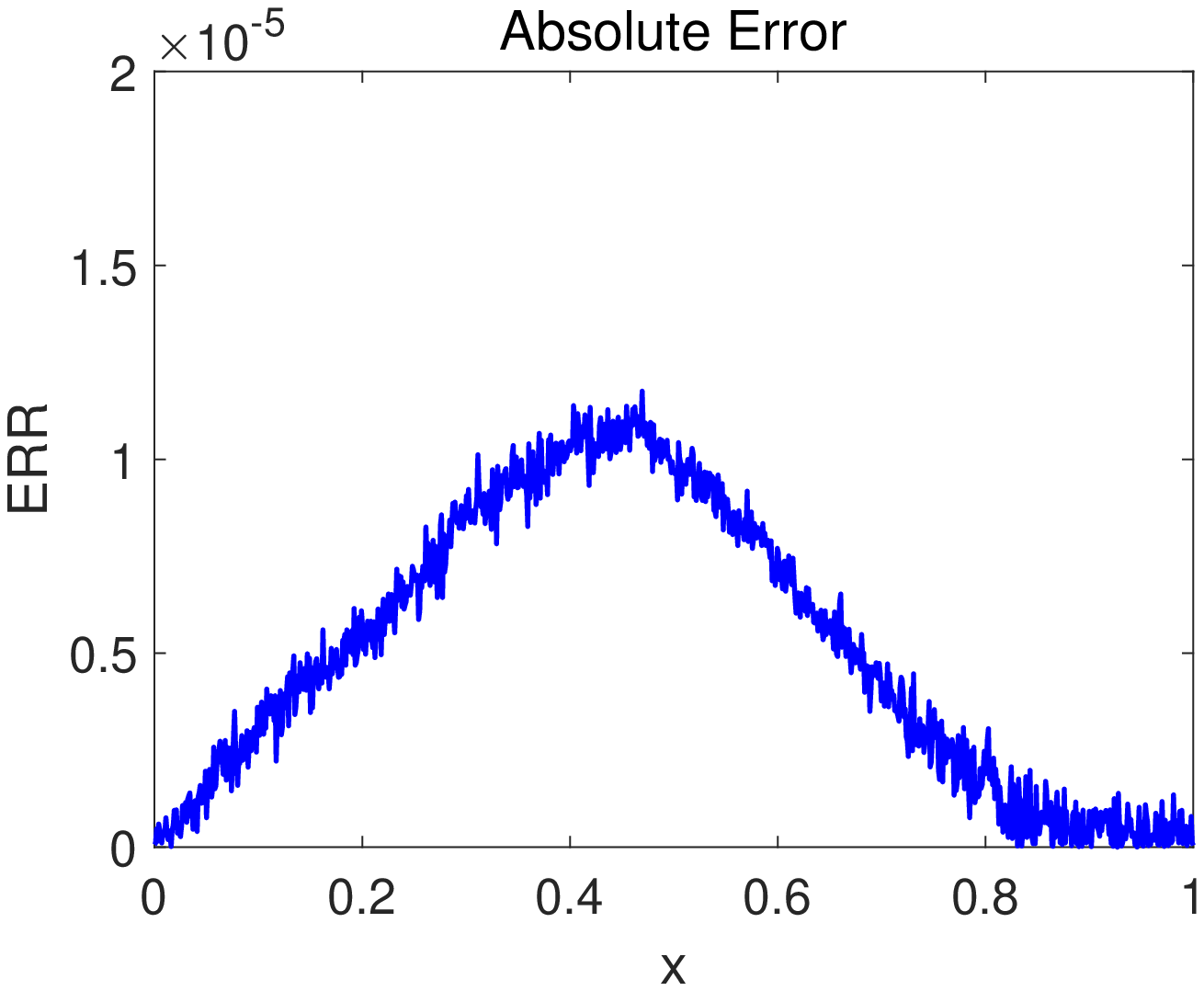}
	}
	\subfigure[point-wise error of MPINN]{
	    \label{0p01ABS2MPINN}
		\includegraphics[scale=0.4]{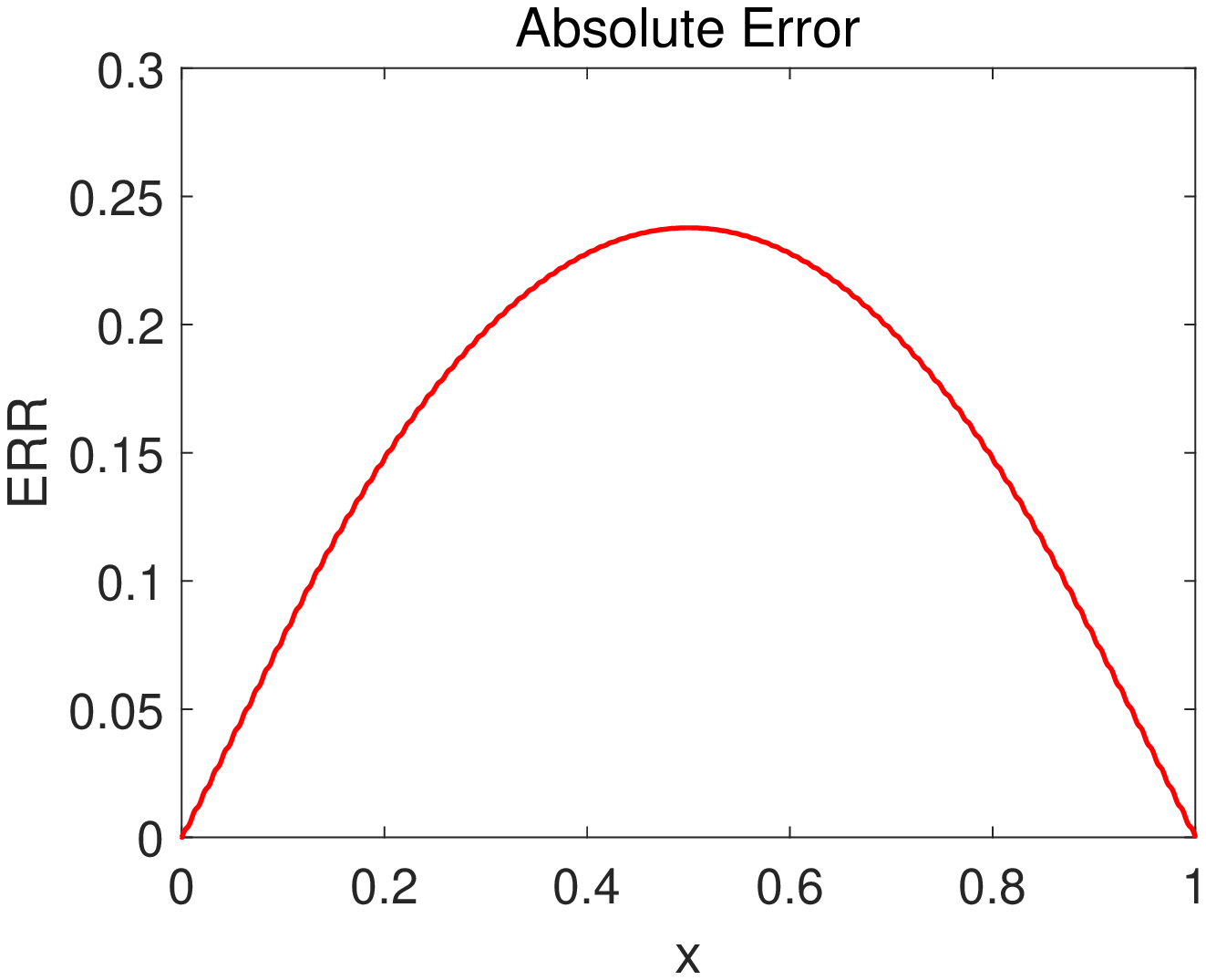}
	}
        \subfigure[point-wise error of LDLM1]{
	    \label{0p01ABS2LDLM1}
		\includegraphics[scale=0.4]{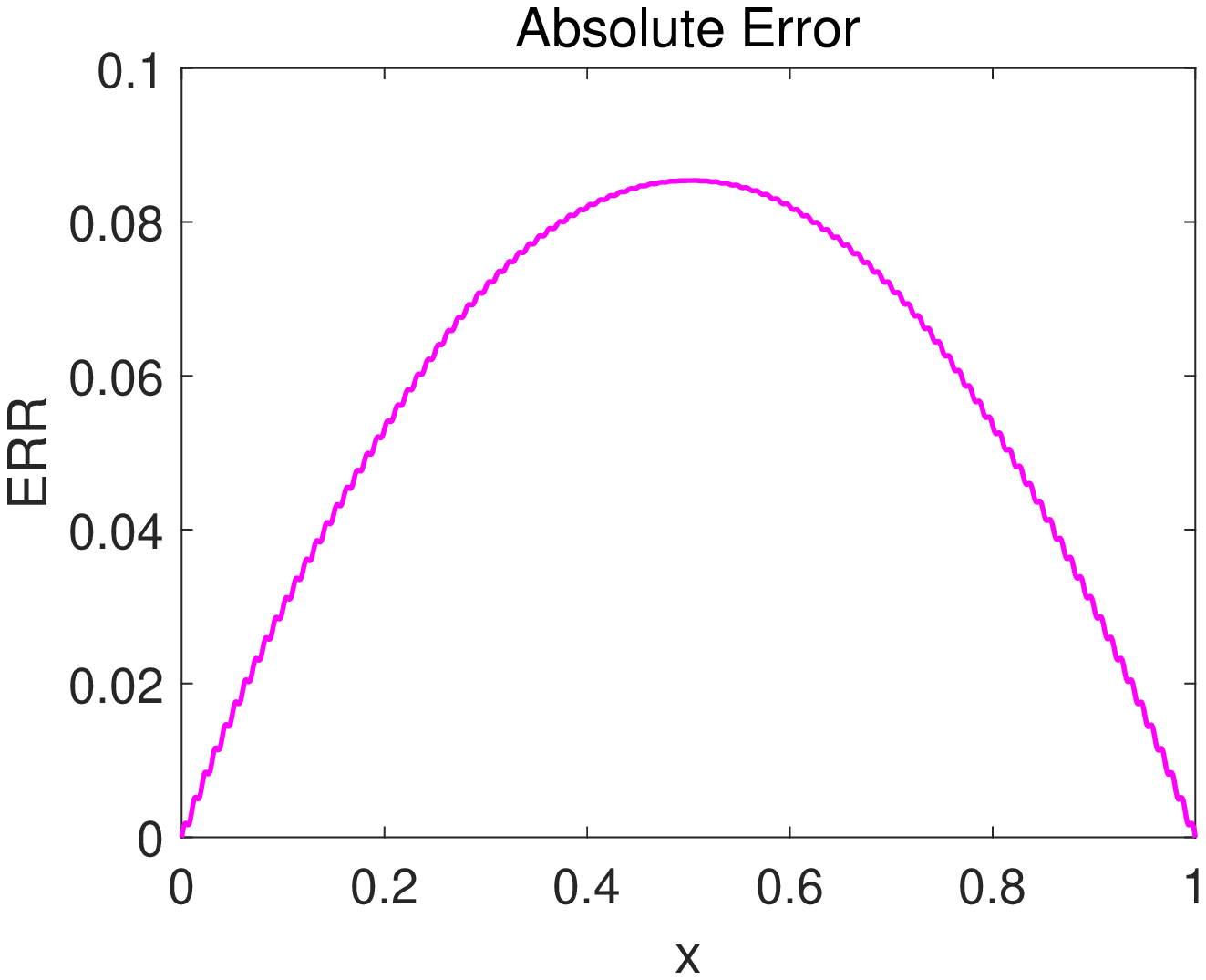}
	}
        \subfigure[point-wise error of LDLM2]{
	    \label{0p01ABS2LDLM2}
		\includegraphics[scale=0.4]{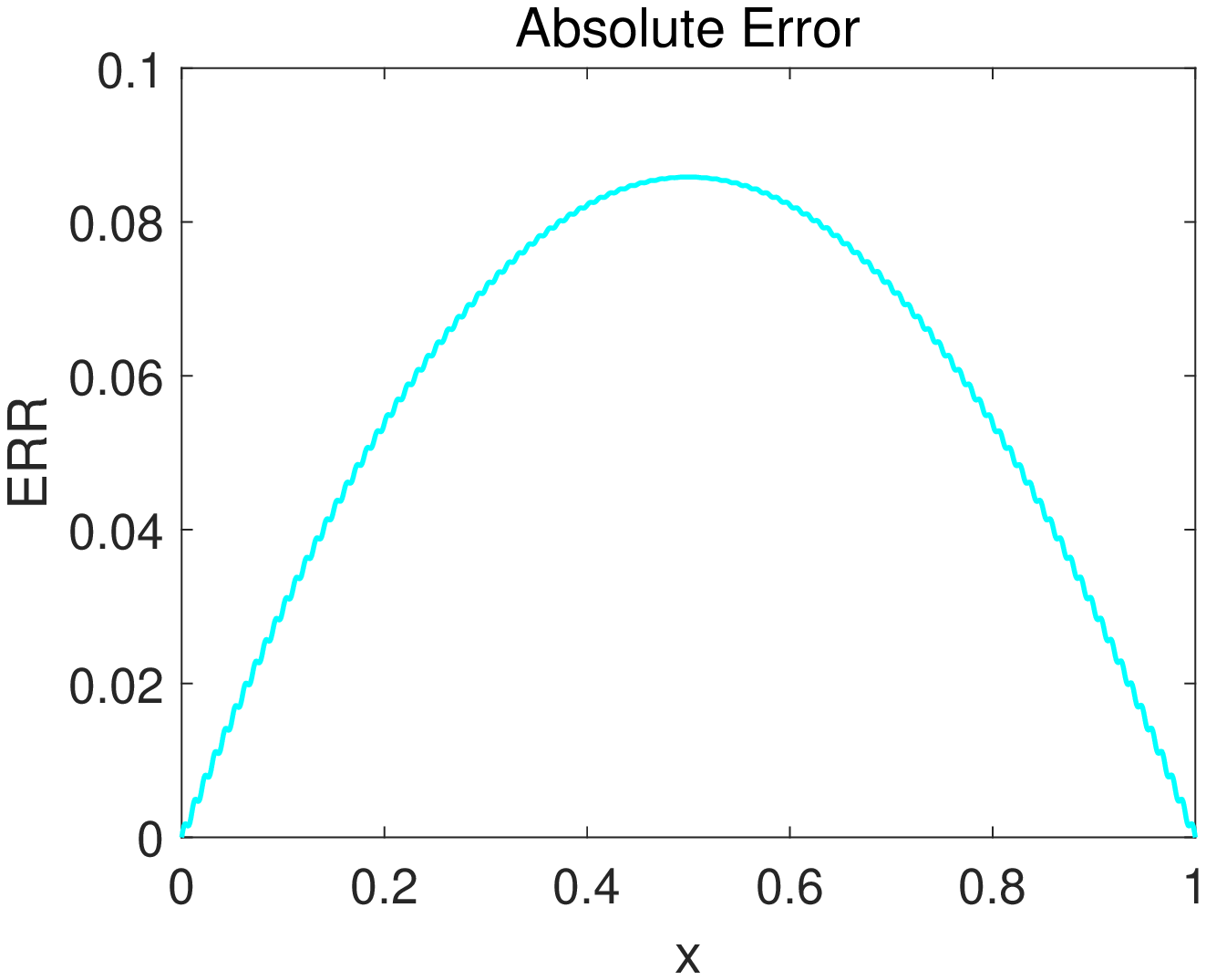}
	}
	\subfigure[REL]{
	    \label{0p01err}
		\includegraphics[scale=0.38]{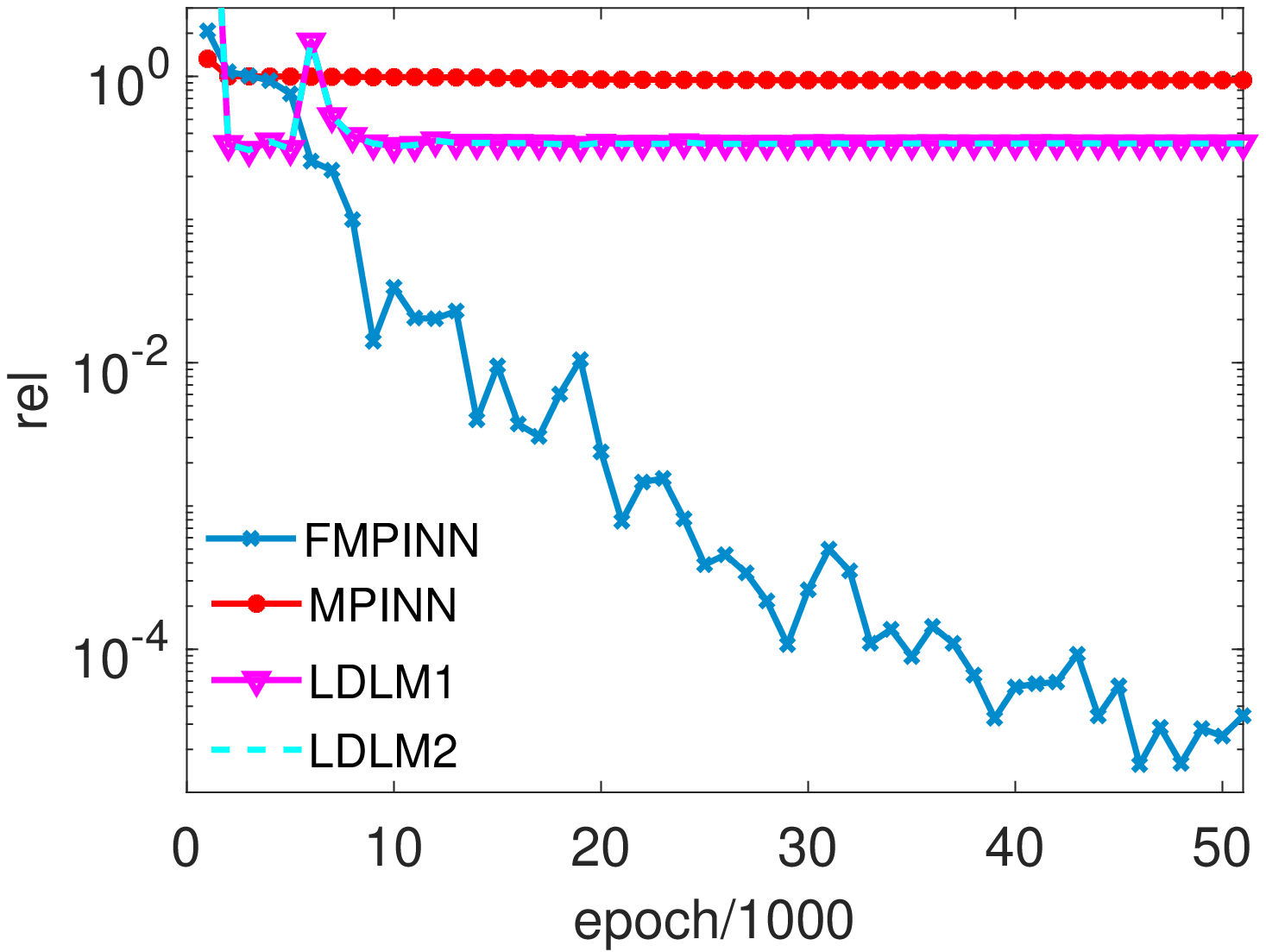}
	}
	\caption{Rough coefficient, loss of flux term and testing results for Example \ref{DiffusionEq_1d_01} when $\varepsilon=0.01$.}
	\label{E1_0p01}
\end{figure}

\begin{figure}
	\centering
	\subfigure[Rough coefficient]{
		\label{0p001Aeps}
		\includegraphics[scale=0.4]{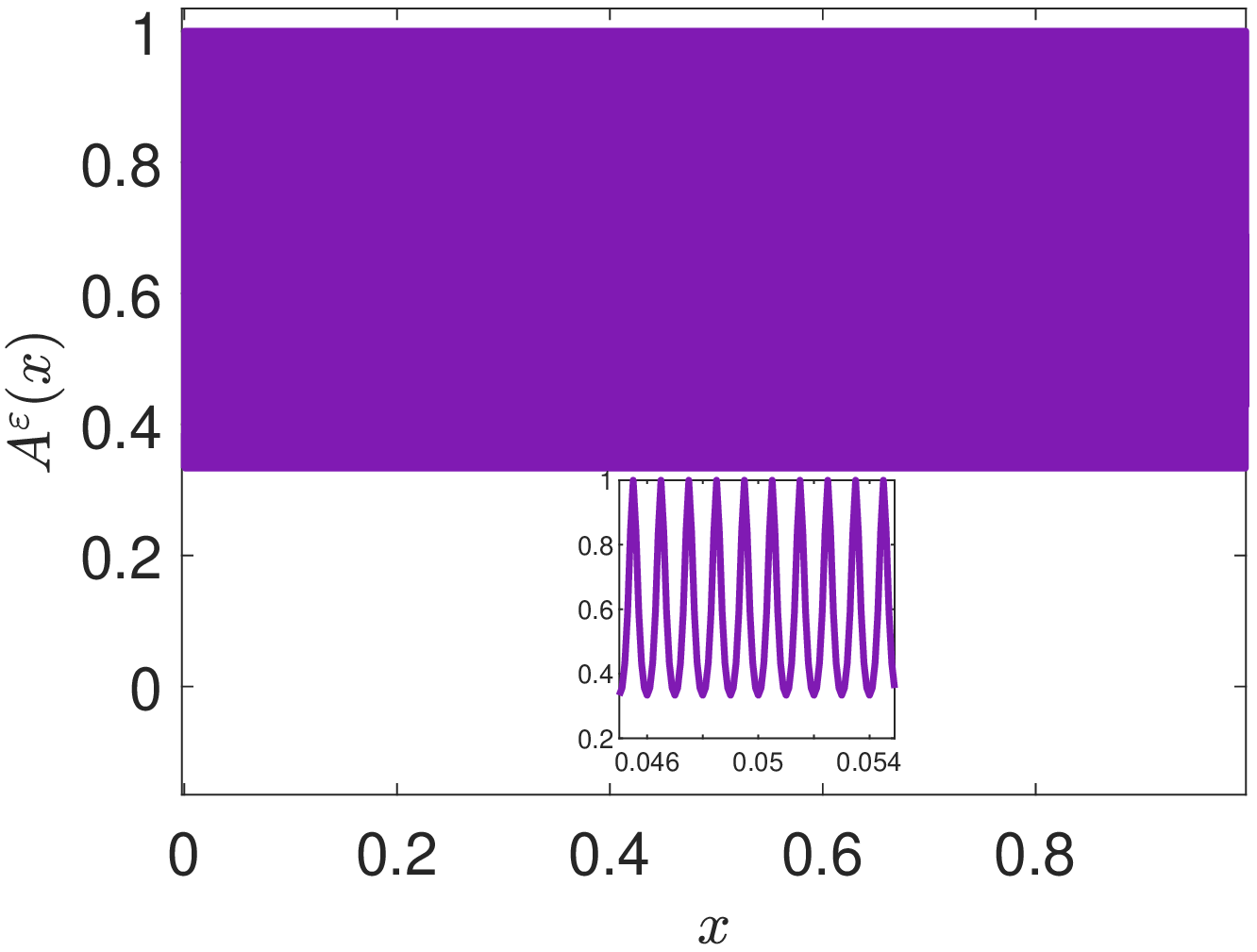}
	}
	\subfigure[loss of flux term for FMPINN]{
		\label{0p001loss2flux}
		\includegraphics[scale=0.4]{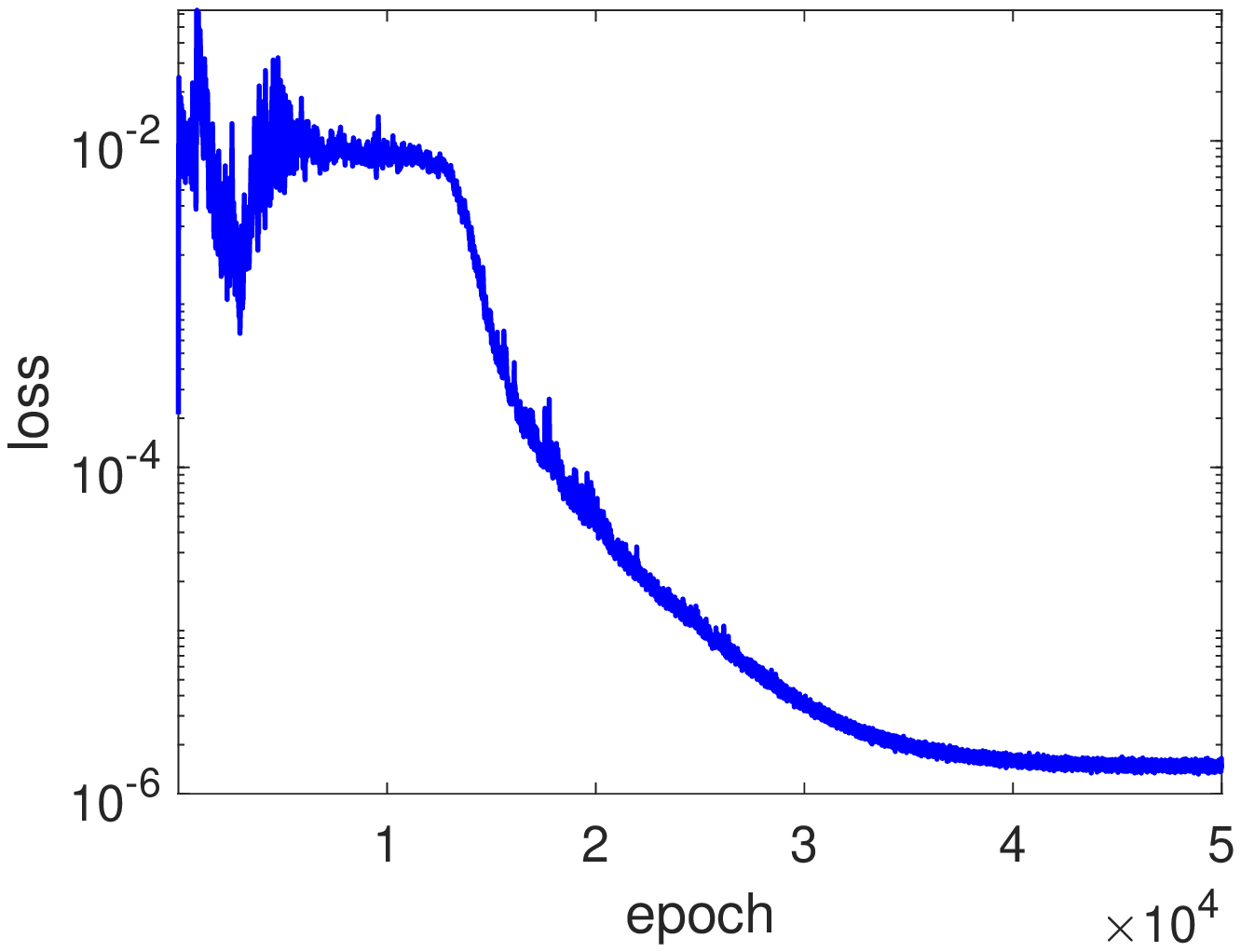}
	}
	\subfigure[Solutions]{
		\label{0p001Solus}
		\includegraphics[scale=0.4]{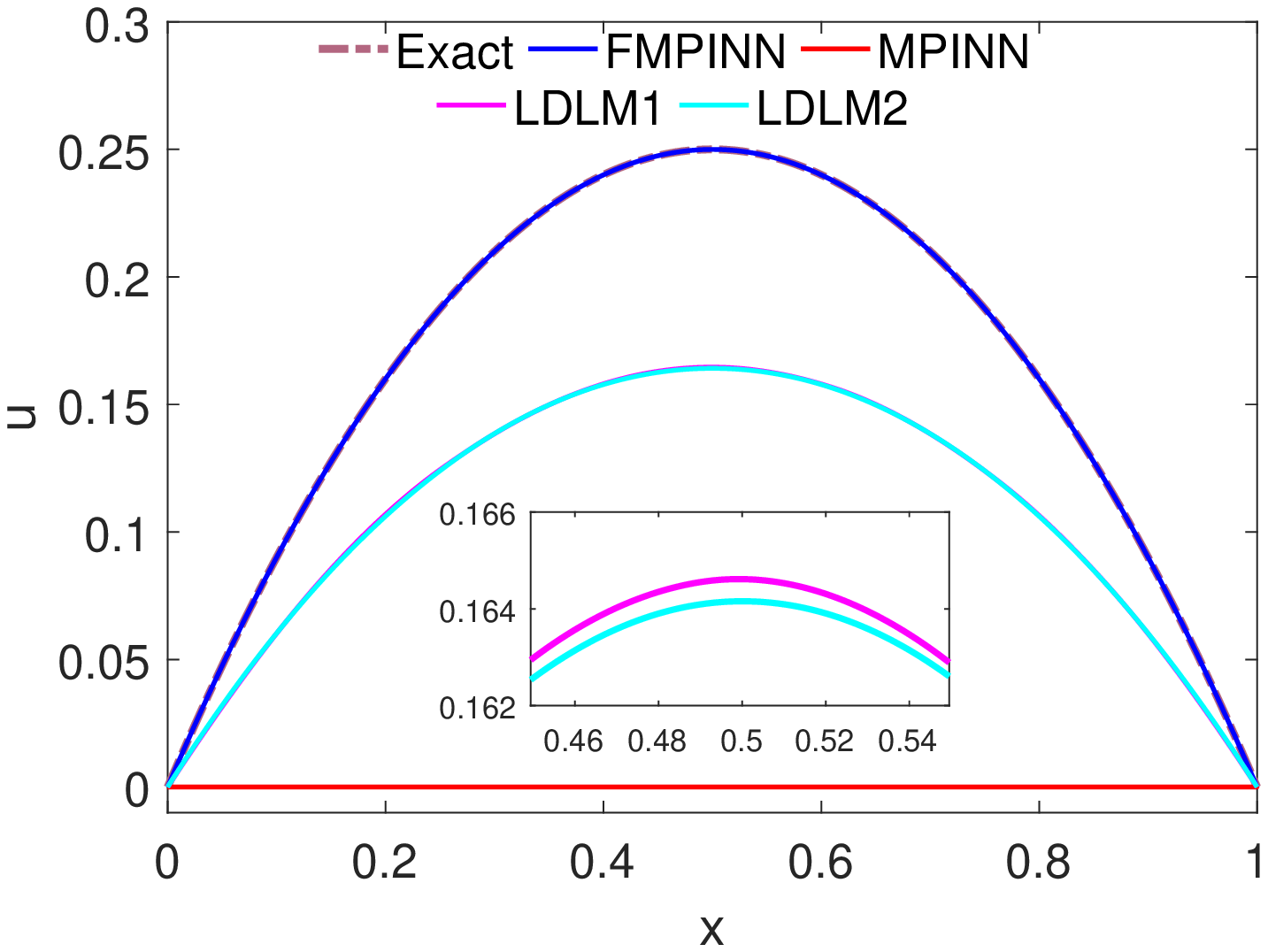}
	}
        \subfigure[point-wise error of FMPINN]{
	    \label{0p001ABS2FMPINN}
		\includegraphics[scale=0.4]{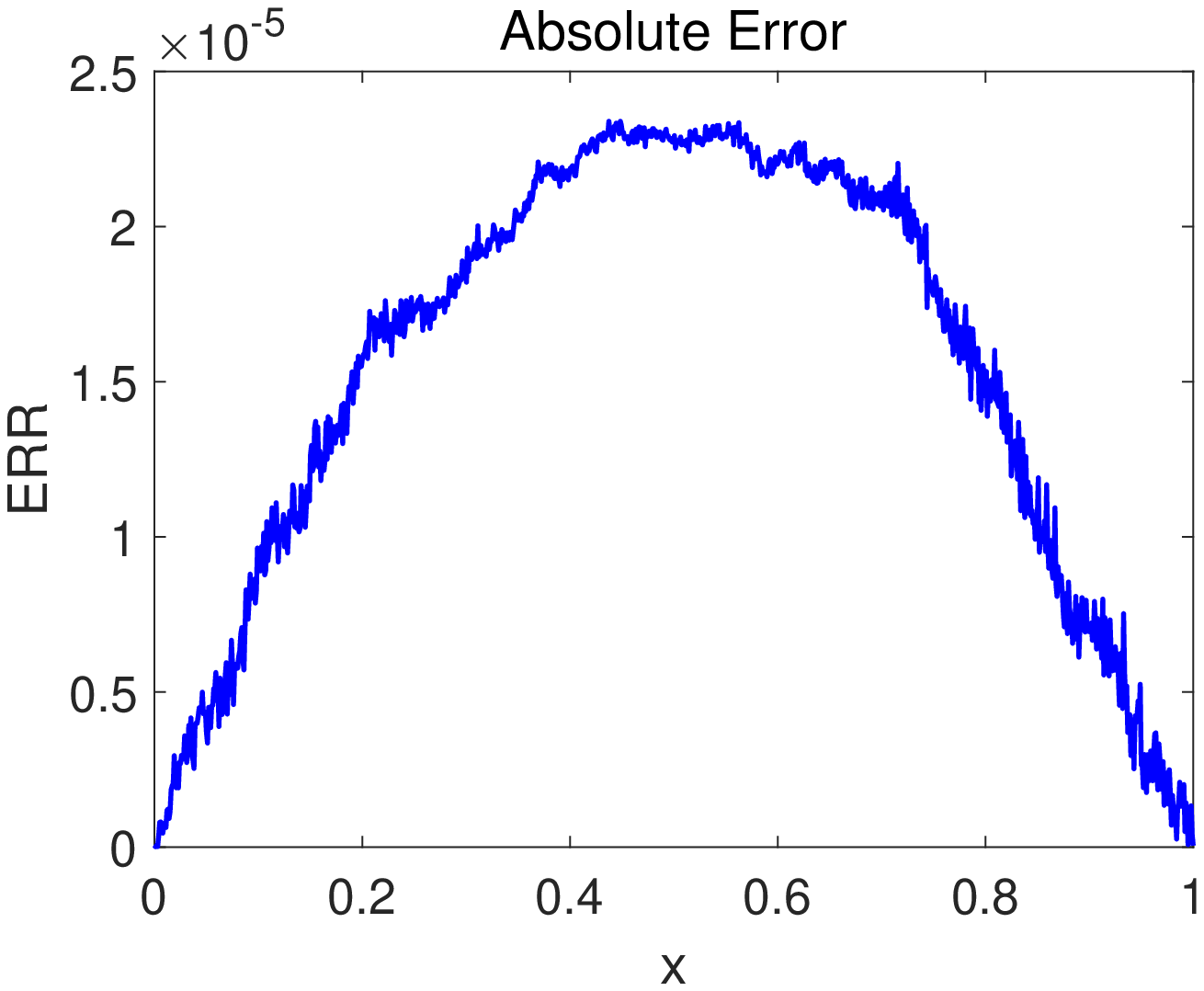}
	}
	\subfigure[point-wise error of MPINN]{
	    \label{0p001ABS2MPINN}
		\includegraphics[scale=0.4]{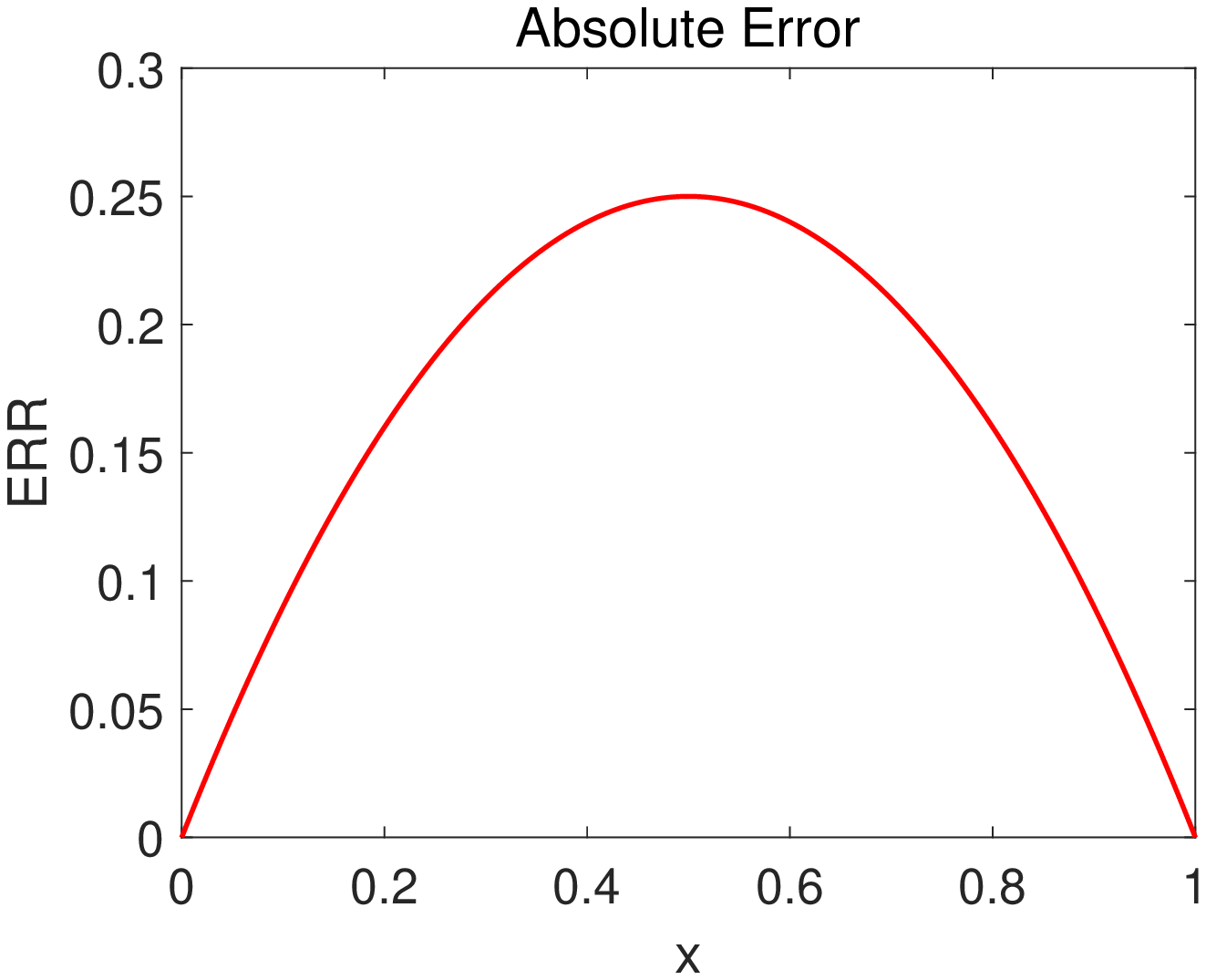}
	}
 	\subfigure[point-wise error of LDLM1]{
	    \label{0p001ABS2LDLM1}
		\includegraphics[scale=0.4]{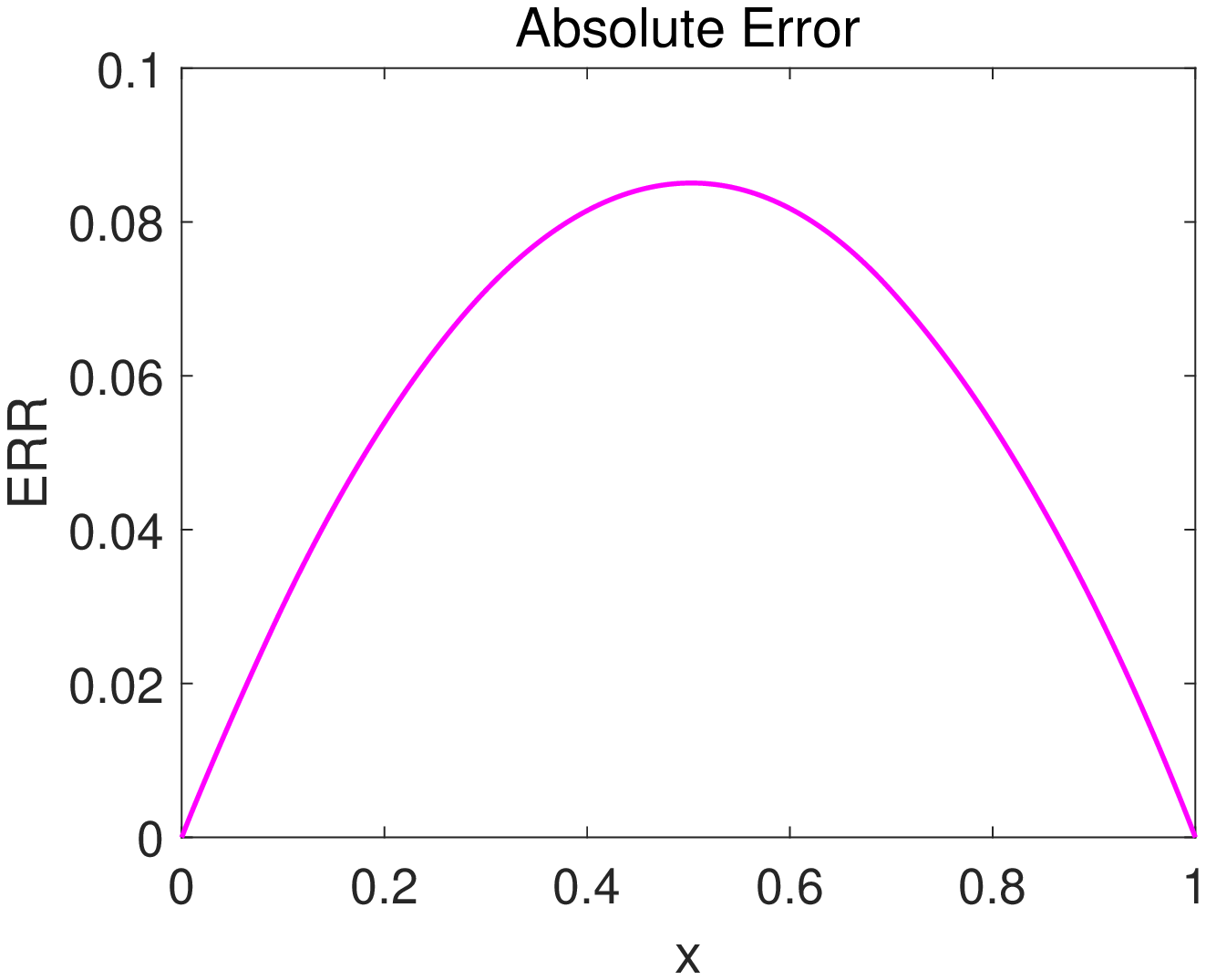}
	}
  	\subfigure[point-wise error of LDLM2]{
	    \label{0p001ABS2LDLM2}
		\includegraphics[scale=0.4]{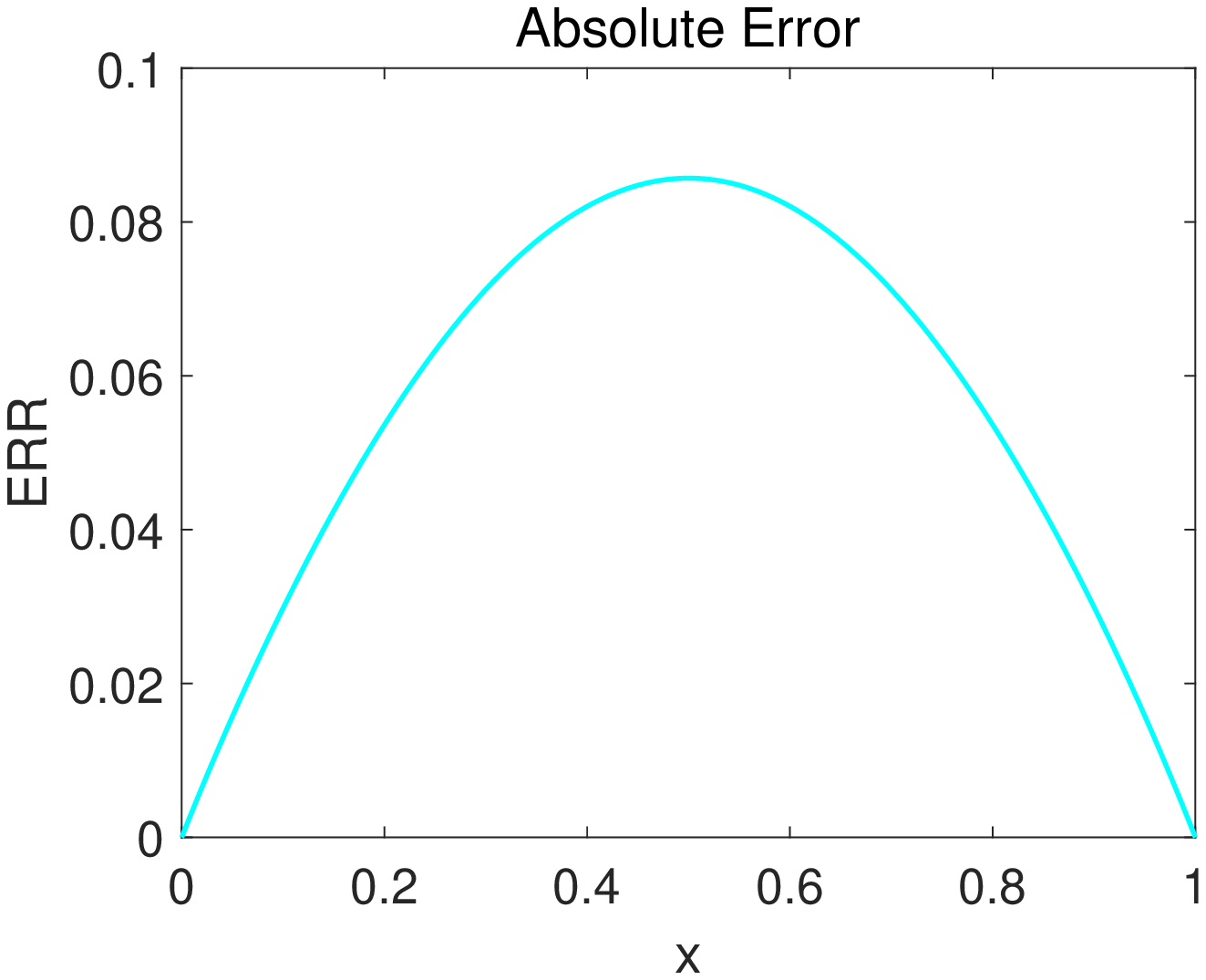}
	}
	\subfigure[REL]{
	    \label{0p001err}
		\includegraphics[scale=0.38]{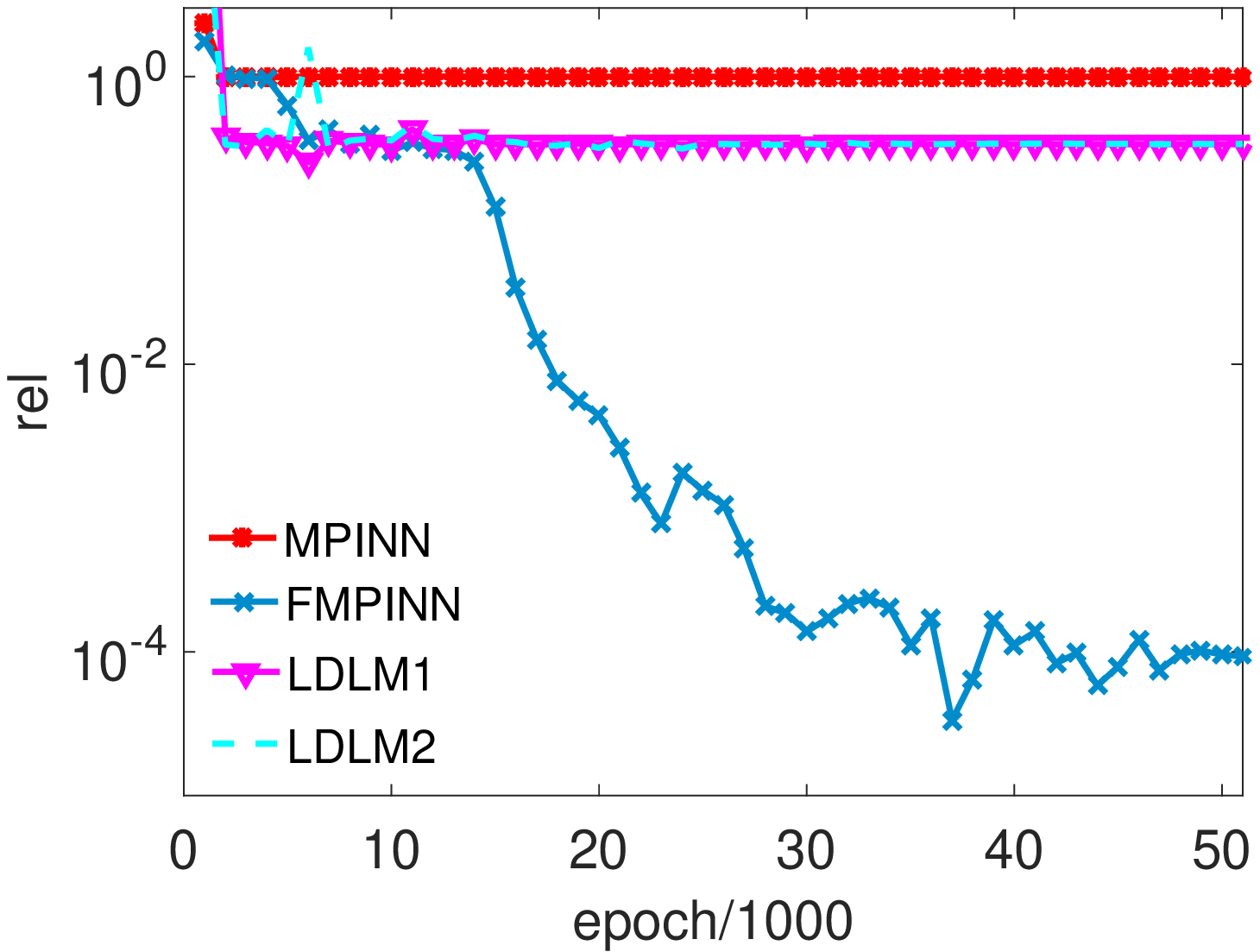}
	}
	\caption{Rough coefficient, loss of flux term and testing results for Example \ref{DiffusionEq_1d_01} when $\varepsilon=0.001$.}
	\label{E1_0p001}
\end{figure}

In this example, we use the FMPINN, MPINN, LDLM1, and LDLM2 models to solve \eqref{eq:multiscale} when $\varepsilon=0.1, 0.01$ and $0.001$, respectively. The size of hidden layer for each subnetwork of FMPINN and MPINN is set as $(30,40,30,30,30)$ and the balance parameter $\beta$ in \eqref{loss2ellptic} is set as 10. The hidden layer's size for LDLM is set as $(300,400,300,300,300)$. Their parameters' numbers are comparable. At each training step, we randomly sample 3000 points inside the $[0,1]$ and 500 boundary points as a training dataset. In addition, the testing dataset includes 1000 equidistant samples from $[0, 1]$. All models are trained for 50000 epochs. We depict the related experiment results in Figs. \ref{E1_0p1}, \ref{E1_0p01} and \ref{E1_0p001}, respectively. Meantime, the final relative errors and total running time are listed in Table \ref{results2multiscale}. 
\begin{table}[H]
	\centering
	\caption{The relative error and running time of FMPINN, MPINN, LDLM1 and LDLM2 for Example \ref{DiffusionEq_1d_01}}
	\label{results2multiscale}
	\begin{tabular}{|c|c|c|c|c|c|c|c|c|c|}
		\hline
		            &\multicolumn{4}{|c|}{REL} &\multicolumn{4}{|c|}{Total time(s)}              \\ \hline
		 $\varepsilon$ &FMPINN   &MPINN   &LDLM1  &LDLM2  &FMPINN  &MPINN   &LDLM1   &LDLM2 \\ \hline
              0.1       &2.92e-6  &2.60e-7 &0.3227 &0.3389 &680.734 &865.849 &345.791 &373.537\\ \hline
		0.01        &3.43e-5  &0.94    &0.3397 &0.3406 &689.729 &868.199 &351.451 &377.089\\ \hline
		0.001       &9.28e-5  &0.99    &0.3389 &0.3398 &691.458 &875.297 &358.435 &388.273\\ \hline
	\end{tabular}
\end{table}

Based on these figures, the FMPINN model can perfectly capture the oscillation of the exact solution for $\varepsilon=0.1, 0.01$ and $0.001$, but LDLM models are not convergent for these cases. At the same time, the performance of MPINN competes with that of FMPINN when $\varepsilon=0.1$. However, the MPINN model fails to solve the multi-scale problem for $\varepsilon=0.1$ and $0.01$. Compared to $\varepsilon=0.01$, the rough coefficient  $A^{\varepsilon}$ with $\varepsilon=0.001$ have more oscillation in the interval $[0,1]$, but the FMPINN still can keep its remarkable performance. According to the point-wise errors in Figs. \ref{0p1ABS2FMPINN}, \ref{0p01ABS2FMPINN} and \ref{0p001ABS2FMPINN} and the relative error in Figs. \ref{0p1err}, \ref{0p01err} and \ref{0p001err}, we can conclude that the FMPINN is able to approximate high-precisely the exact solution of \eqref{eq:multiscale} in one-dimensional space. In addition, the total time in Table \ref{results2multiscale} shows the running time of FMPINN is less than that of MPINN for 50000 training epochs. 

\begin{figure}[H]
	\centering
	\subfigure[Loss for flux term]{
		\label{study2beta_lossAdu}
		\includegraphics[scale=0.4]{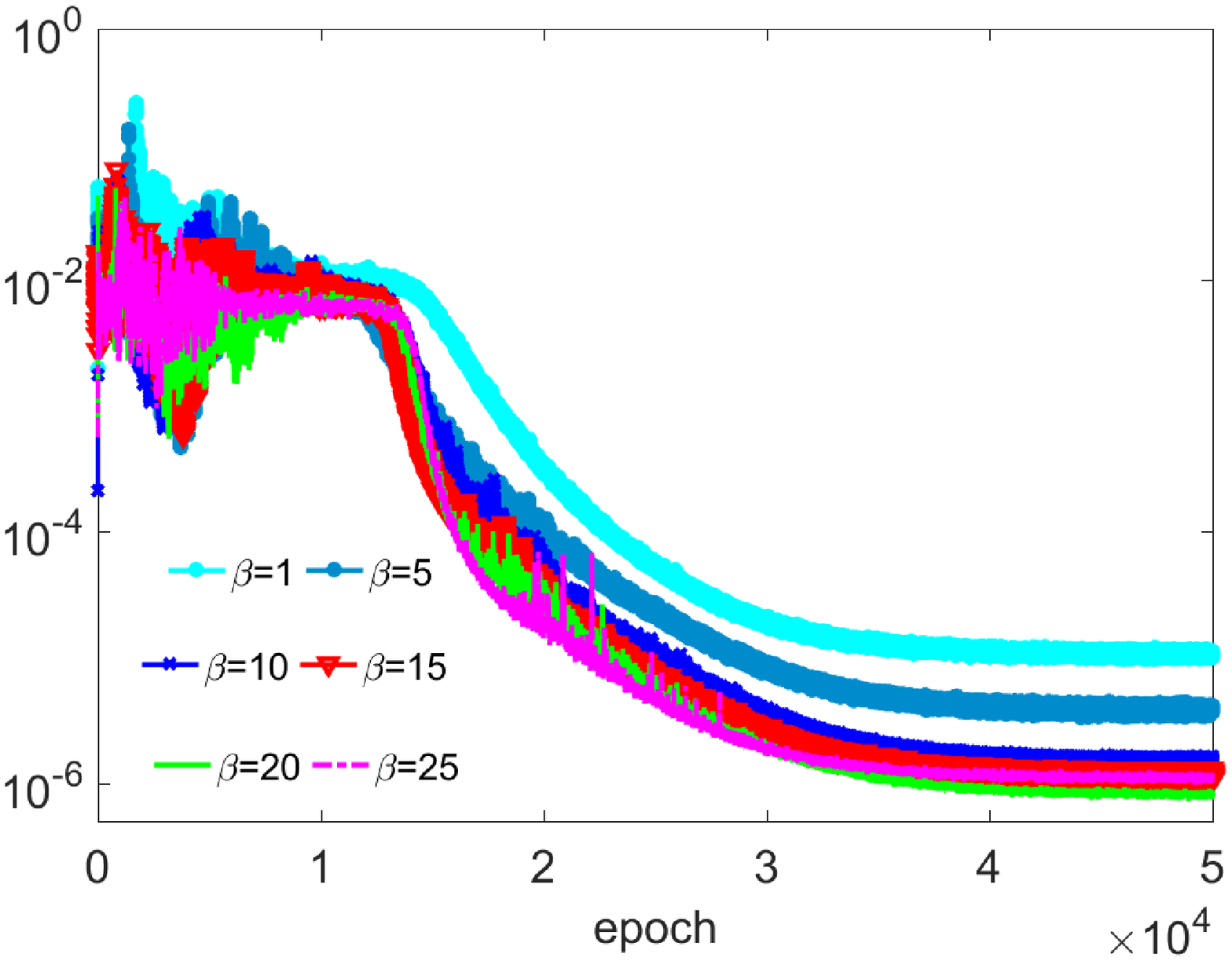}
	}
	\subfigure[REL]{
		\label{study2beta_rel}
		\includegraphics[scale=0.475]{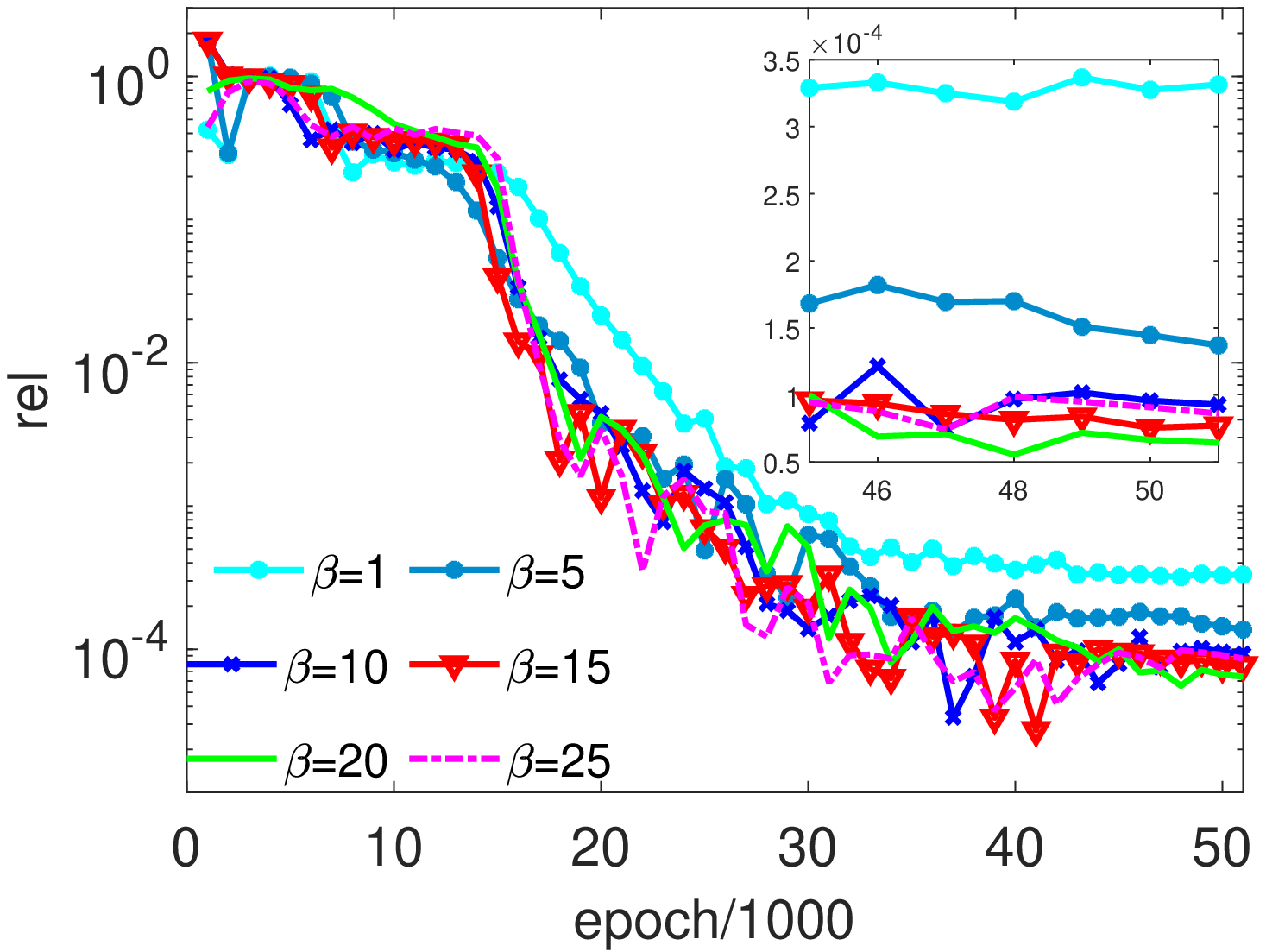}
	}
	\caption{The loss of flux term VS training epoch and the relative error VS testing epoch for Example \ref{DiffusionEq_1d_01} when $\varepsilon=0.001$.}
	\label{Beta_study}
\end{figure}

\emph{Influence of hyper-parameter $\beta$:} In the previous tests, the parameter $\beta$ was initially set to 10. Now, we study the influence of $\beta$ for our FMPINN model. In these tests, we take $\varepsilon=0.001$ in \eqref{DiffusionEq_1d_01_aeps}, and consider values of $\beta$ equal to $1, 5, 10, 15, 20$ and $25$, while keeping all other parameters fixed. All models with different $\beta$ values are trained for 50000 epochs. Fig. \ref{Beta_study} plots the results of flux loss for the training process as well as the relative error for testing. Additionally, the final relative errors obtained from the tests are listed in Table \ref{Table2study_beta}.
\begin{table}[H]
	\centering
	\caption{ The relative error of FMPINN model VS various $\beta$ for Example \ref{DiffusionEq_1d_01}.}
	\label{Table2study_beta}
	\begin{tabular}{|c|c|c|c|c|c|c|c|}
		\hline
		$\beta$  &1        &5       &10       &15       &20       & 25     \\  \hline
		  REL    &3.31e-4  &1.37e-4 &9.27e-5  &7.75e-5  &6.40e-5  &8.60e-5 \\  \hline
	\end{tabular}
\end{table}

According to the above results in Fig. \ref{Beta_study} and Table \ref{Table2study_beta}, it can be observed that the FMPINN model exhibits remarkable and stable performance across different values of $\beta$. The performances of the FMPINN model for $\beta=1$ and $\beta=5$ are slightly weaker than that of other cases. The loss of flux term is also stable and consistent with the trendlines of REL. Therefore, for the subsequent tests, we will continue to set $\beta=10$.

\begin{example}\label{Diffusion1D3scale}
Let us attempt to solve the following three-scale problem  with Dirichlet boundary in $\Omega =[0, 1]$. In which,
\begin{equation*}
  A^{\varepsilon}(x) = \left(2+\cos\left(2\pi\frac{x}{\varepsilon_1}\right)\right)\left(2+\cos\left(2\pi\frac{x}{\varepsilon_2}\right)\right)
\end{equation*}
with two small parameter $\varepsilon_1,\varepsilon_2>0$ such that $\varepsilon_1^{-1},\varepsilon_2^{-1}\in\mathbb{N}^+$ and  an exact solution is given by 
	\begin{equation}\label{Diffusion1D3scaleU}
      u^{\varepsilon}(x) = x-x^2+\frac{\varepsilon_1}{4\pi}\sin\left(2\pi\frac{x}{\varepsilon_1}\right) +\frac{\varepsilon_2}{4\pi}\sin\left(2\pi\frac{x}{\varepsilon_2}\right).
	\end{equation}
Clearly, $u^\varepsilon(0)=u^\varepsilon(1)=0$. One can obtain the force side after careful computation, we omit it here. 
\end{example}

\begin{figure}[H]
	\centering
	\subfigure[Rough coefficient]{
		\label{3scale2Aeps}
		\includegraphics[scale=0.4]{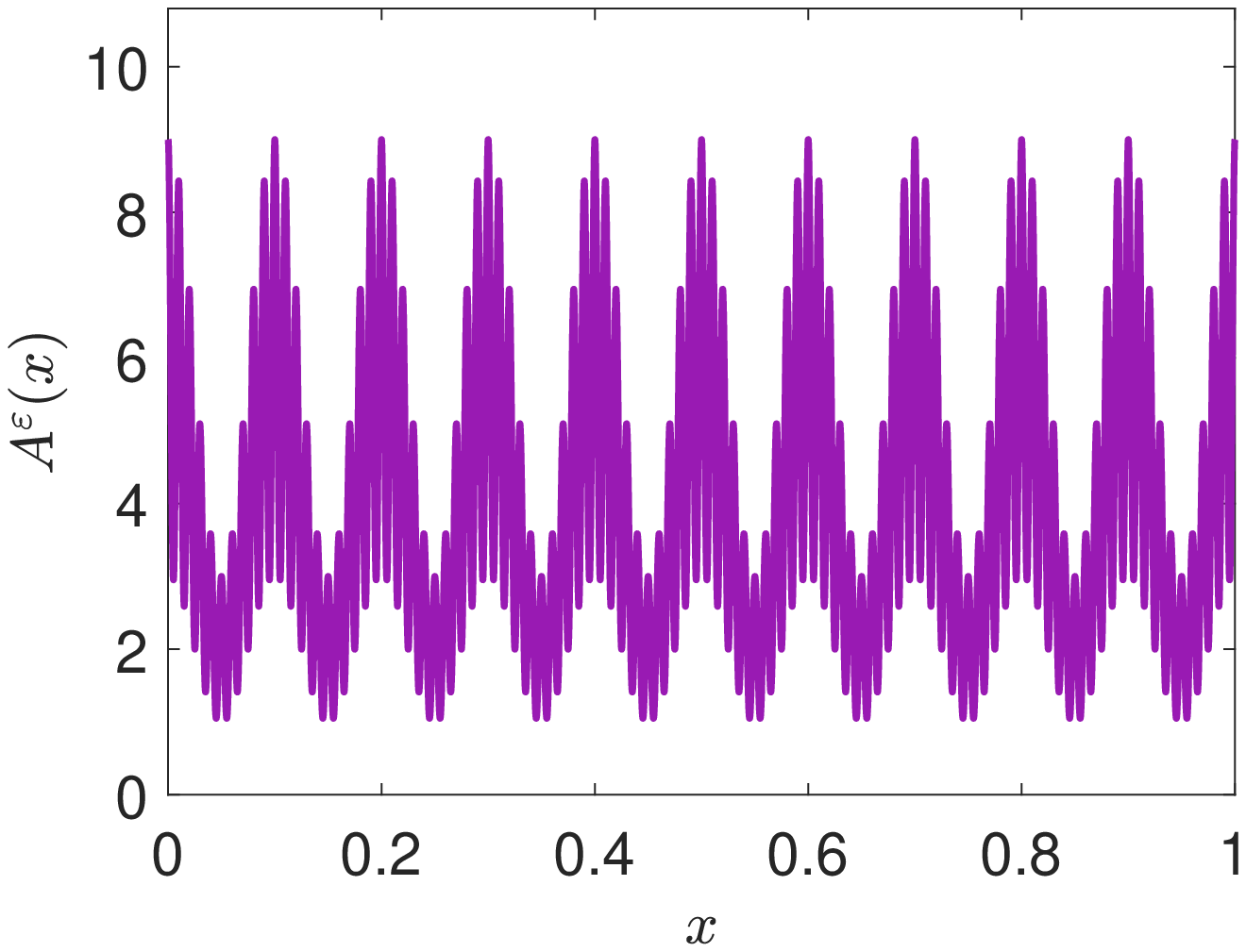}
	}
	\subfigure[loss of flux term for FMPINN]{
		\label{3scale2loss_flux}
		\includegraphics[scale=0.4]{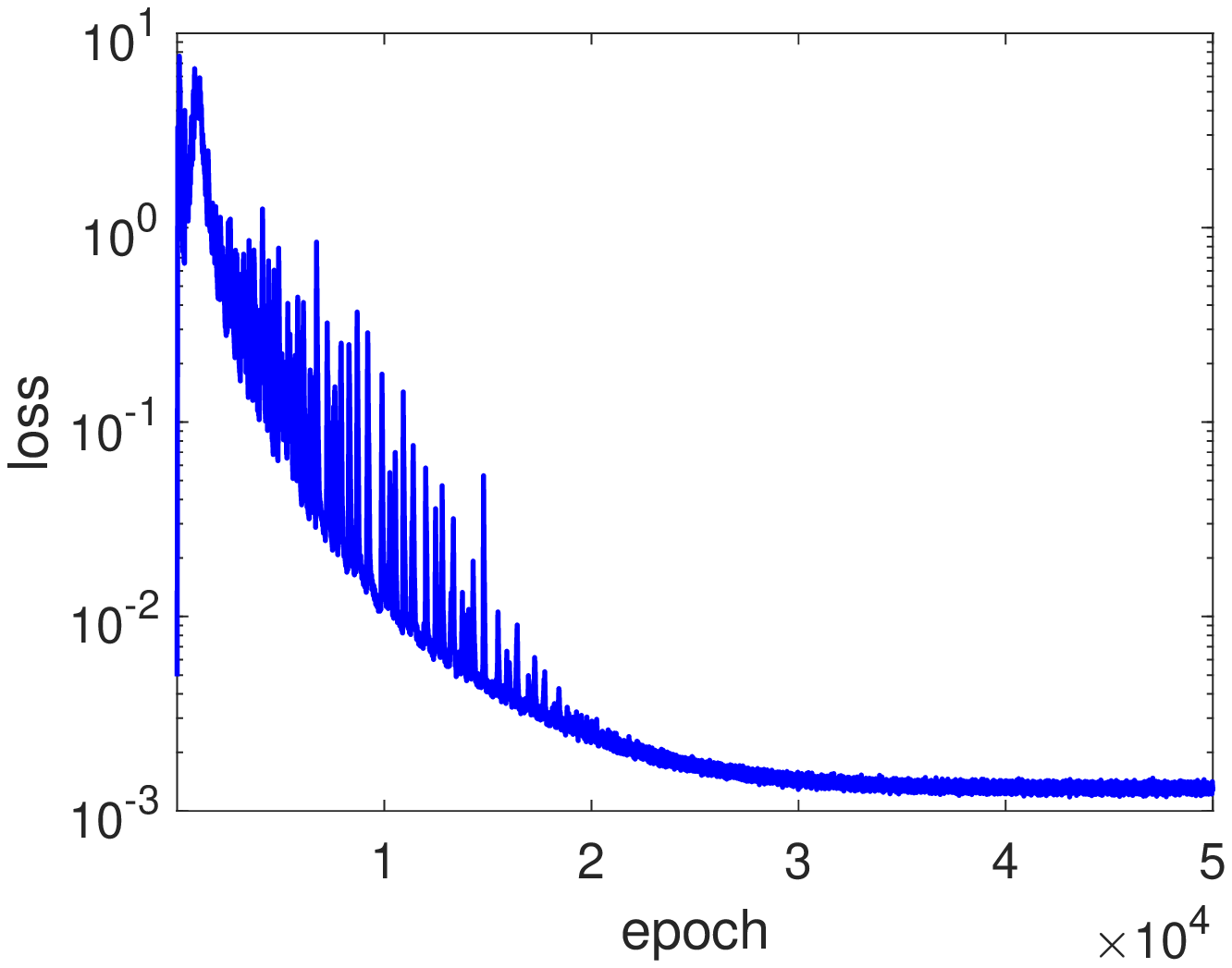}
	}
	\subfigure[Solutions]{
		\label{3scale2Solus}
		\includegraphics[scale=0.4]{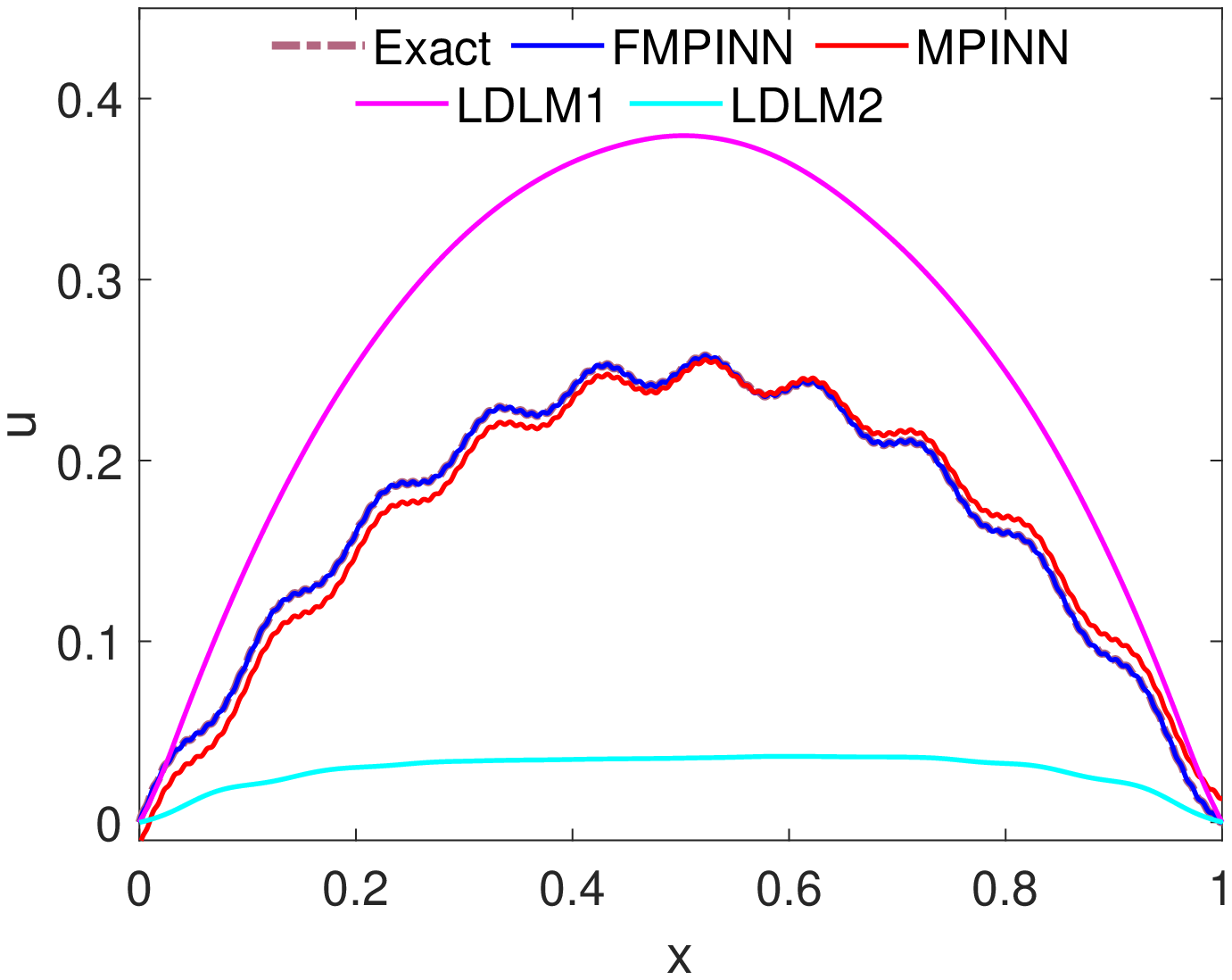}
	}
	\subfigure[point-wise error of FMPINN]{
		\label{3scale2absFMPINN}
		\includegraphics[scale=0.4]{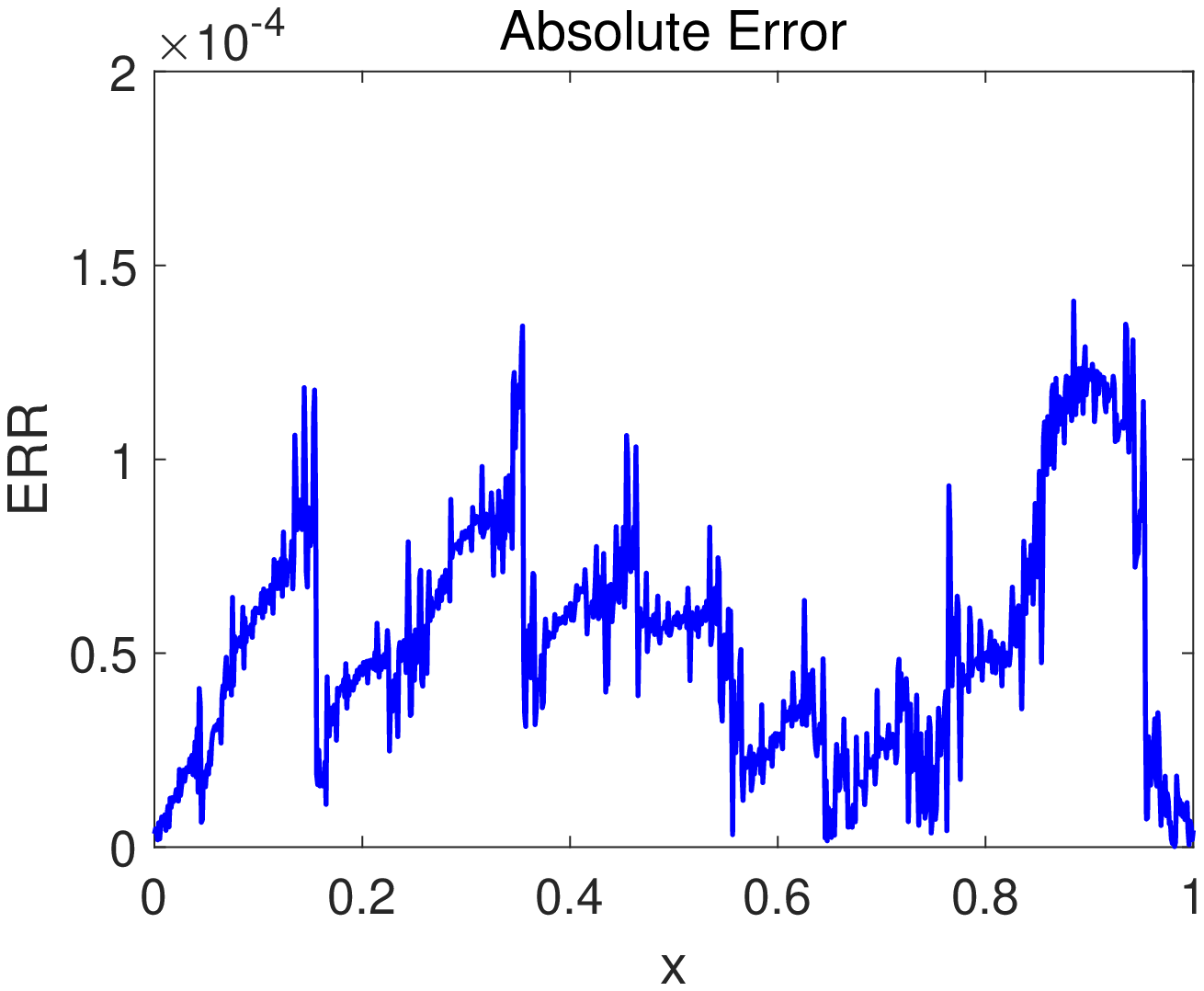}
	}
        \subfigure[point-wise error of MPINN]{
		\label{3scale2absMPINN}
		\includegraphics[scale=0.4]{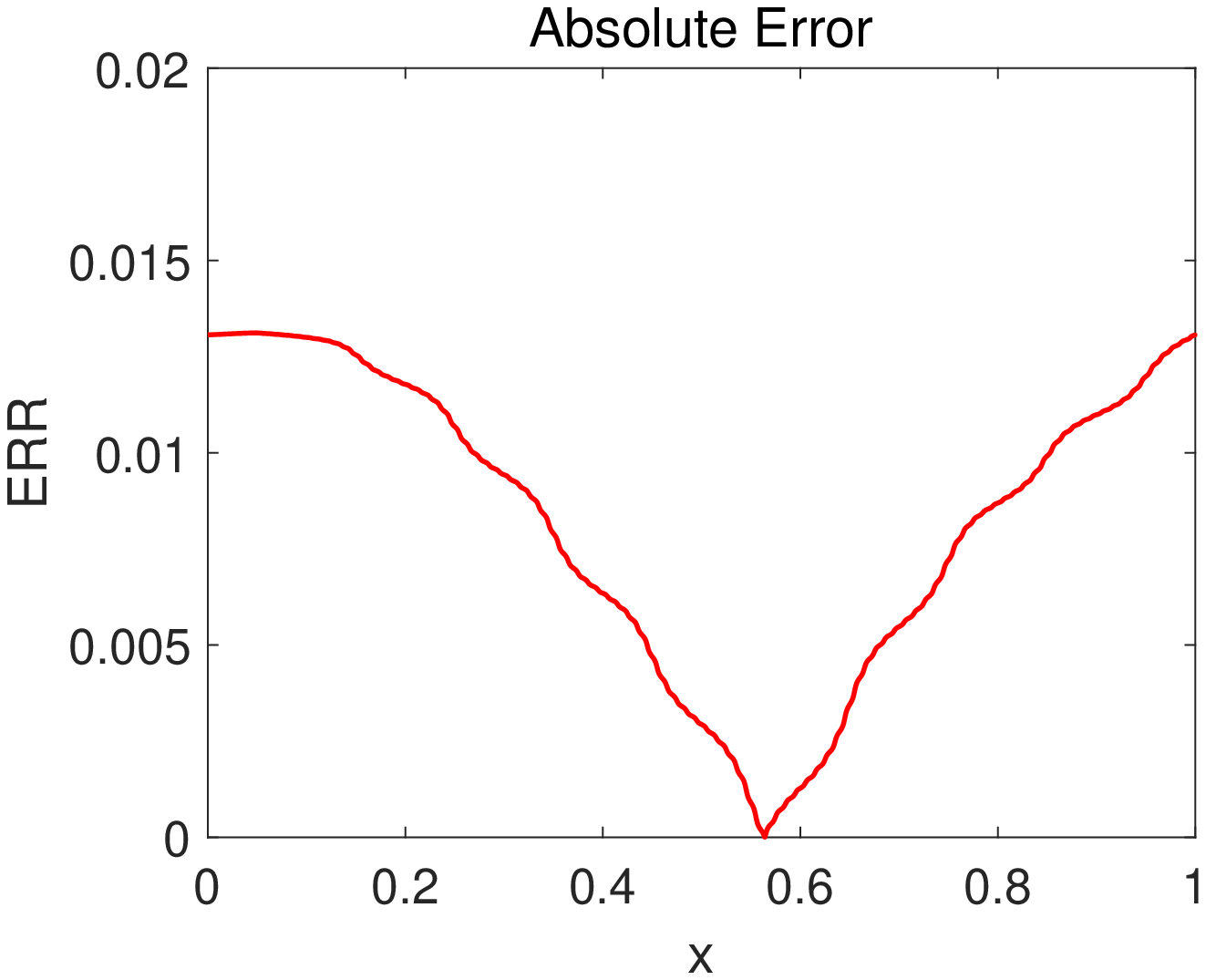}
	}
         \subfigure[point-wise error of LDLM1]{
		\label{3scale2absLDLM1}
		\includegraphics[scale=0.4]{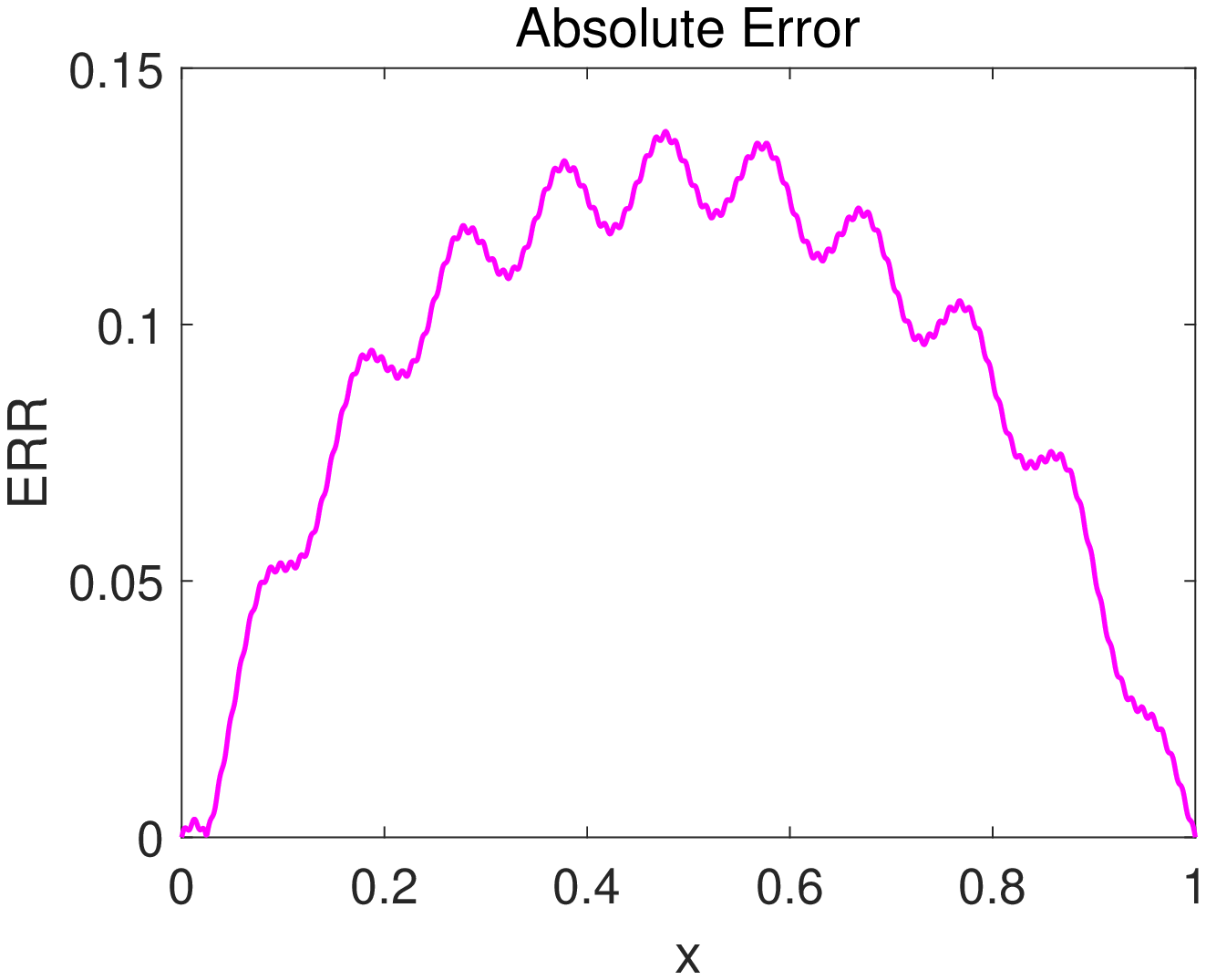}
	}
         \subfigure[point-wise error of LDLM2]{
		\label{3scale2absLDLM2}
		\includegraphics[scale=0.4]{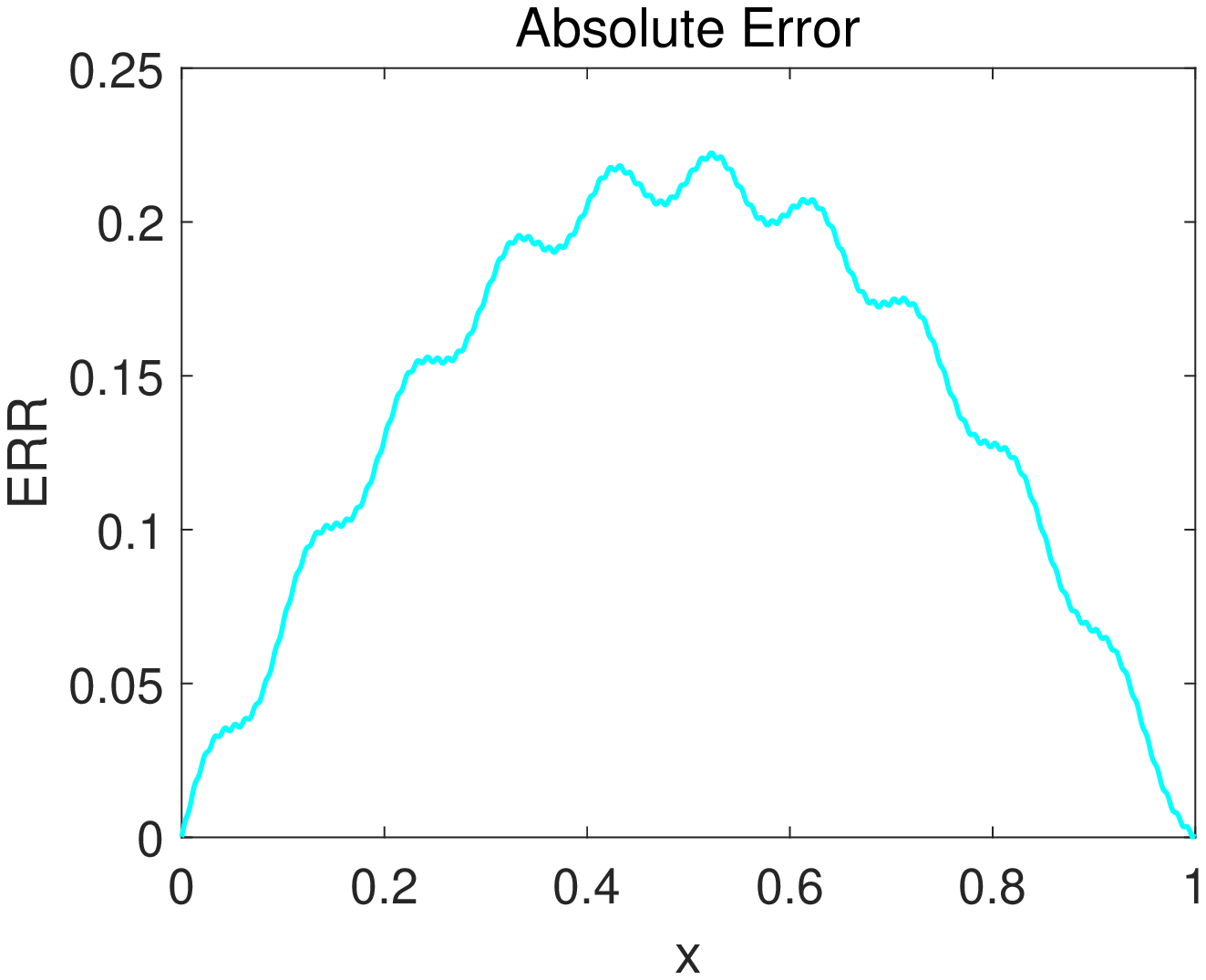}
	}
	\subfigure[REL]{
		\label{3scale_rel}
		\includegraphics[scale=0.38]{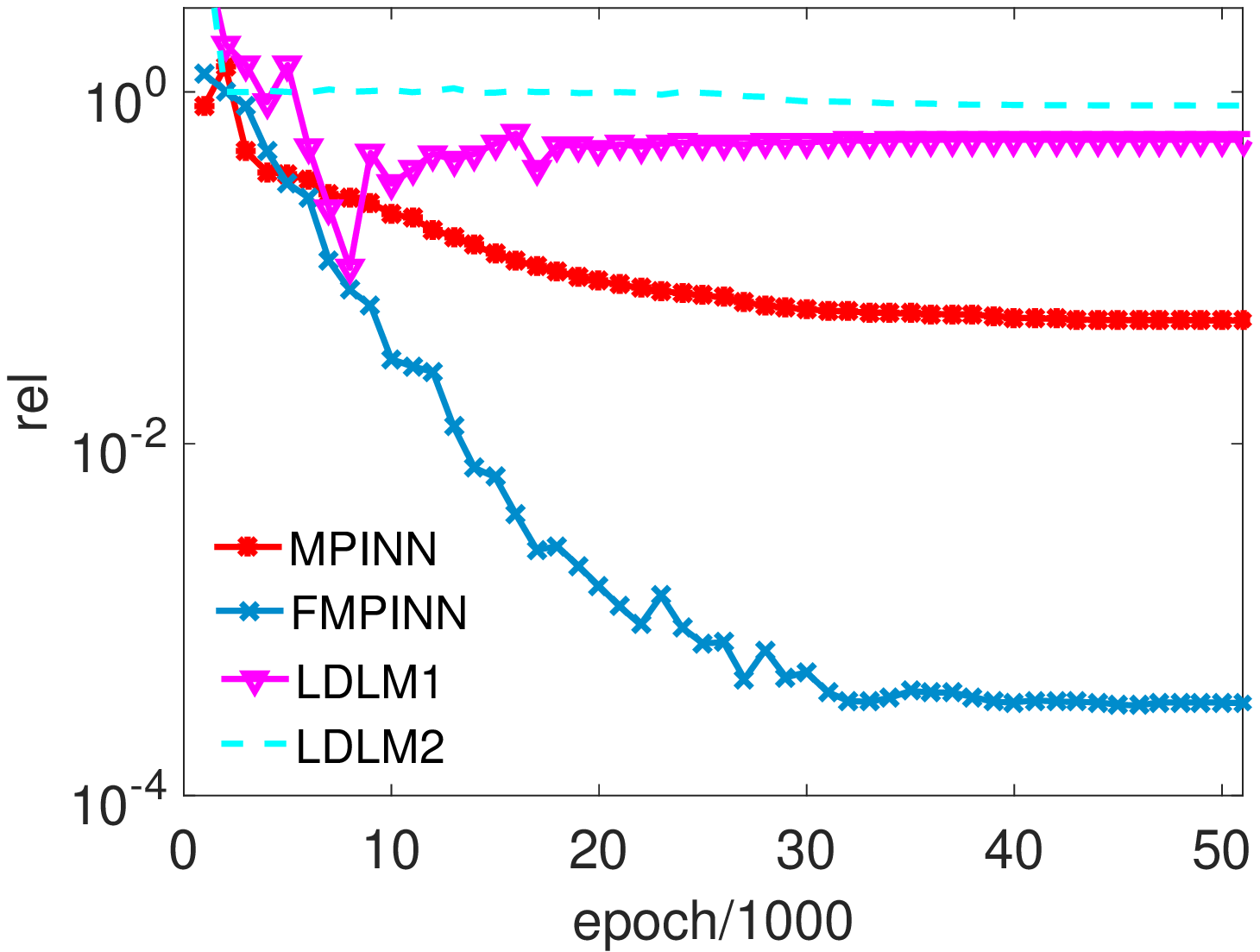}
	}
	\caption{Rough coefficient, loss of flux term and testing results for Example \ref{Diffusion1D3scale} when $\varepsilon_1=0.1$ and $\varepsilon_2=0.01$.}
	\label{3scale_Diffusion}
\end{figure}

We solve the above three scale problem when $\varepsilon_1=0.1$ and $\varepsilon_2=0.01$ by employing the aforementioned FMPINN, MPINN, LDLM1 and LDLM2 models, respectively. Their all settings are same as the Example \ref{DiffusionEq_1d_01}. The training dataset includes 3000 interior random points and 500 boundary random points, and the testing dataset includes 1000 equidistant samples. The related experiment results are listed in Table \ref{results23scale} and plotted in Fig. \ref{3scale_Diffusion}, respectively. 

\begin{table}[H]
	\centering
	\caption{The relative error and consumed time of FMPINN, MPINN, LDLM1, and LDLM2 for Example \ref{Diffusion1D3scale}.}
	\label{results23scale}
	\begin{tabular}{|c|c|c|c|c|}
		\hline
		Method         &FMPINN   &MPINN    &LDLM1    &LDLM2     \\  \hline
		REL            &3.36e-4  &5.02e-2  &0.5341   &0.8372    \\  \hline
		Total time(s)  &696.537  &965.076 &377.18   &399.575   \\  \hline
	\end{tabular}
\end{table}

Fig. \ref{3scale_Diffusion} shows that the FMPINN model still is well able to capture all oscillations of the exact solution for the three-scale problem, the MPINN model also captures the profile of the solution of \eqref{eq:multiscale} with $\varepsilon_1=0.1$ and $\varepsilon_2=0.01$. However, the LDLM1 and LDLM2 all fail to fit the solution. Figs. \ref{3scale2absFMPINN} -- \ref{3scale2absLDLM2} not only show the point-wise errors of FMPINN for major points that are close to zero but also reveal the point-wise error of FMPINN is very smaller than that of the MPINN and the LDLM models are all bad. Additionally, Fig. \ref{3scale_rel} and Table \ref{results23scale} illustrate that the REL of FMPINN is superior to that of MPINN by more than two orders of magnitude, and its running time is 696.537 seconds and less than that of MPINN.

From the above results, we conclude that the FMPINN model is remarkable to address the \eqref{eq:multiscale} with rough coefficient in one-dimensional space, it generally outperforms the MPINN and LDLM models.

\begin{example}\label{Example2d_twoscales}
	 We consider the following two-dimensional problem for \eqref{eq:multiscale}  with Dirichlet boundary in regular domains $\Omega=[-1,1]\times[-1,1]$. In this example, we choose the $f(x_1,x_2)=5$ and provide the following two-scales coefficient with scale separation
	\begin{equation}\label{Diffusion2d_twoscales}
	A^{\varepsilon}(x_1,x_2) =\frac{1.5+\sin(2\pi x_1/\varepsilon)}{1.5+\sin(2\pi x_2/\varepsilon)}+\frac{1.5+\sin(2\pi x_2/\varepsilon)}{1.5+\cos(2\pi x_1/\varepsilon)}+\sin(4x_1^2x_2^2)+1.
	\end{equation}
where  $\varepsilon>0$ is a small parameter such that $\varepsilon^{-1}\in\mathbb{N}^+$. Since the corresponding exact solution can not be expressed explicitly in this example, then a reference solution $u^{\varepsilon}(x_1,x_2)$ is set as the finite element solution computed by numerical homogenization method \cite{owhadi2014polyharmonic}  on a square grid $[-1,1]\times[-1,1]$ with mesh-size $h=1/128$.
\end{example}

\begin{figure}
	\centering 
	\subfigure[Rough coefficient $A^{\varepsilon}$]{
		\label{2d_twoscales:a}     
		\includegraphics[scale=0.425]{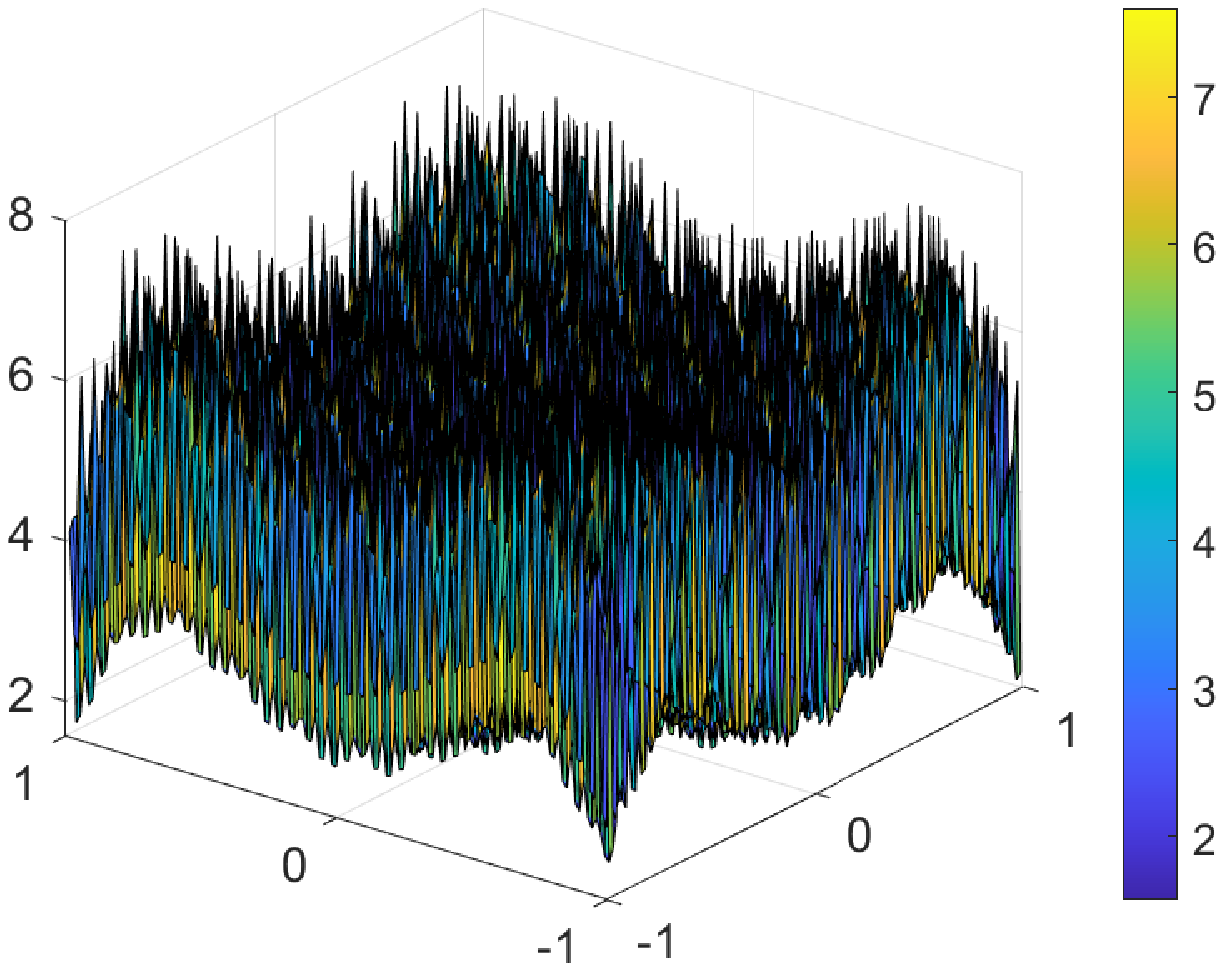}  
	}    
	\subfigure[reference solution]{
		\label{2d_twoscales:solu}     
		\includegraphics[scale=0.425]{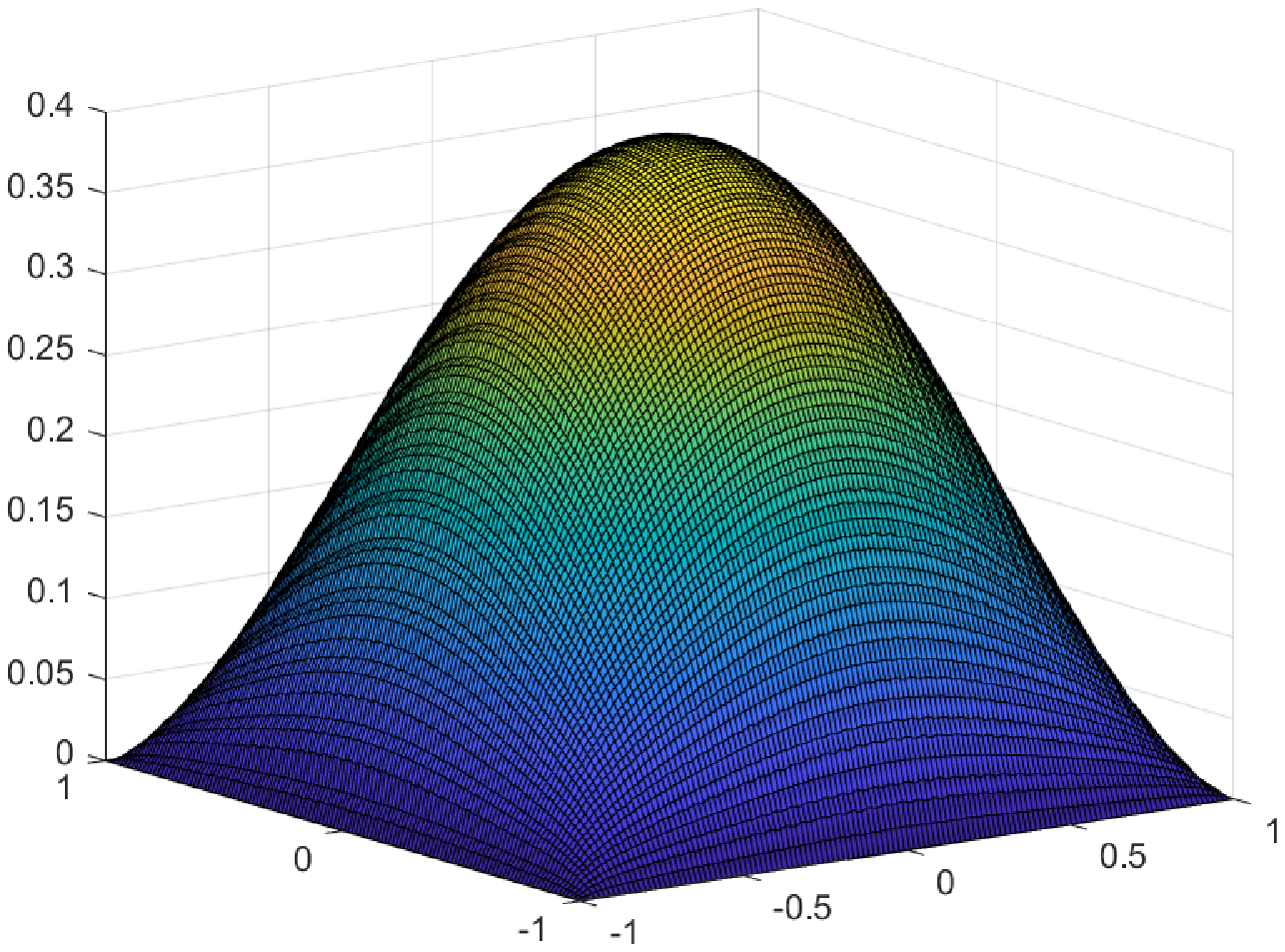}  
	}  
	\subfigure[loss of flux term for FMPINN]{
		\label{2d_twoscales:loss2flux}     
		\includegraphics[scale=0.425]{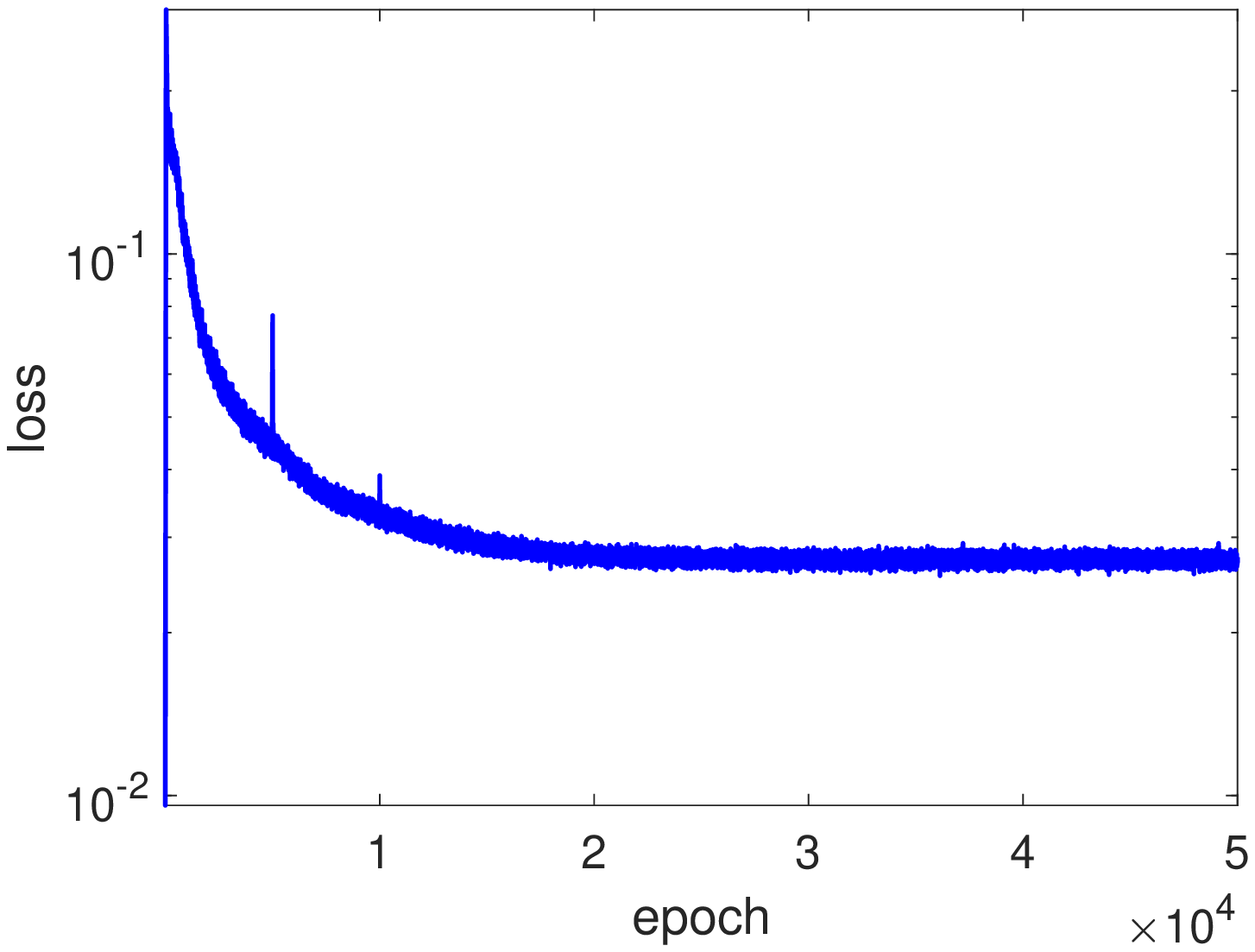}  
	}        
	\subfigure[point-wise error of FMPINN]{ 
		\label{2d_twoscales:b}     
		\includegraphics[scale=0.425]{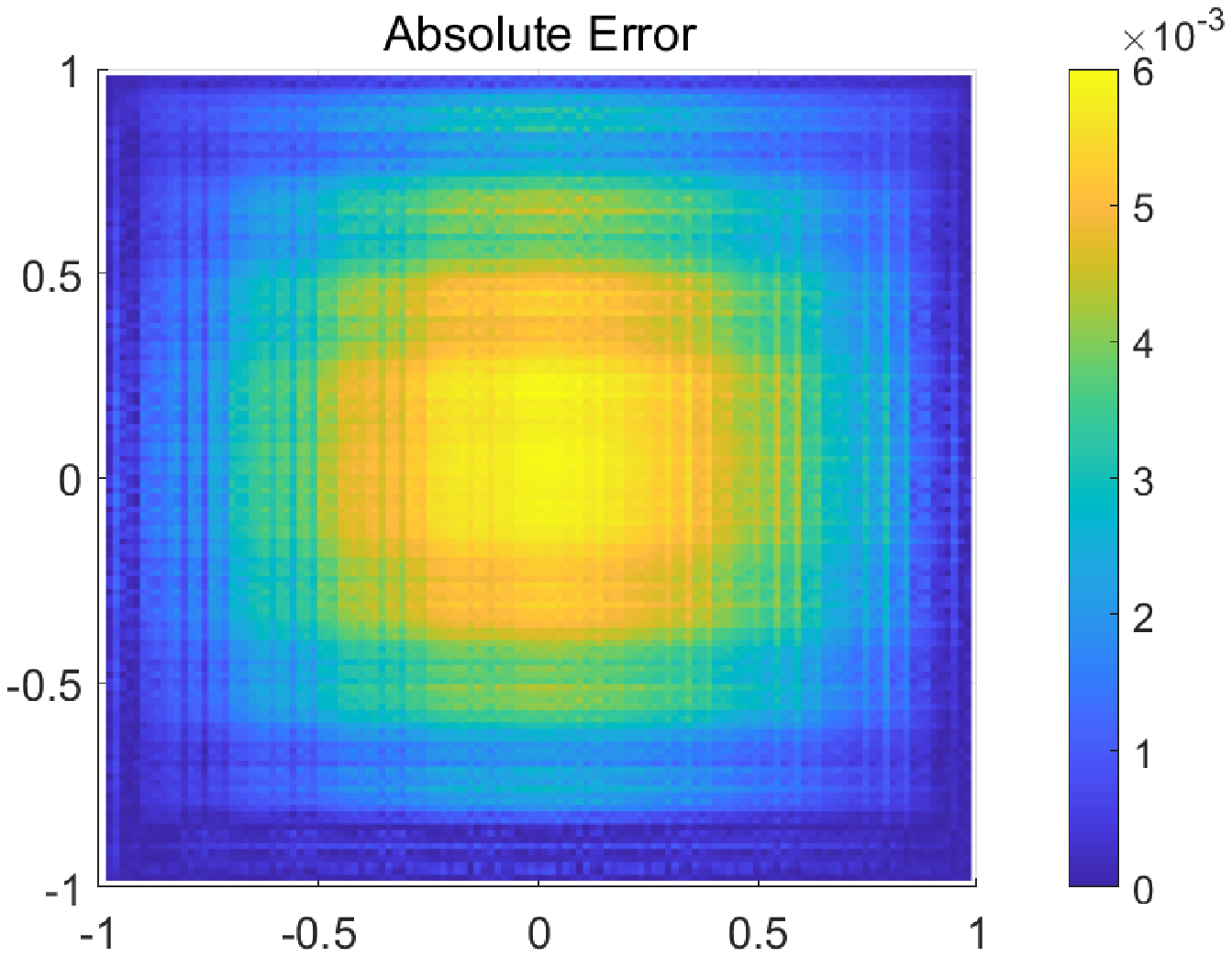}     
	} 
	\subfigure[point-wise error of MPINN]{ 
		\label{2d_twoscales_PERR2MPINN}     
		\includegraphics[scale=0.425]{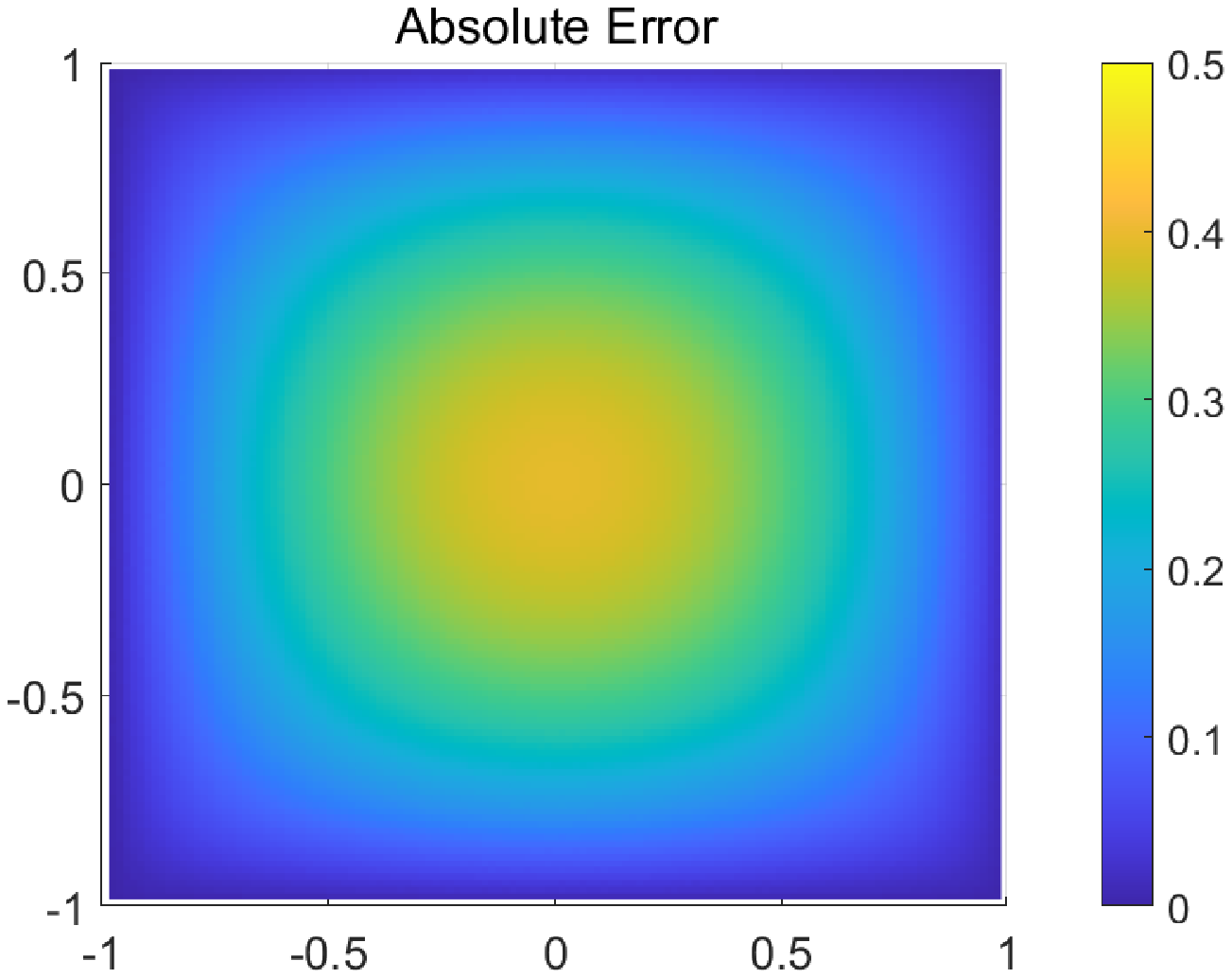}     
	}   
 	\subfigure[point-wise error of LDLM1]{ 
		\label{2d_twoscales_PERR2LDLM1}     
		\includegraphics[scale=0.425]{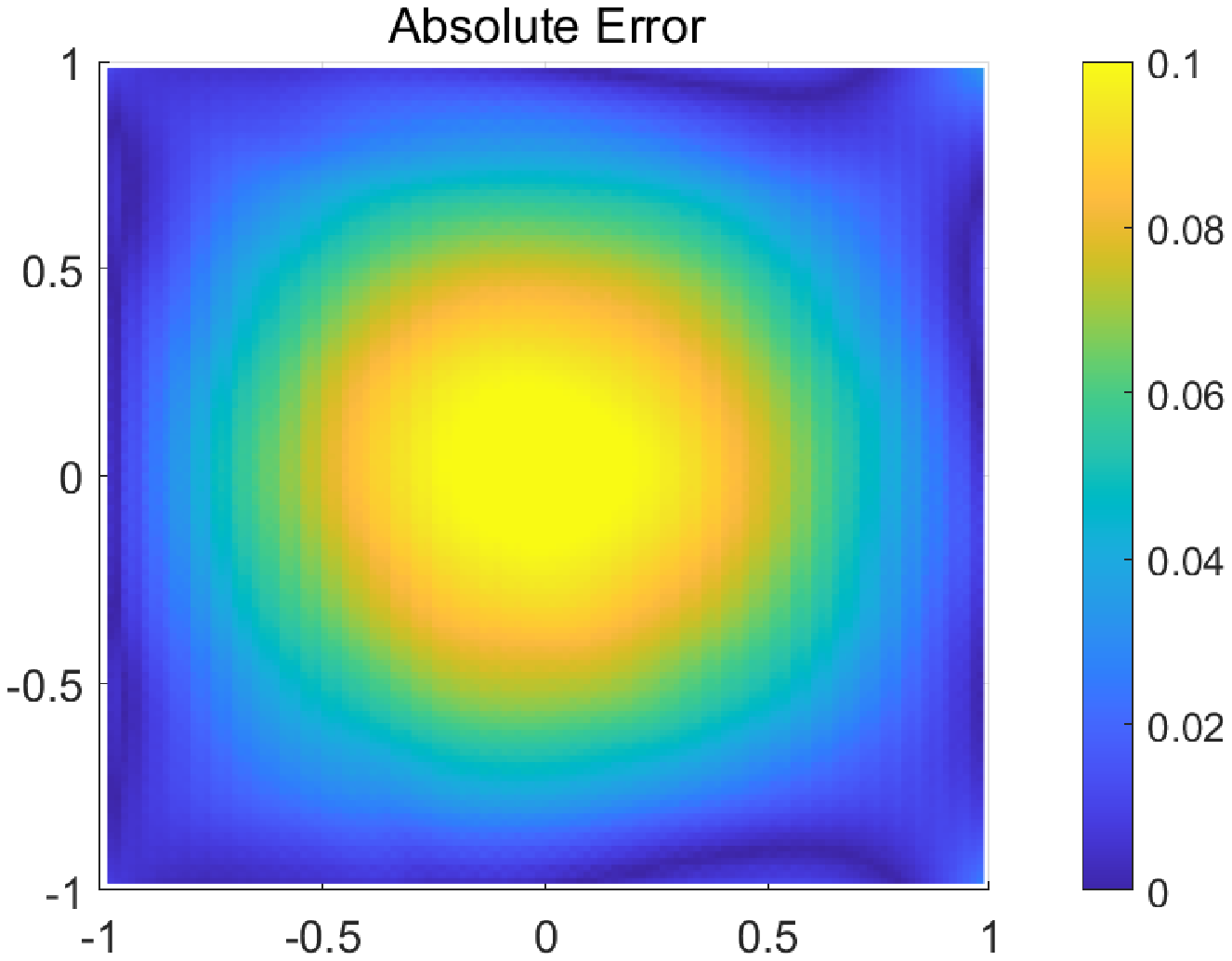}     
	}  
 	\subfigure[point-wise error of LDLM2]{ 
		\label{2d_twoscales_PERR2LDLM2}     
		\includegraphics[scale=0.425]{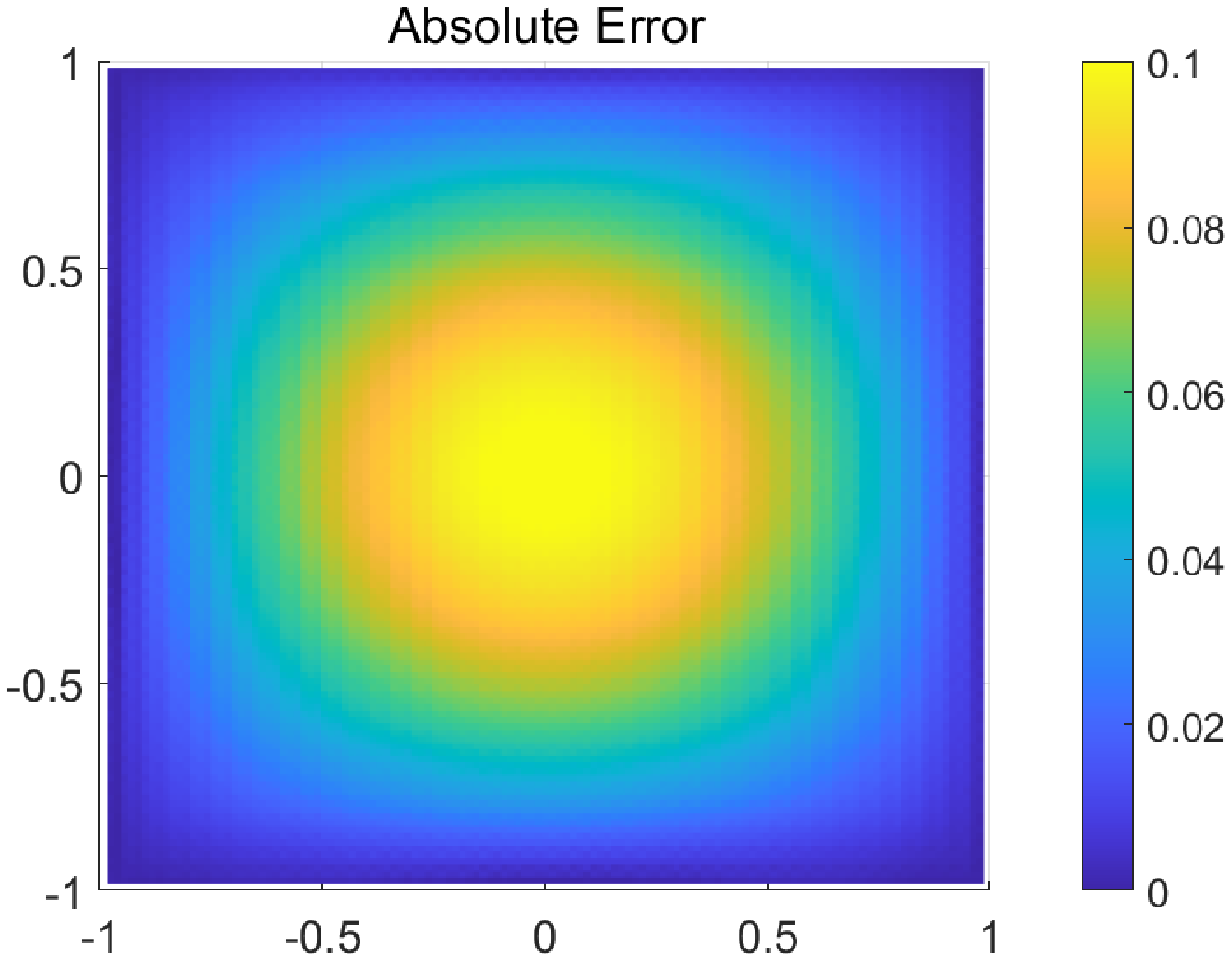}     
	}  
	\subfigure[REL] { 
		\label{2d_twoscales:e}     
		\includegraphics[scale=0.425]{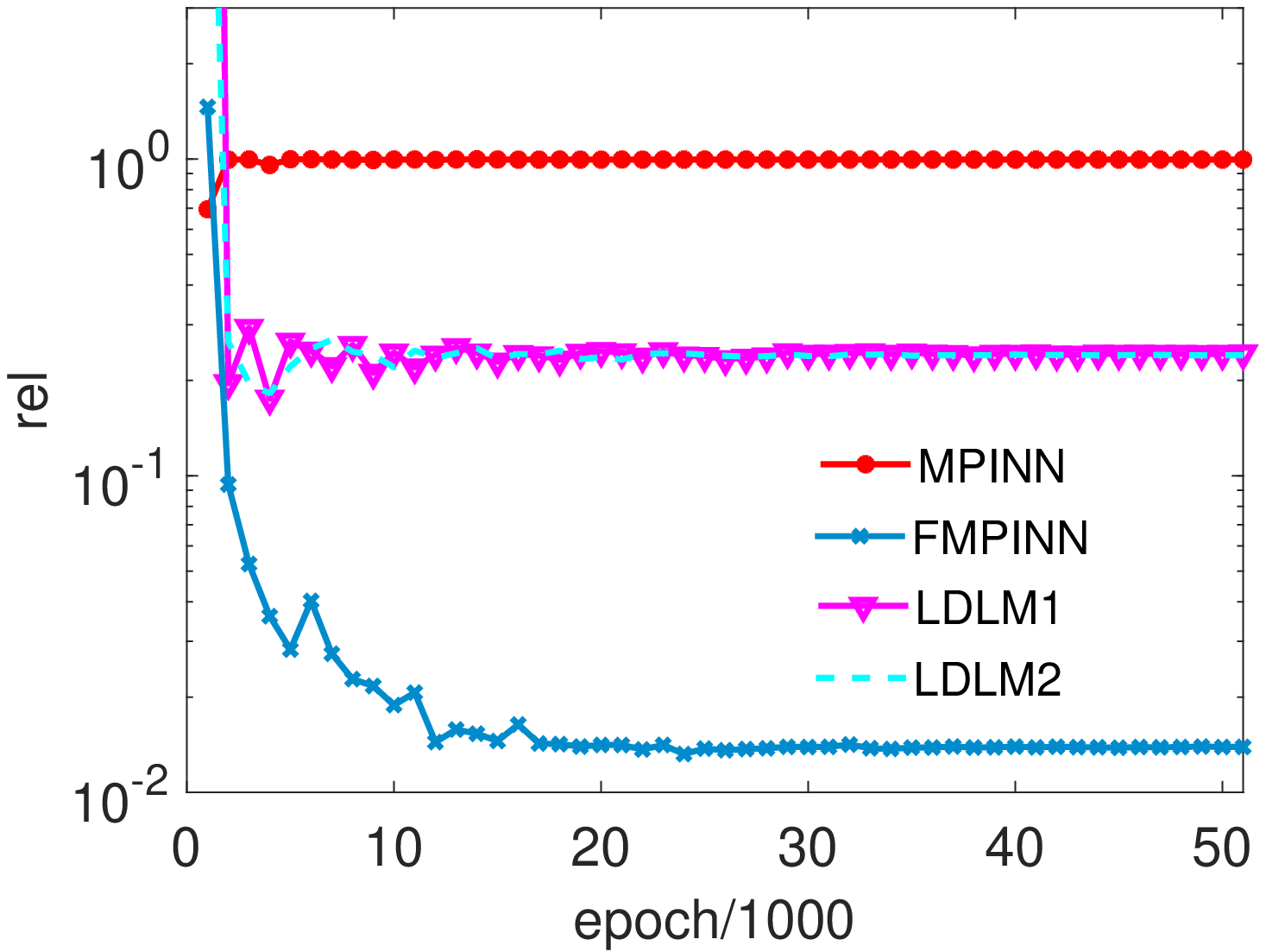}     
	}
\caption{Rough coefficient, reference solution, loss of flux term and testing results for Example \ref{Example2d_twoscales}} 
	\label{2d_twoscales}         
\end{figure}

We solve the above two scale problem when $\varepsilon=0.05$ by employing the aforementioned FMPINN, MPINN, LDLM1 and LDLM2 models, respectively.
The size of hidden layer for each subnetwork of FMPINN and MPINN is set as (40, 60, 40,40,40) and the hidden layers' size for LDLMs is set as $(400, 250, 250, 200, 200)$. At each training step, the training dataset includes 5000 points sampled inside the $\Omega$ and 2000 boundary points sampled from the $\partial\Omega$, respectively. In order to test our models, the testing dataset is the collection of all grid points in domain $[-1,1]\times[-1,1]$ with mesh-size $h=1/128$. The related experiment results are listed in Table \ref{results2d_twoscales} and plotted in Fig. \ref{2d_twoscales}, respectively. 

\begin{table}[H]
	\centering
	\caption{The relative error and consumed time of FMPINN, MPINN, LDLM1 and LDLM2 for Example \ref{Example2d_twoscales}.}
	\label{results2d_twoscales}
	\begin{tabular}{|c|c|c|c|c|}
		\hline
		Method         &FMPINN     &MPINN    &LDLM1   &LDLM2    \\  \hline
		REL            &0.0139     &0.99     &0.2431  &0.2401   \\  \hline
		Total time(s)  &2098.258   &3885.934 &626.685 &689.619  \\  \hline
	\end{tabular}
\end{table}

In this example, the $A^{\varepsilon}(x_1,x_2)$ have two different frequency components and is quite oscillating(seeing Fig. \ref{2d_twoscales:a}), then DNN will encounter some troubles to address multi-scale PDEs\eqref{eq:multiscale}.  According to the results of point-wise error (Figs.\ref{2d_twoscales:b} - \ref{2d_twoscales_PERR2LDLM2}) and relative errors(Fig.\ref{2d_twoscales:e}), the performance of our FMPINN model is still superior to that of  the MPINN, LDLM1 and LDLM2 models, and can obtain a favorable approximation to multi-scale problems\eqref{eq:multiscale}. In addition, the test REL curve in Fig. \ref{2d_twoscales:e} indicates the FMPINN model is stable in the whole training cycle and its tendency is consistent with the curve of loss for flux term in Fig. \ref{2d_twoscales:loss2flux}. Clearly, the running time of our FMPINN model is about half of that of the MPINN model, which means the FMPINN model is efficient in solving multi-scale PDEs\eqref{eq:multiscale} with two scales coefficient.

\begin{example}\label{Example2d_01}
	 We consider the following two-dimensional problem for \eqref{eq:multiscale} with Dirichlet boundary  in regular domains $\Omega=[-1,1]\times[-1,1]$. In this example, we choose the $f(x_1,x_2)=1$ and provide a multi-frequency coefficient
	\begin{equation}\label{Diffusion2d_Aeps}
	A^{\varepsilon}(x_1,x_2) =\Pi_{i=1}^{5} \bigg{(}1+0.5\cos\left(2^i\pi(x_1+x_2)\right)\bigg{)}\bigg{(}1+0.5\sin\left(2^i\pi(x_2-3x_1)\right)\bigg{)}.
	\end{equation}
Same as the Example \ref{Example2d_twoscales}, a reference solution $u^{\varepsilon}(x_1,x_2)$ is set as the finite element solution computed by numerical homogenization method \cite{owhadi2014polyharmonic} on a square grid $[-1,1]\times[-1,1]$ with mesh-size $h=1/128$. 
\end{example}

\begin{figure}[H]
	\centering 
	\subfigure[Rough coefficient $A^{\varepsilon}$]{
		\label{2d:a}     
		\includegraphics[scale=0.425]{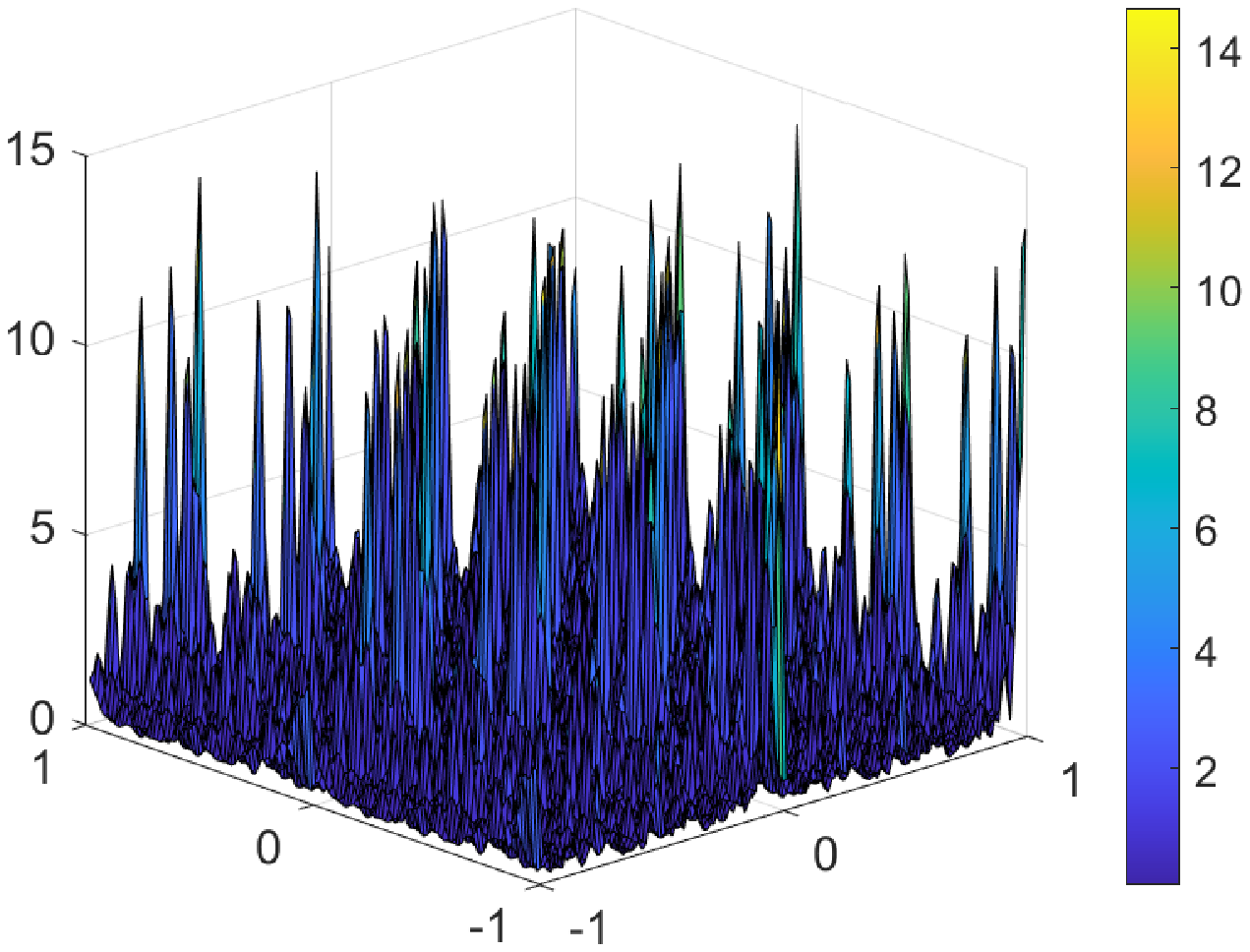}  
	}    
	\subfigure[reference solution]{
		\label{2d:solu}     
		\includegraphics[scale=0.425]{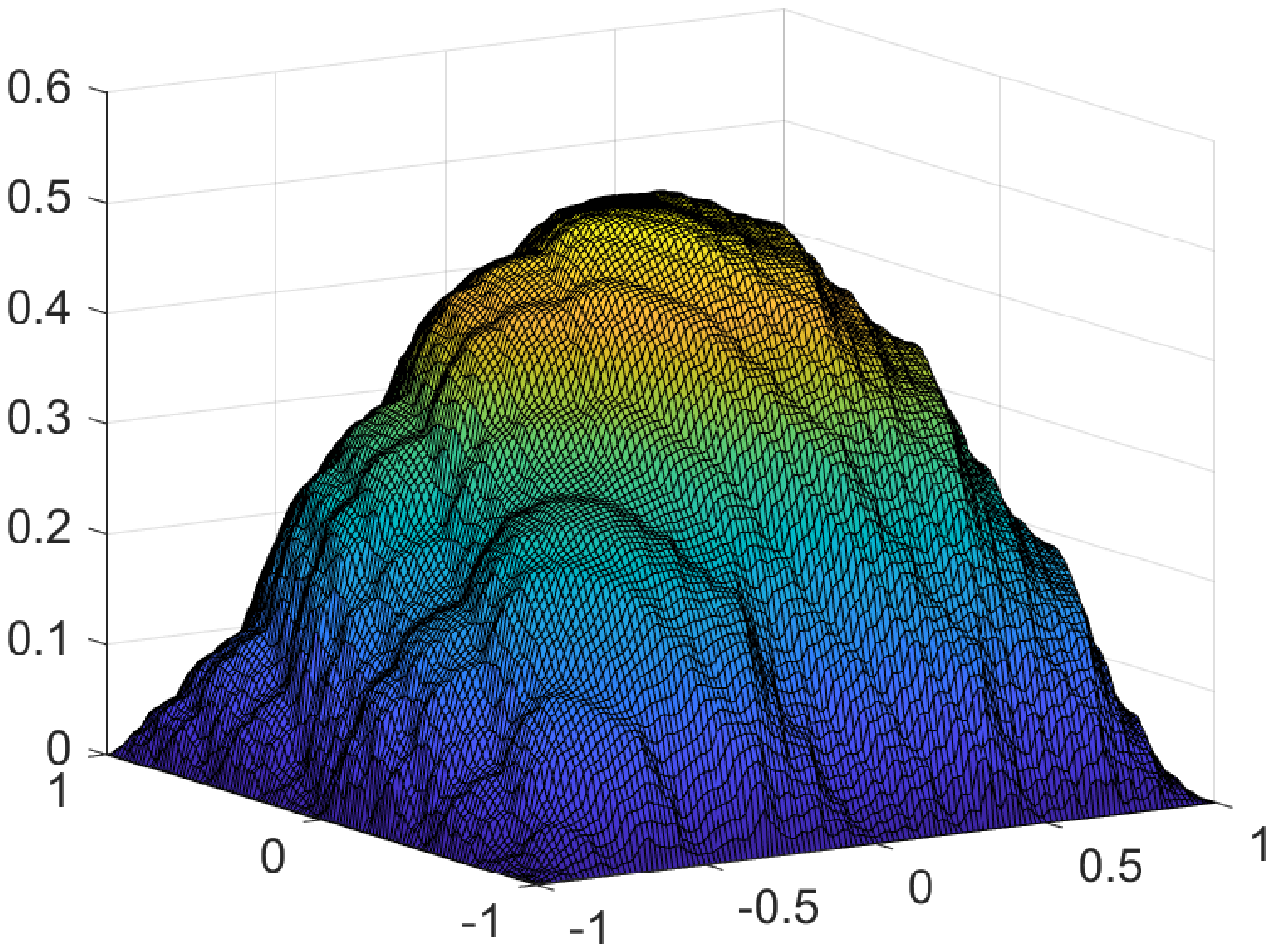}  
	}  
	\subfigure[loss of flux term for FMPINN]{
		\label{2d:loss2flux}     
		\includegraphics[scale=0.425]{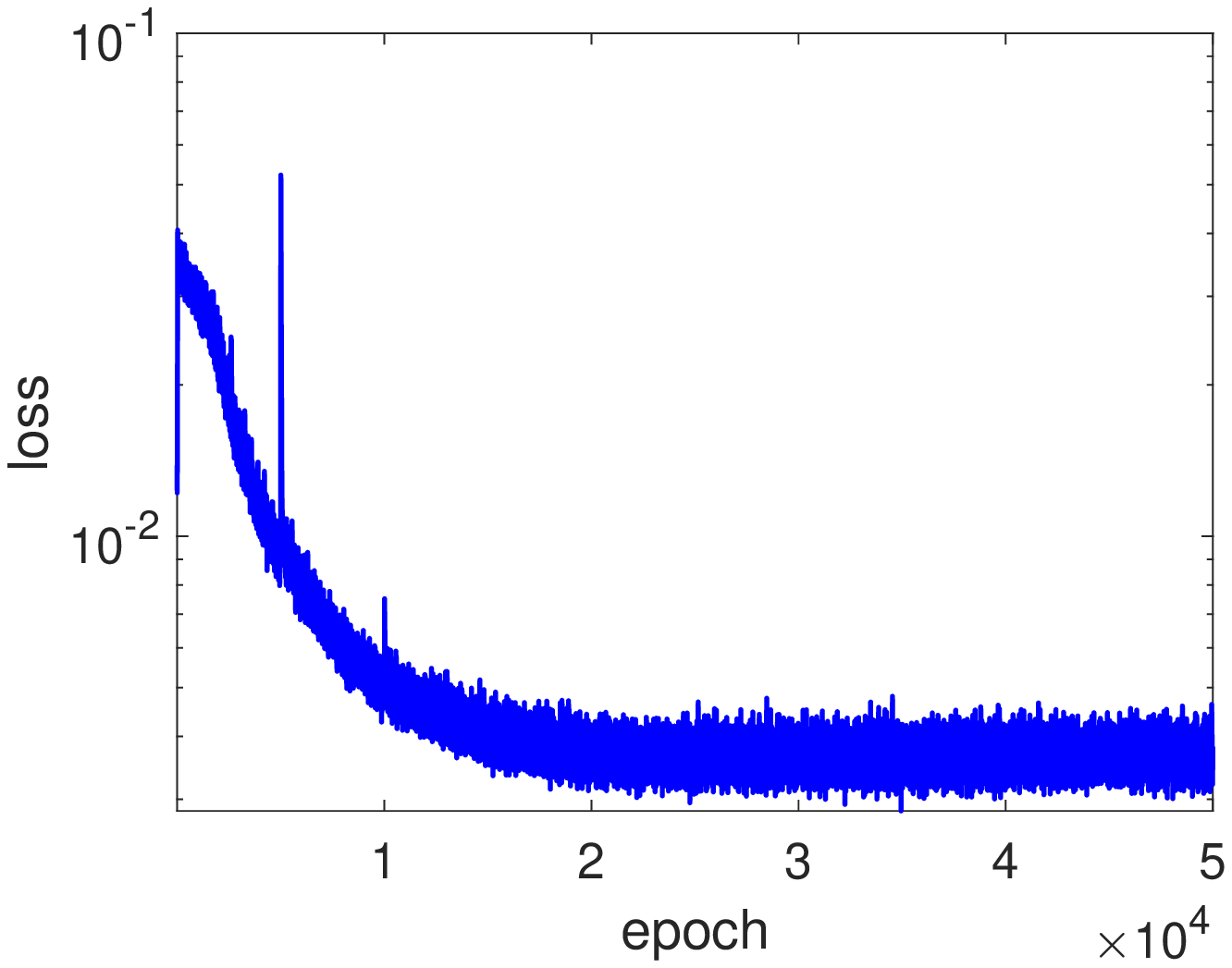}  
	}        
	\subfigure[point-wise error of FMPINN]{ 
		\label{2d:b}     
		\includegraphics[scale=0.425]{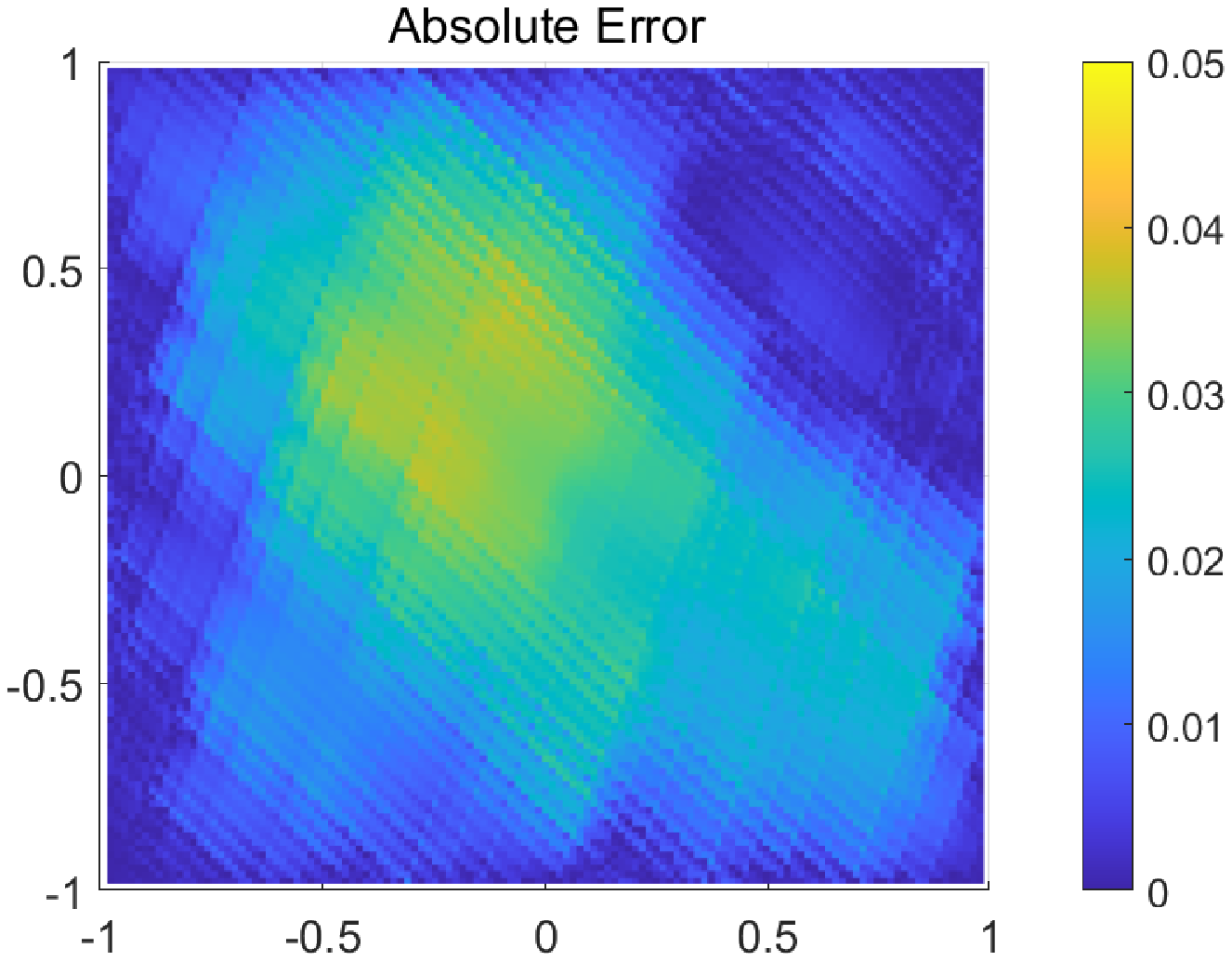}     
	} 
	\subfigure[point-wise error of MPINN]{ 
		\label{2dPERR2MPINN}     
		\includegraphics[scale=0.425]{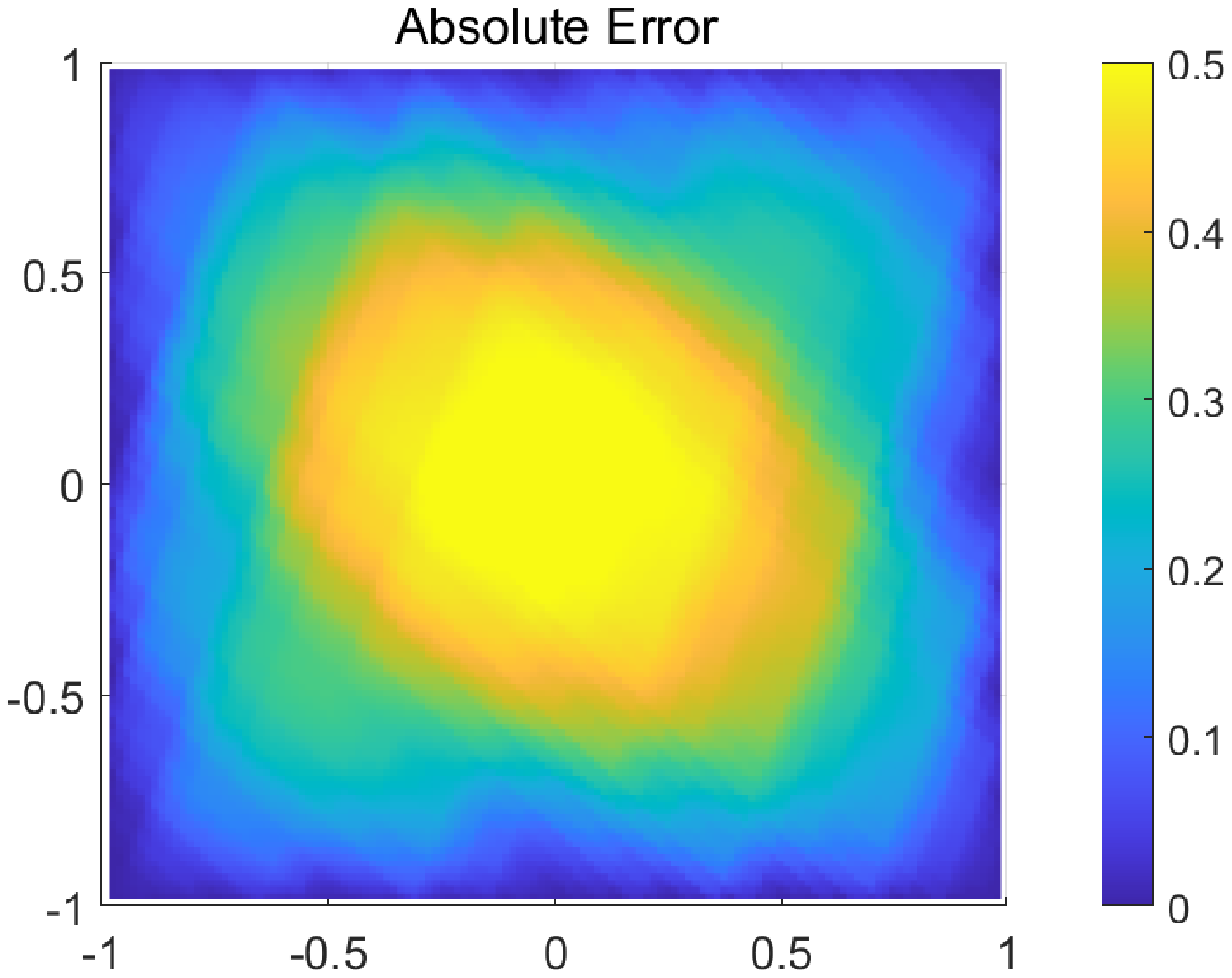}     
	}   
 	\subfigure[point-wise error of LDLM1]{ 
		\label{2dPERR2LDLM1}     
		\includegraphics[scale=0.425]{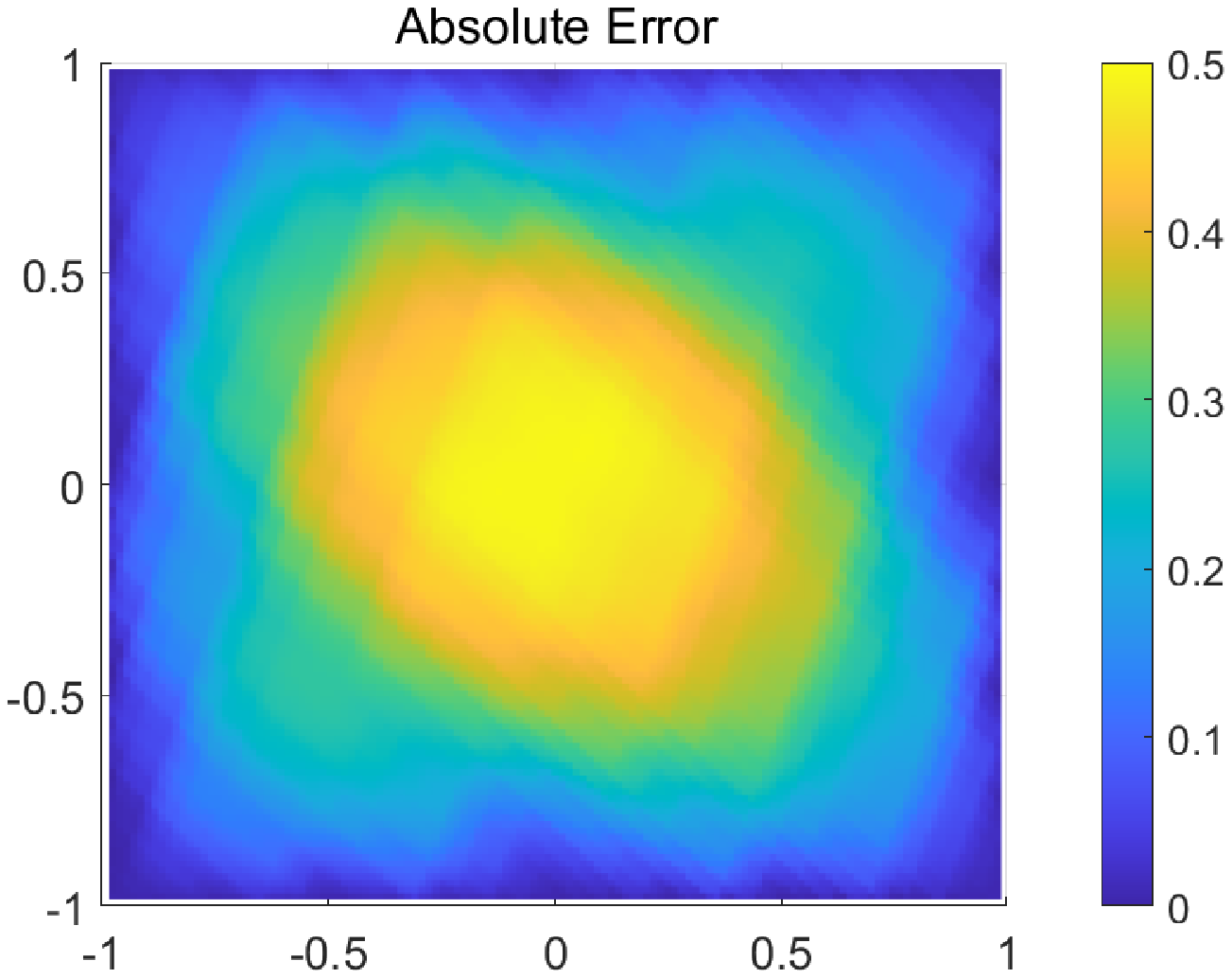}     
	}  
 	\subfigure[point-wise error of LDLM2]{ 
		\label{2dPERR2LDLM2}     
		\includegraphics[scale=0.425]{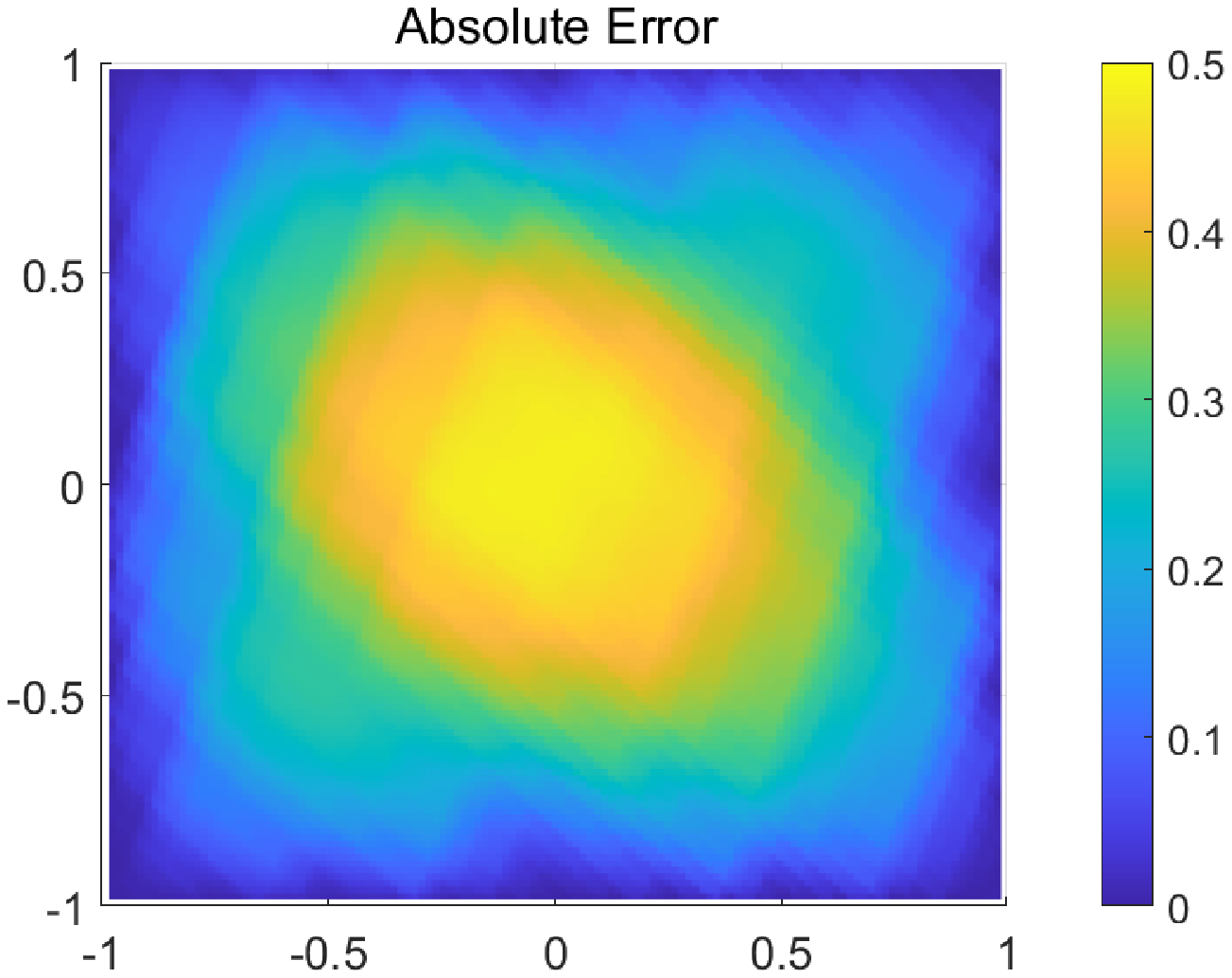}     
	}  
	\subfigure[REL] { 
		\label{2dPDE1q6:e}     
		\includegraphics[scale=0.425]{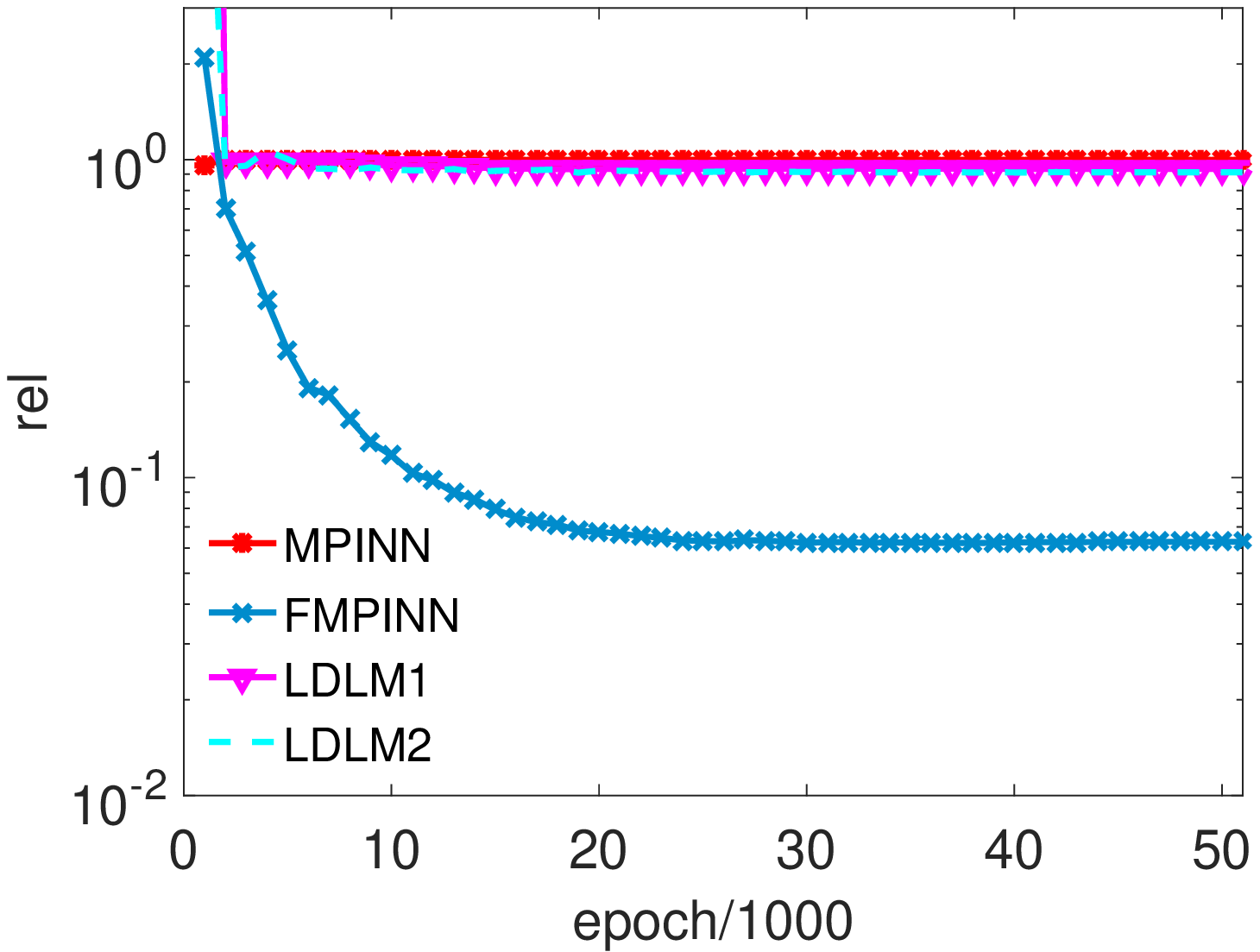}     
	}
\caption{Rough coefficient, reference solution, loss of flux term and testing results for Example \ref{Example2d_01}} 
	\label{2dPDE3q6}         
\end{figure}

By meticulously implementing the previously mentioned FMPINN, MPINN, and LDLMs models with the specified setups, we obtain the approximated solution of \eqref{eq:multiscale} with \eqref{Diffusion2d_Aeps}. The setup for all models is identical to that of Example  \ref{Example2d_twoscales}. During each training step, the training dataset comprises 5000 points randomly sampled from $\Omega$ and 2000 boundary points sampled from the boundary $\partial\Omega$, respectively. Meantime, the testing dataset composes of grid points on the square domain $[-1,1]\times[-1,1]$ with mesh-size $h=1/128$. The related experiment results are listed in Table \ref{results2d} and plotted in Fig. \ref{2dPDE3q6}, respectively. 

\begin{table}[H]
	\centering
	\caption{The relative error and consumed time of FMPINN, MPINN, LDLM1 and LDLM2 for Example \ref{Example2d_01}.}
	\label{results2d}
	\begin{tabular}{|c|c|c|c|c|}
		\hline
		Method         &FMPINN     &MPINN    &LDLM1   &LDLM2    \\  \hline
		REL            &0.0628     &0.99     &0.936   &0.9127   \\  \hline
		Total time(s)  &2013.258   &3985.934 &606.685 &659.619  \\  \hline
	\end{tabular}
\end{table}

In this example, the $A^{\varepsilon}(x_1,x_2)$ is obviously oscillating with six different frequency components(seeing Fig. \ref{2d:a}), it will increase the difficulty for DNN to address multi-scale PDEs\eqref{eq:multiscale}.  The point-wise error (Figs.\ref{2d:b} - \ref{2dPERR2LDLM2}) and the relative errors(Fig.\ref{2dPDE1q6:e}) indicate that our FMPINN model is still favorable to capture the solution of multi-scale problems with complex multi-frequency coefficient, but the MPINN, LDLM1 and LDLM2 models all performance poorly  for approximating the solution of \eqref{eq:multiscale}. Additionally, the test REL curve in Fig. \ref{2dPDE1q6:e} and the curve of loss for flux term in Fig. \ref{2d:loss2flux} are all flat indicates the FMPINN model is stable in the whole training cycle. Moreover, the running time of our FMPINN model is less than that of the MPINN model in solving multi-scale PDEs\eqref{eq:multiscale} for coefficient \eqref{Diffusion2d_Aeps}.

\begin{example}\label{E3d}
We next study the performance of our FMPINN model to solve the elliptic equation \eqref{eq:multiscale} with Dirichlet boundary in a cubic domain $\Omega=[0,1]\times[0,1]\times[0,1]$. In which, we take
 \begin{equation}\label{coeff_3D}
A^{\varepsilon}(x_1,x_2,x_{3})=2 + \sin\left(\frac{2\pi x_1}{\varepsilon}\right)\sin\left(\frac{2\pi x_2}{\varepsilon}\right)\sin\left(\frac{2\pi x_3}{\varepsilon}\right).
\end{equation}
with a small parameter $\varepsilon>0$ such that $\varepsilon^{-1}\in\mathbb{N}^+$. Also, we let the force side $f(x_1, x_2, x_3)=20$ and the boundary function $g(x_1,x_2, x_3)=0$ on $\partial \Omega$. 
\end{example}

\begin{figure}
    \centering 
    \subfigure[Rough coefficient $A^{\varepsilon}$ for $z=0.3125$]{
        \label{3d:a}     
        \includegraphics[scale=0.4]{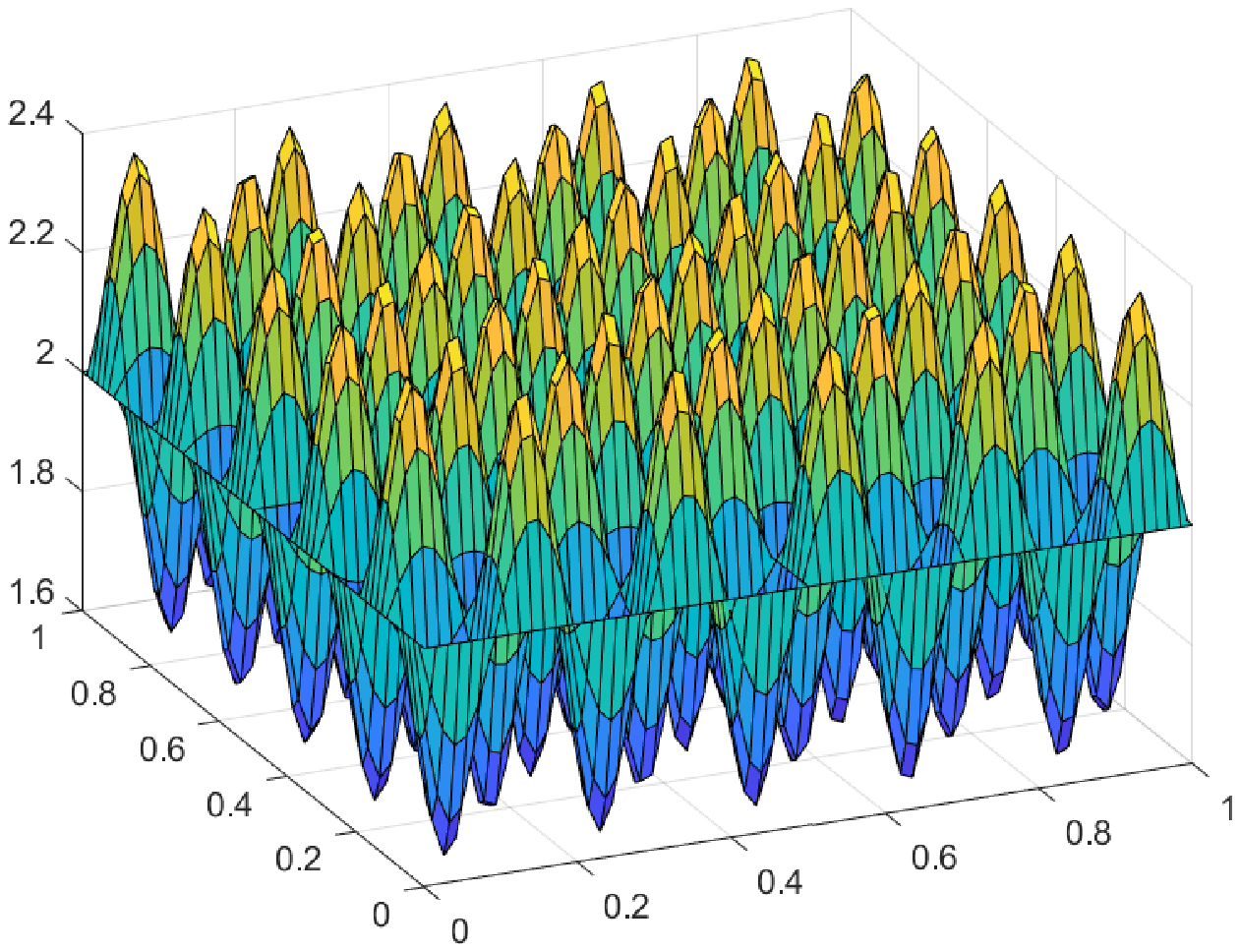}  
    }    
    \subfigure[Reference solution for $z=0.3125$]{
        \label{3d:solu}     
        \includegraphics[scale=0.4]{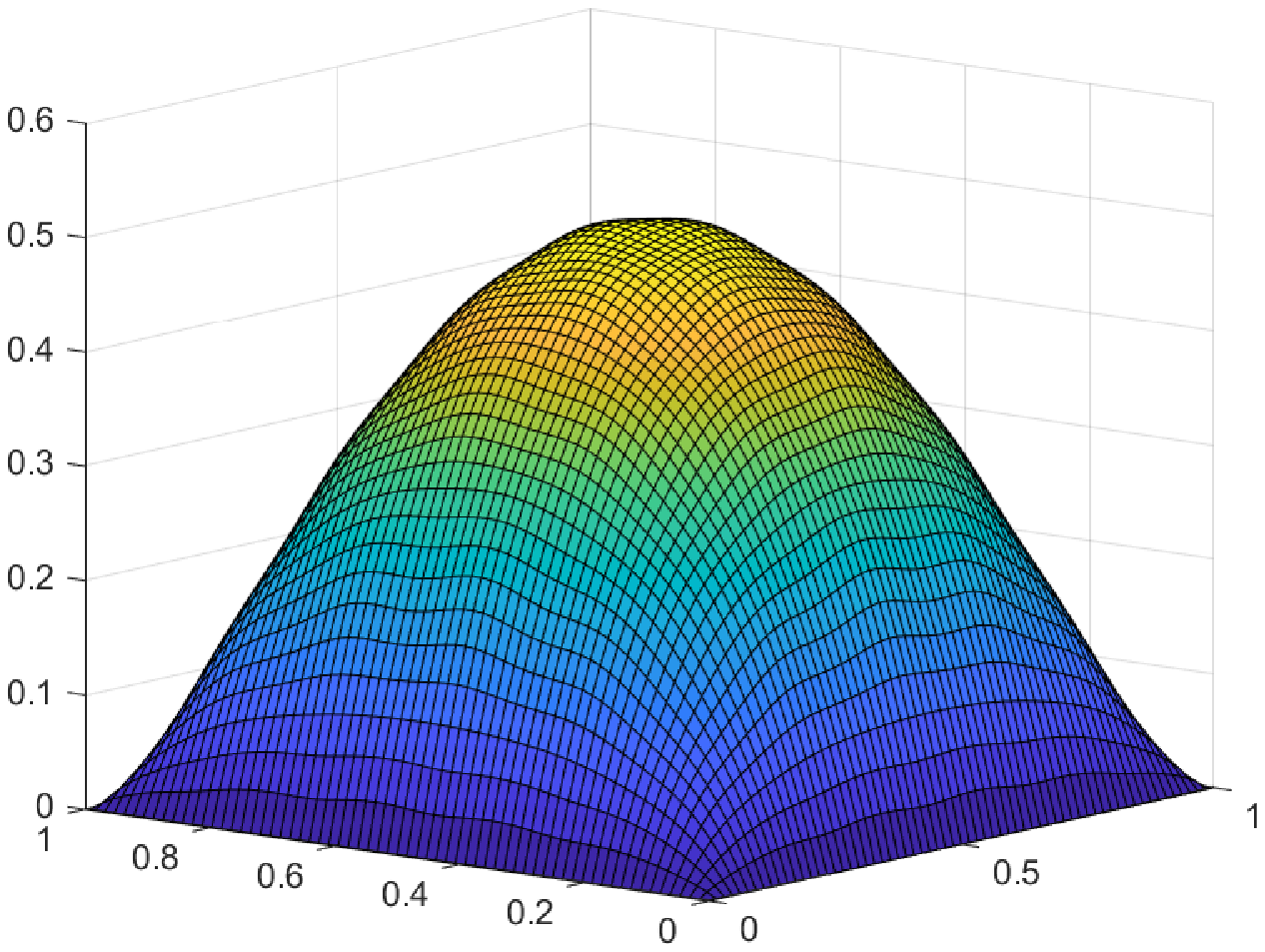}  
    }  
    \subfigure[loss of flux term for FMPINN]{
        \label{3d:loss2flux}     
        \includegraphics[scale=0.4]{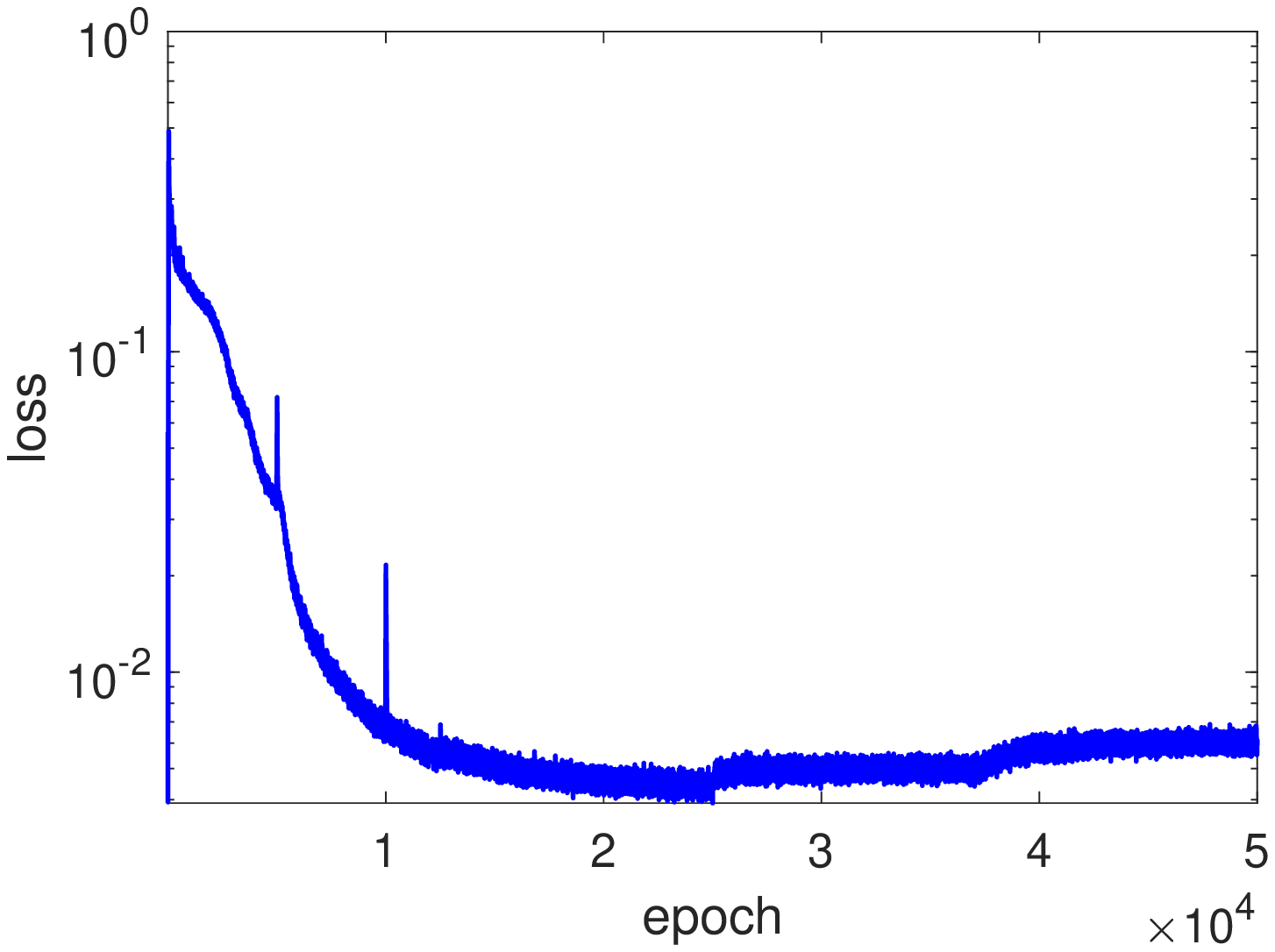}  
    }        
    \subfigure[point-wise error of FMPINN]{ 
        \label{3d:b}     
        \includegraphics[scale=0.4]{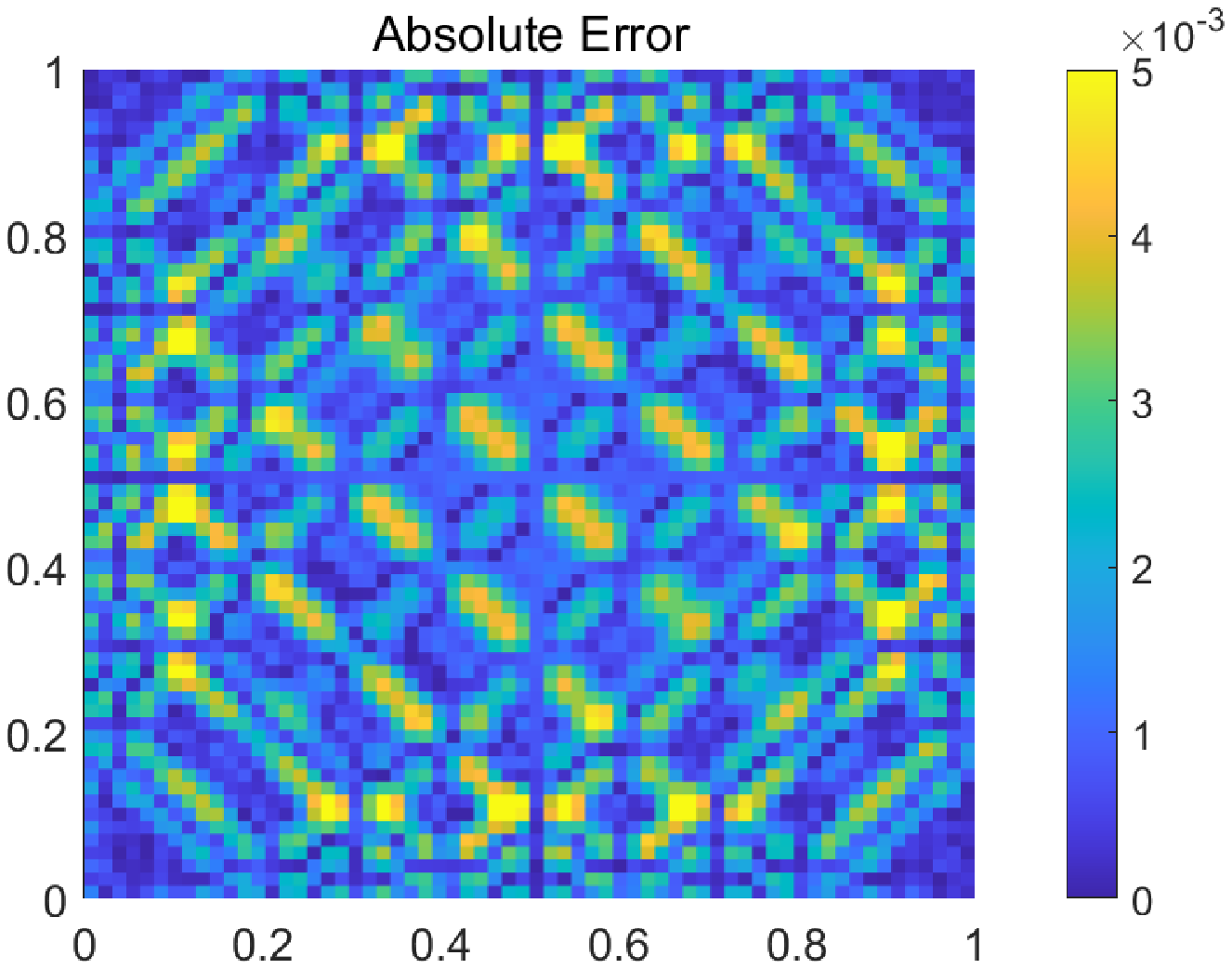}     
    } 
    \subfigure[point-wise error of MPINN]{ 
        \label{3dPERR2MPINN}     
        \includegraphics[scale=0.4]{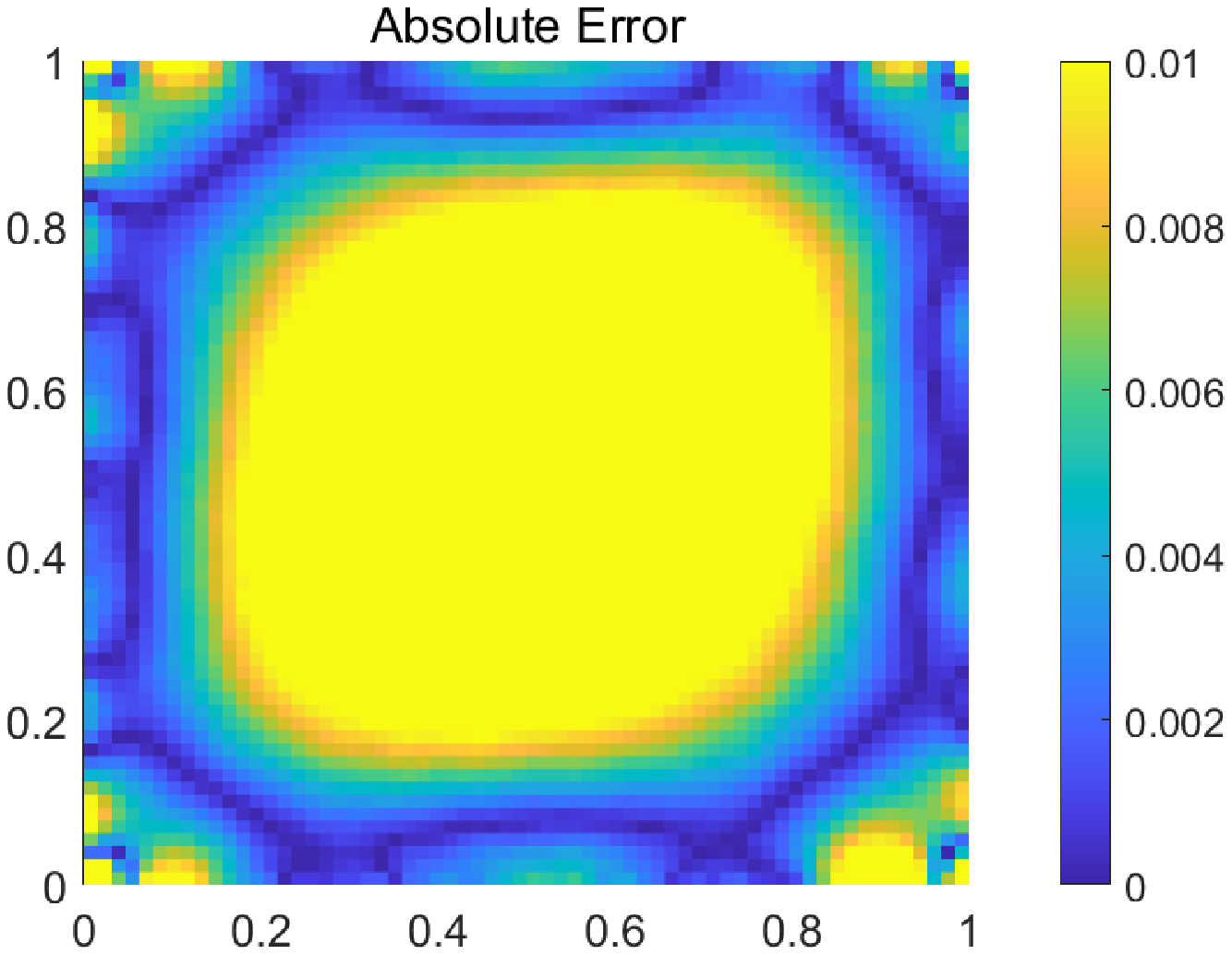}     
    }   
    \subfigure[point-wise error of LDLM1]{ 
        \label{3dPERR2LDLM1}     
        \includegraphics[scale=0.4]{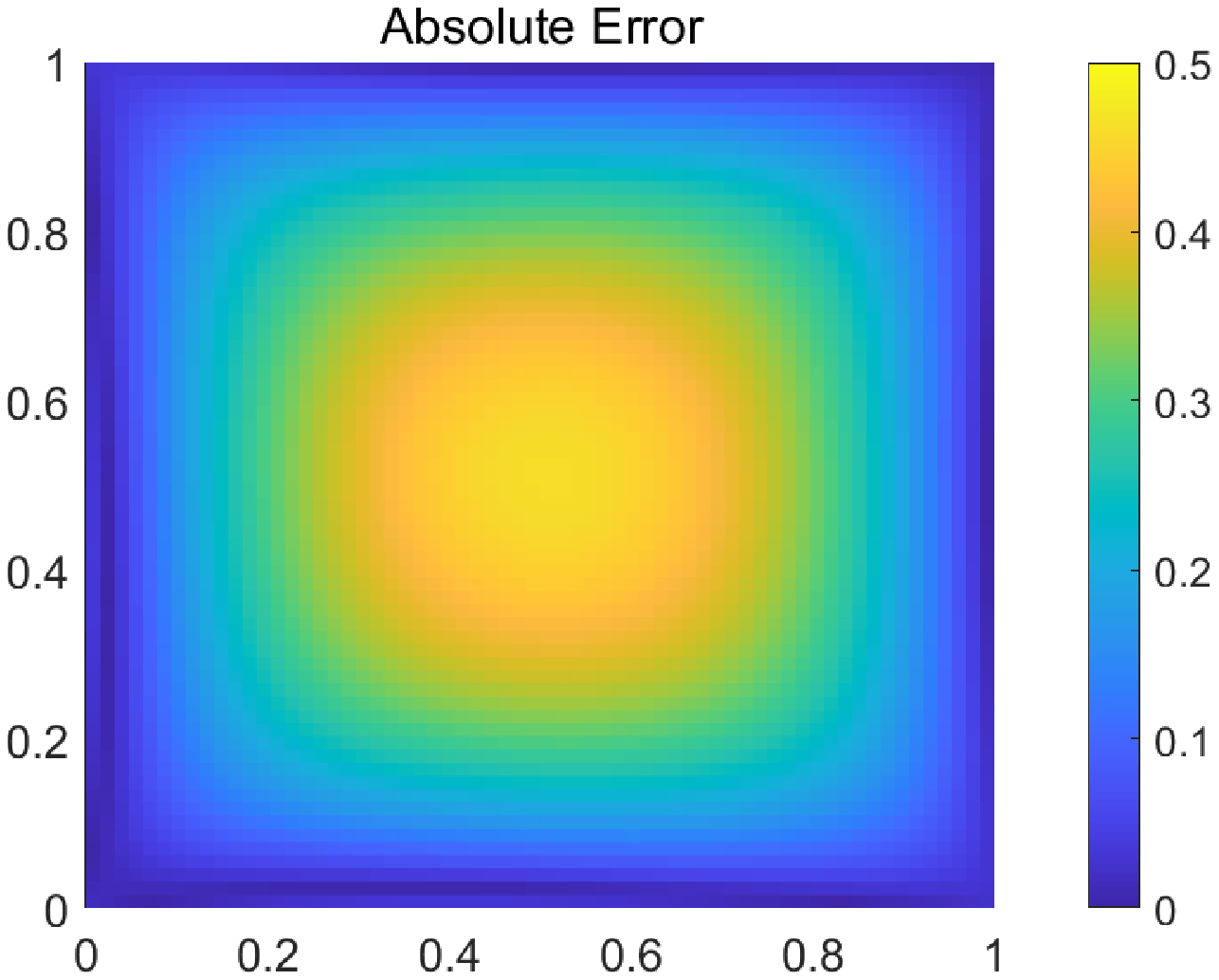}     
    } 
    \subfigure[point-wise error of LDLM2]{ 
        \label{3dPERR2LDLM2}     
        \includegraphics[scale=0.4]{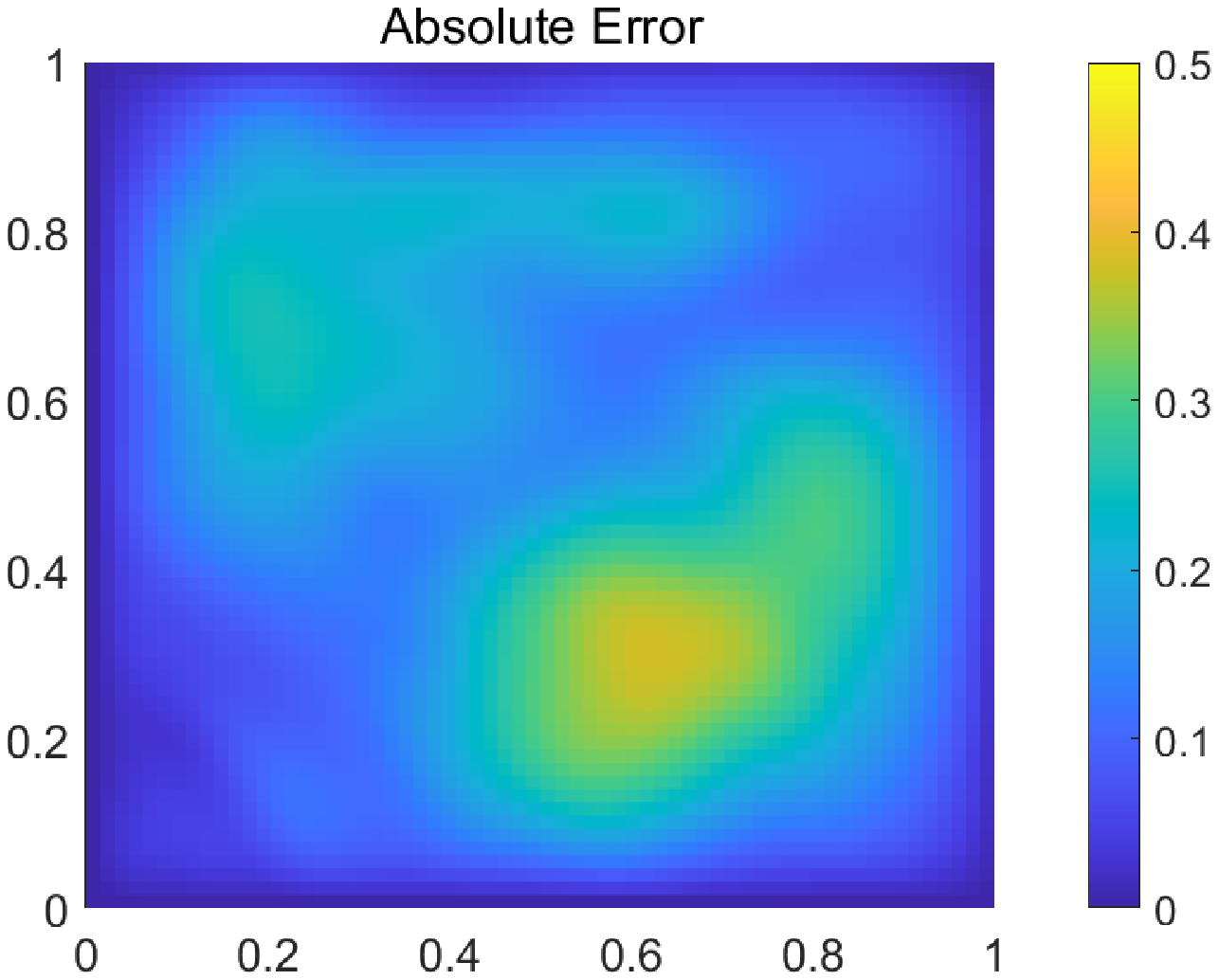}     
    } 
    \subfigure[REL] { 
        \label{3dPDE1q6:e}     
        \includegraphics[scale=0.38]{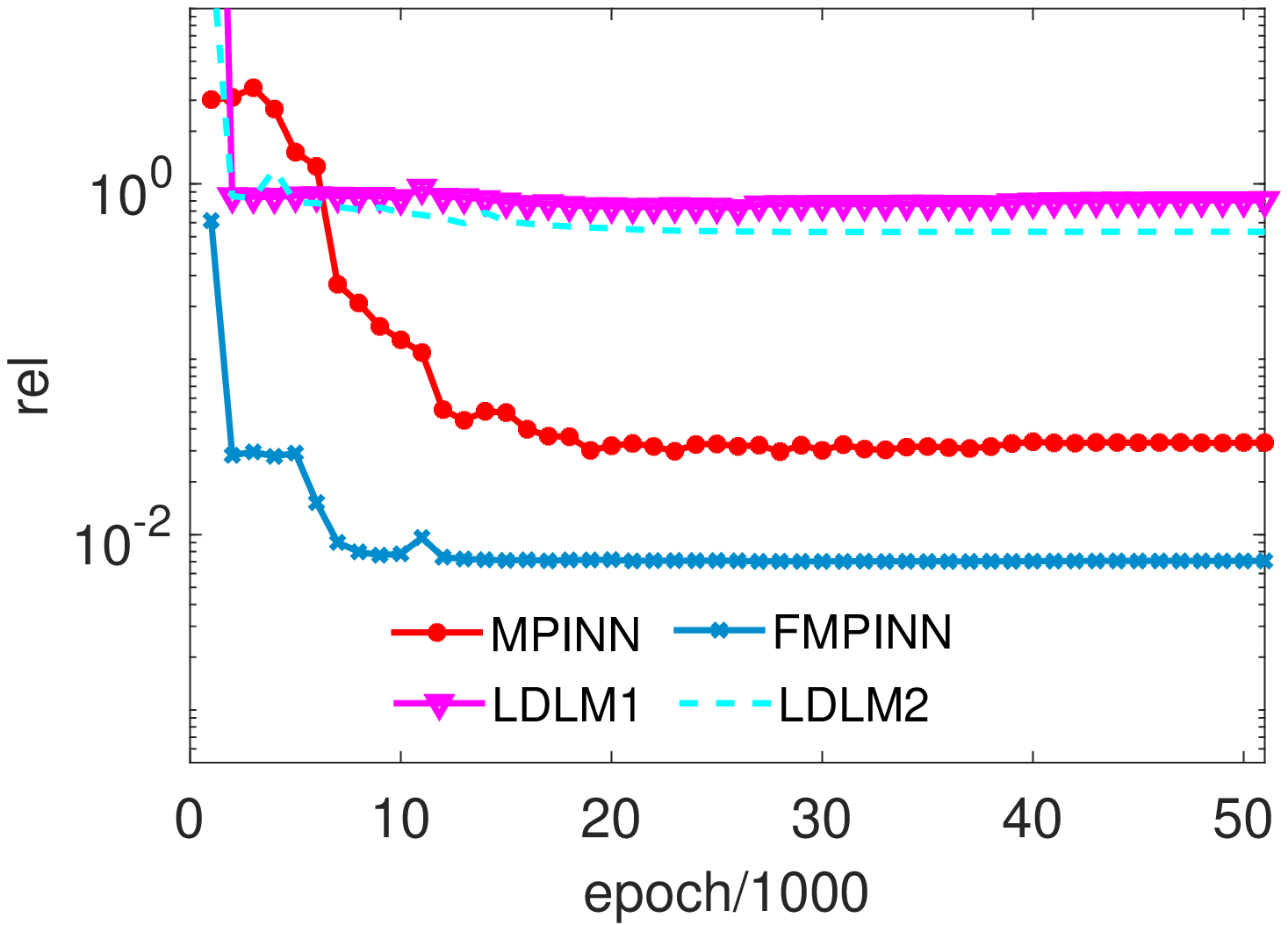}     
    }
    \caption{Rough coefficient, exact solution, loss of flux term and testing results for Example \ref{E3d}} 
    \label{3dPDE}         
\end{figure}

We ultilize the FMPINN, MPINN, LDLM1 and LDLM2 models to approximate the solution of three-dimensional multi-scale problem \eqref{eq:multiscale} with rough coefficient \eqref{coeff_3D} when $\varepsilon=0.1$, the setups the four models are  same as the Example \ref{Example2d_01}. The training dataset includes 7500 interior points and 1000 boundary points randomly sampled from $\Omega$ and $\partial \Omega$, respectively. To facilitate the process, a reference solution $u^{\varepsilon}(x_1,x_2, x_3)$ is established as the numerical solution obtained using the finite difference method on the domain $[-1,1]\times[-1,1]\times[-1,1]$ with a mesh-size $h=1/64$. The test dataset is formed by including all grid points within the domain $[-1,1]\times[-1,1]$ with a mesh-size $h=1/64$, while keeping the value of $z$ fixed at 0.3125. We list the total running time and REL in Table \ref{results3D} and plot the related results in Fig. \ref{3dPDE}. 

\begin{table}[H]
	\centering
	\caption{The relative error and running time of FMPINN, MPINN, LDLM1 and LDLM2 for Example \ref{E3d}.}
	\label{results3D}
	\begin{tabular}{|c|c|c|c|c|}
		\hline
		Method         &FMPINN    &MPINN    &LDLM1    &LDLM2    \\  \hline
		REL            &0.0071   &0.0335   &0.8048   &0.5326    \\  \hline
		Total time(s)  &5179.601  &9271.072&1065.541  &1195.233 \\  \hline
	\end{tabular}
\end{table}

Based on the results in Fig.\ref{3dPDE}, we can see that our FMPINN model still outperforms the MPINN and LDLMs model for multi-scale problems in three-dimensional space. The point-wise absolute error and the relative error of the former one are much smaller than that of the latter three, the precision of FMPINN is very good with a small absolute point-wise error. Additionally, the REL curve and the loss curve of the flux term are all flat in the later period of the training process, which means the performance of FMPINN is stable. The running time of FMPINN is 5179.601 seconds and less 3800 seconds than MPINN's.

\begin{example}\label{E8d}
	We consider the following eight-dimensional problem for \eqref{eq:multiscale} with Dirichlet boundary in regular  domains $\Omega=[0,1]^8$. In which, we take
	\begin{equation*}
     \begin{aligned}
         A(x_1,x_2,\cdots,x_{8})=1+\frac{1}{8}\bigg{[}&\cos(2\pi x_1)+\cos(4\pi x_2)+\cos(8\pi x_3)+\cos(16\pi x_4)+\\
         &\cos(16\pi x_5)+\cos(8\pi x_6)+\cos(4\pi x_7)+\cos(2\pi x_8)\bigg{]}.
     \end{aligned}
	\end{equation*}
	Meantime, an exact solution satisfied \eqref{eq:multiscale}  is given by
 \begin{equation*}
     u(x_1,x_2,\cdots,x_{8})=\prod_{j=1}^{8}\sin(\pi x_j)
	\end{equation*}
 The functions $f(x_1,x_2,\cdots,x_{8})$ in $\Omega$ and $g(x_1,x_2,\cdots,x_{8})$ on $\partial \Omega$ are easy to obtain according to the rough coefficient and exact solution, we omit it.
\end{example}

\begin{figure} [H]
	\centering    
	\subfigure[loss of flux term for FMPINN]{
		\label{8dPDE1:c}     
		\includegraphics[scale=0.33]{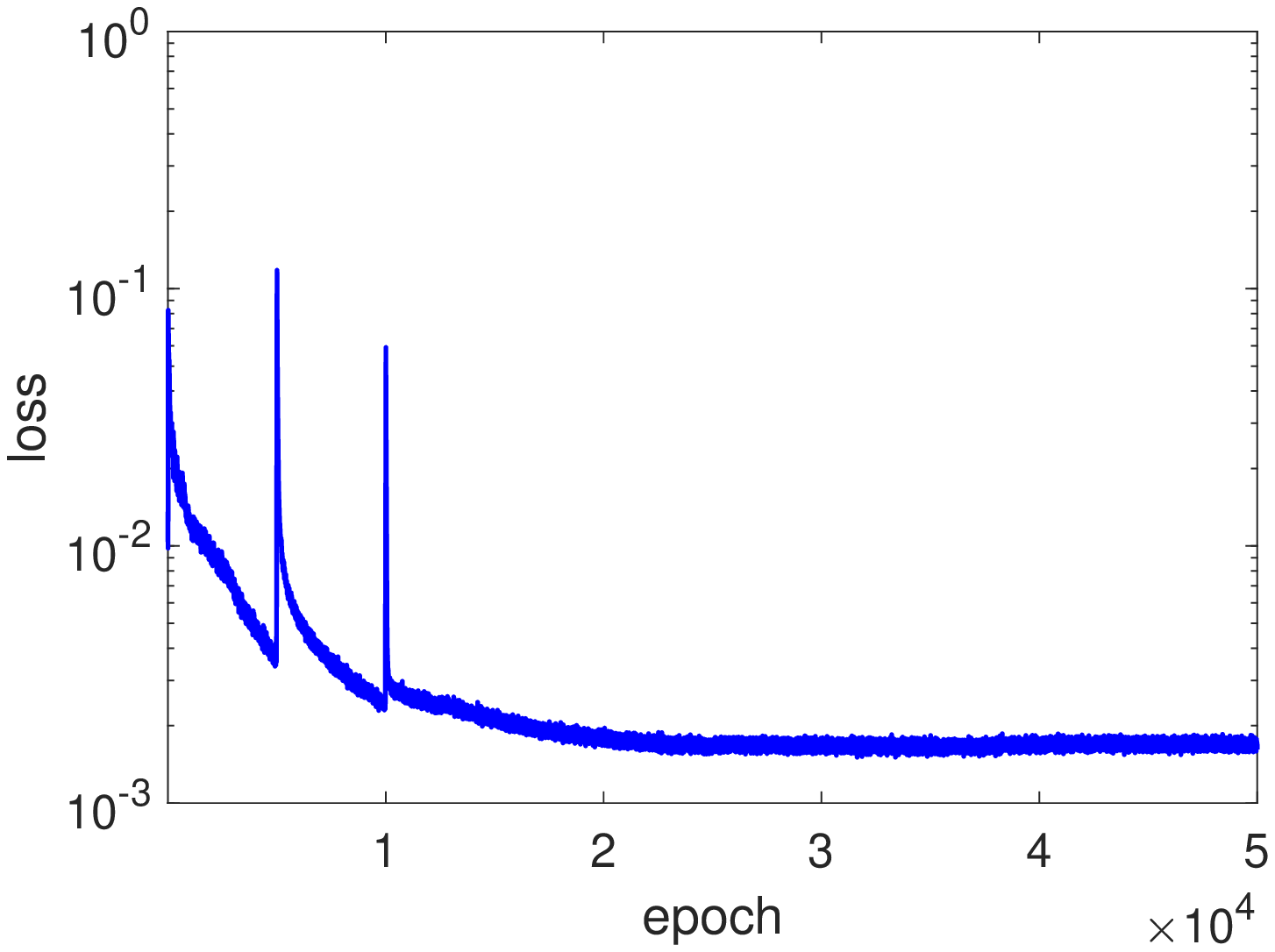}  
	}  
	\subfigure[point-wise error for FMPINN]{
		\label{8dPDE1_PERR2FMPINN}     
		\includegraphics[scale=0.33]{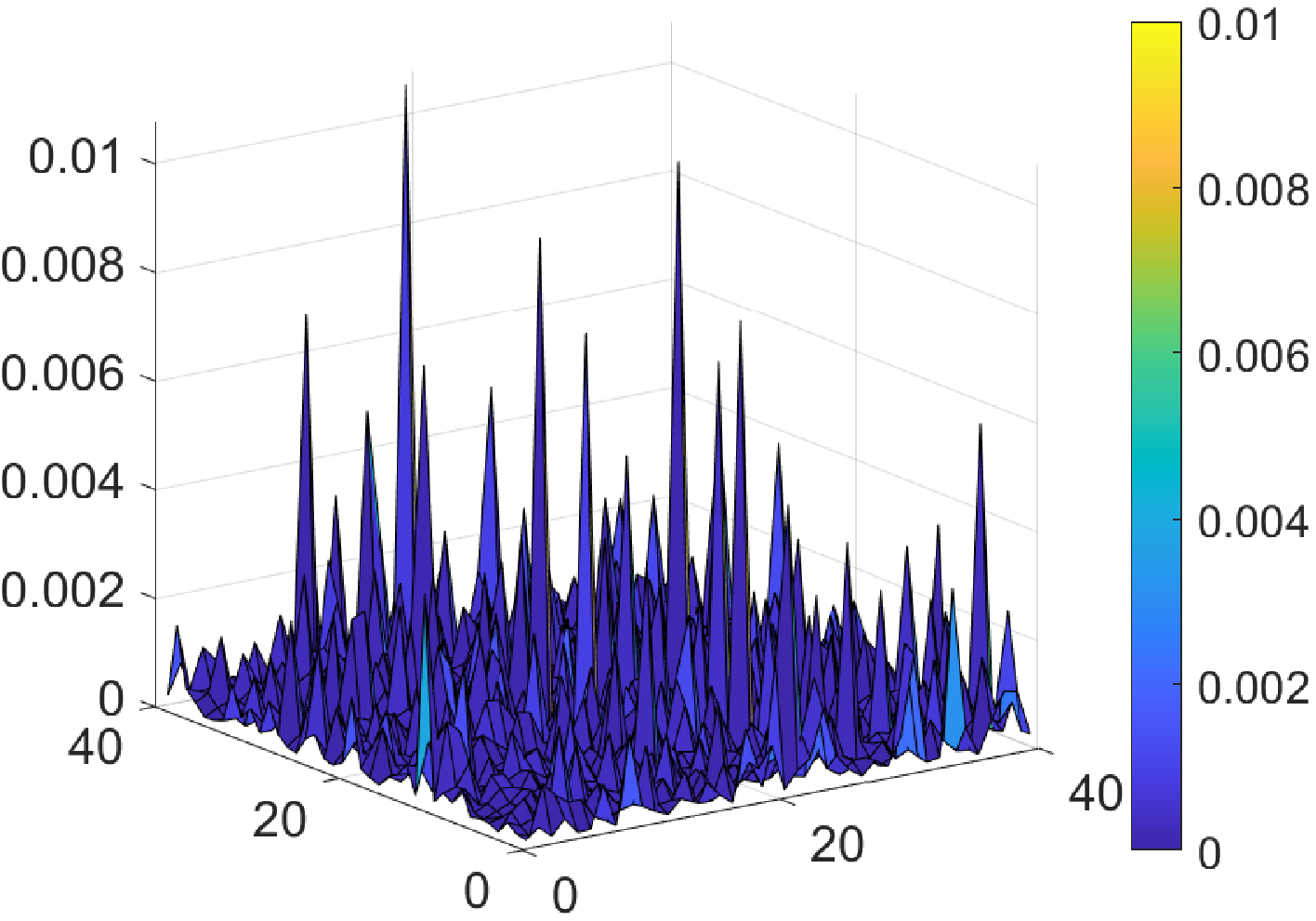}  
	}    
 	\subfigure[point-wise error for LDLM1]{
		\label{8dPDE1_PERR2LDLM_REQU}     
		\includegraphics[scale=0.33]{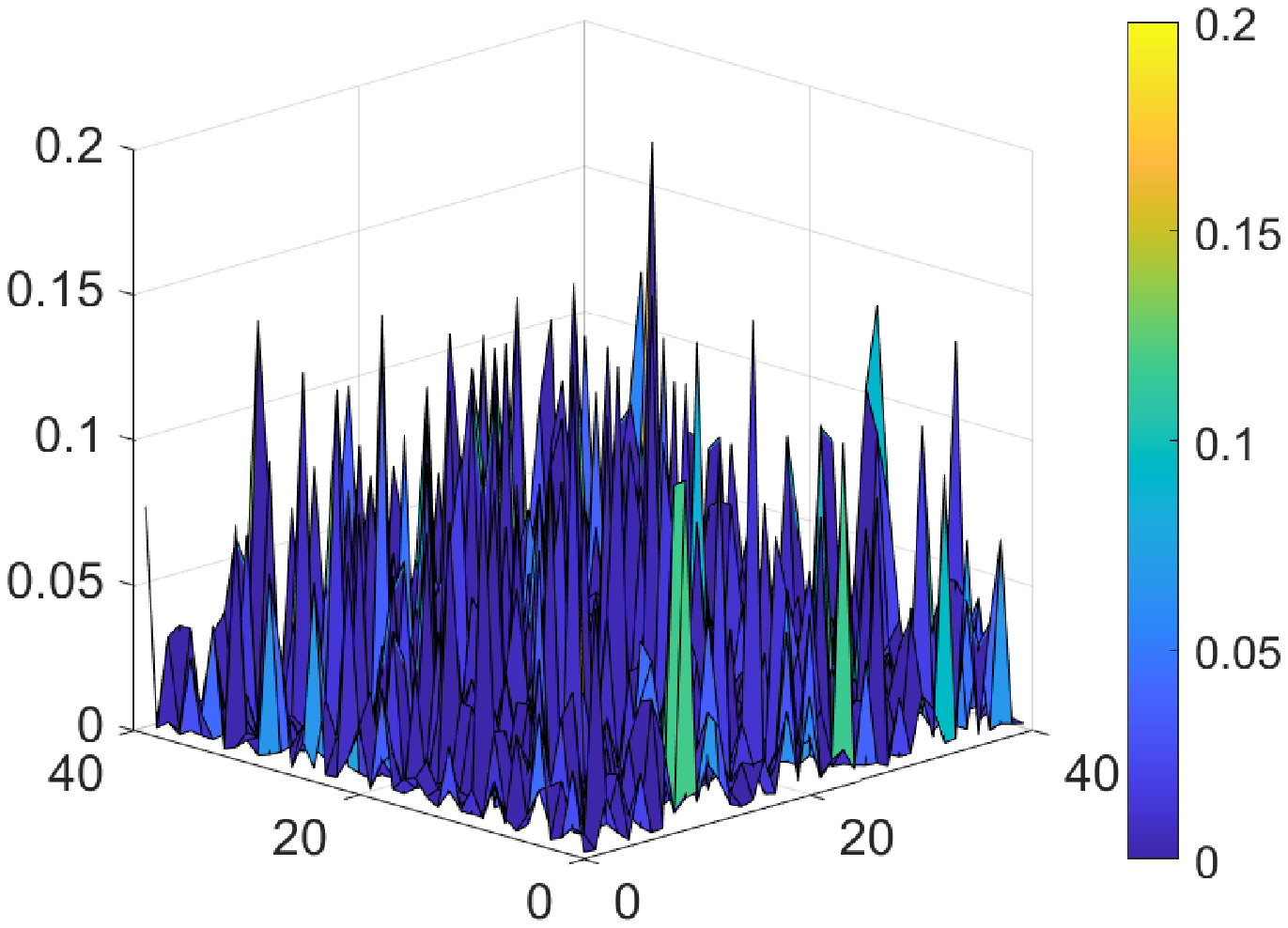}  
	}    
 	\subfigure[point-wise error for LDLM2]{
		\label{8dPDE1_PERR2LDLM_SINCOS}     
		\includegraphics[scale=0.33]{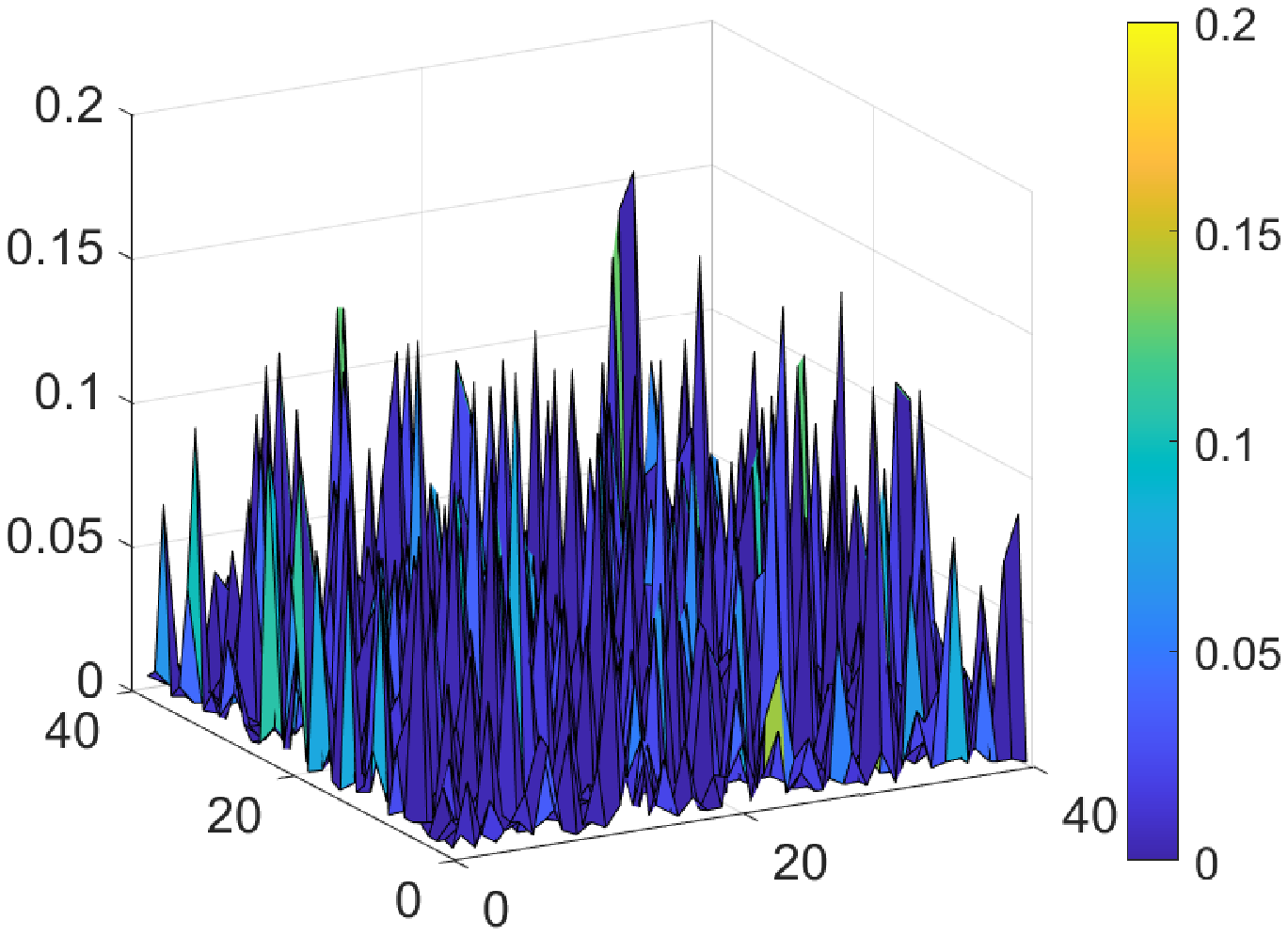}  
	}    
	\subfigure[REL]{ 
		\label{8dPDE1:f}     
		\includegraphics[scale=0.33]{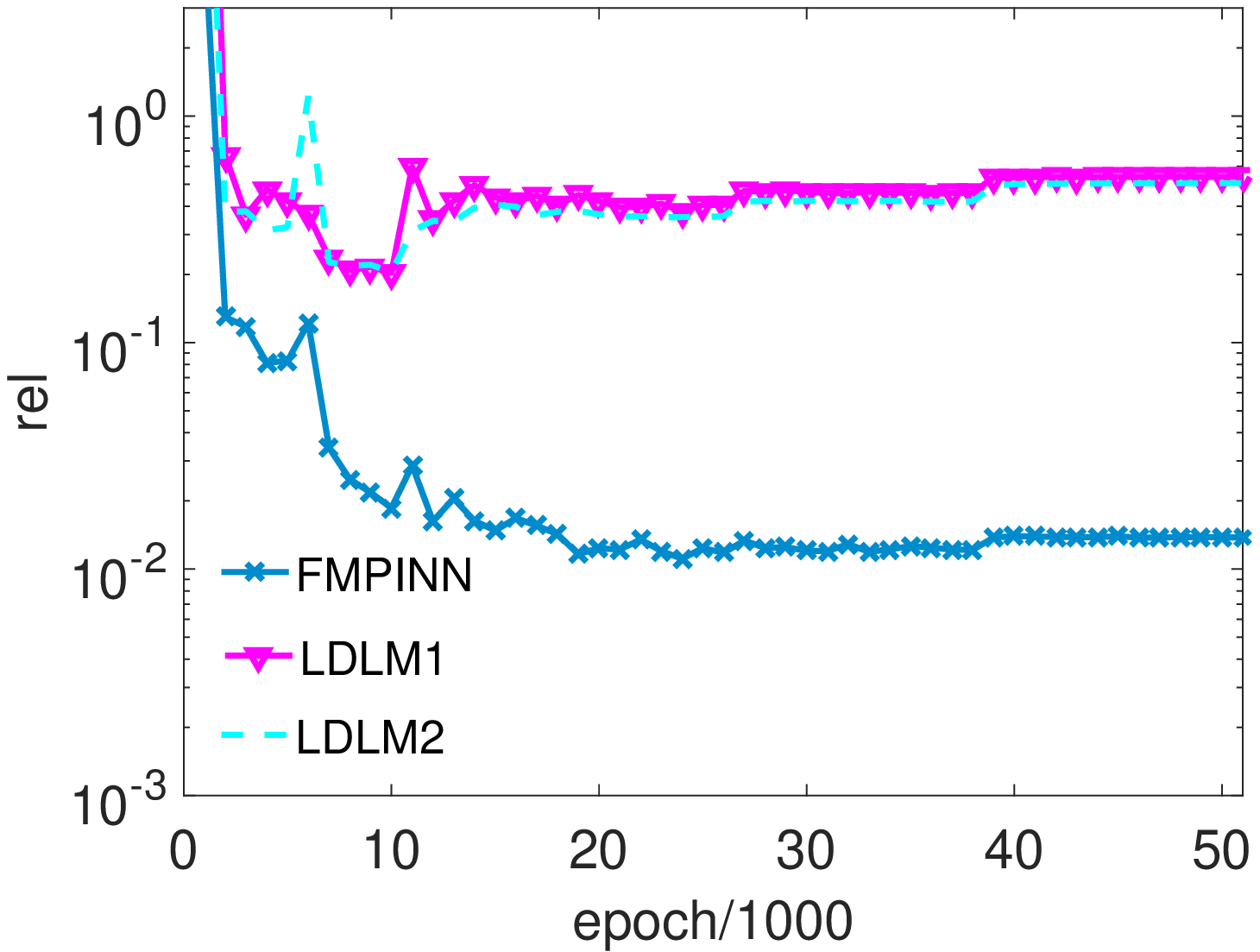}     
	}
	\caption{Loss of flux term and testing results for Example \ref{E8d}.}
	\label{8DResults}         
\end{figure}

In this example, we only perform the FMPINN, LDLM1 and LDLM2 model to solve the \eqref{eq:multiscale} in eight-dimensional space, because the huge computation requirement of MPINN has exceeded the limitation of memory for our station. The size of hidden layers for each subnetwork of FMPINN is set as (60, 80, 60, 60, 60) and the hidden layers' size for LDLM is set as $(400, 500, 300, 300, 300)$. At each training step, we construct the training dataset by sampling 20000 interior points inside the $\Omega$ and 5000 boundary points from the $\partial\Omega$. A testing dataset is given that included 1600 random points distributed in $\Omega$. The related experiment results are plotted in Fig.\ref{8DResults} and listed in Table \ref{results8D}. Additionally, the point-wise error for the FMPINN model evaluated on 1600 sample points is projected into a rectangular region with mesh size $40\times 40$. Noting that the mapping is only aimed at visualizing, it is independent of the actual coordinates of those points.

\begin{table}[H]
	\centering
	\caption{ The relative error and running time of FMPINN, LDLM1 and LDLM2  for Example \ref{E8d}.}
	\label{results8D}
	\begin{tabular}{|c|c|c|c|}
		\hline
		Method         &FMPINN       &LDLM1    &LDLM2     \\  \hline
		REL            &0.01378      &0.5394   &0.5054    \\  \hline
		Total time(s)  &21035.011    &1757.92  &1942.742  \\  \hline
	\end{tabular}
\end{table}

For an eight-dimensional problem, the FMPINN still can obtain a satisfactory solution for \eqref{eq:multiscale} with small point-wise absolute error and relative error. However, the LDLM1 and LDLM2 both fail to approximate the solution of \eqref{eq:multiscale}. Additionally, the loss of flux term and overall REL show that the FMPINN model is also stable during the training process. The running time of LDLMs is less thah that of FMPINN in Table \ref{results8D}, but their performance are obviously weaker that the latter's.

\section{Conclusion}\label{sec:05}
Physics-informed neural networks (PINN) have gained significant popularity in solving both forward and inverse problems. However, the normal PINN with a multi-scale DNN framework is unable to solve multiscale PDEs with rough coefficients. Inspired by the mixed finite element method, this work designs a Fourier-based mixed PINN(dubbed FMPINN) by combining a dual (flux) technique and Fourier decomposition to solve a class of elliptic multi-scale PDEs. By incorporating the loss of the flux term into the loss function, our model achieves improved stability and robustness. To handle multi-frequency contents, a Fourier activation function has been used to address the input data transformed radially by different frequency factors, and a sub-network is designed to match the target function, this strategy can improve clearly the accuracy and convergence rate for the FMPINN method. Compared to the previous works of PINN, this novel method skillfully casts the original problem into two first-order systems, it will overcome the shortcomings of the computational burden for high-order derivatives in DNN and the ill-condition of neural tangent kernel matrix resulting from the rough coefficient. Computational results show this novel method is feasible and efficient to solve this multi-scale equation with an inhomogeneous coefficient in various dimensional spaces. In the future, we aim to extend this novel network architecture, incorporating Fourier theory and lower-order mixed schemes, to tackle more complex multiscale problems.

\section*{Declaration of interests}
The authors declare that they have no known competing financial interests or personal relationships that could have appeared to influence the work reported in this paper.

\section*{Credit authorship contribution Statement}
Xi'an Li: Conceptualization, Methodology, Investigation, Validation, Writing - Original Draft. 
Jinran Wu: Investigation, Writing - Review \& Editing.
You-Gan Wang: Writing - Review \& Editing, Xin Tai: Writing - Review \& Editing, Jianhua Xu: Writing - Review \& Editing.

\section*{Acknowledgements}
The authors wish to thank Prof. Dr. Zhi-Qin John Xu for valuable suggestions which improved the quality of the paper. 
 
\bibliography{References}
\bibliographystyle{model1-num-names}
\end{document}